\documentclass[12pt]{amsart}
\usepackage{amssymb}
\usepackage{bbm}

\usepackage{tikz}
\usetikzlibrary{decorations.pathreplacing}
\usetikzlibrary{shapes.misc}
\usepackage{tikz-cd}

\usepackage[pdfencoding=auto,psdextra]{hyperref}
\usepackage{textgreek}
\usepackage{enumerate}

\usepackage[hmargin=1in,vmargin=1.25in]{geometry}
   \title{Dens, nests and the Loehr-Warrington conjecture}

   \author{J. Blasiak}
   \author{M. Haiman}
   \author{J. Morse}
   \author{A. Pun}
   \author{G. H. Seelinger}

   \address[Blasiak]
   {Dept.\ of Mathematics\\
    Drexel University\\
    Philadelphia, PA}
   \email{jblasiak@gmail.com}

   \address[Haiman]
   {Dept.\ of Mathematics\\
    University of California\\
    Berkeley, CA}
   \email{mhaiman@math.berkeley.edu}

   \address[Morse]
   {Dept.\ of Mathematics\\
    University of Virginia\\
    Charlottesville, VA}
   \email{morsej@virginia.edu}

   \address[Pun]
    {Dept.\ of Mathematics\\
     University of Virginia\\
     Charlottesville, VA}
   \email{ayp6e@virginia.edu}

   \address[Seelinger]
   {Dept.\ of Mathematics\\
    University of Michigan\\
    Ann Arbor, MI}
   \email{ghseeli@umich.edu}

   \thanks{Authors were supported by NSF Grants DMS-1855784 (J.~B.)
    and DMS-1855804 (J.~M.).}

   \date{\today}
\newtheorem{thm}{Theorem}[subsection]
\newtheorem{lemma}[thm]{Lemma}
\newtheorem{prop}[thm]{Proposition}
\newtheorem{cor}[thm]{Corollary}

\theoremstyle{definition}
\newtheorem{defn}[thm]{Definition}
\newtheorem{convention}[thm]{Convention}

\theoremstyle{remark}
\newtheorem{example}[thm]{Example}

\newtheorem{remark}[thm]{Remark}

\newcommand{\NN}{{\mathbb N}}
\newcommand{\QQ}{{\mathbb Q}}

\newcommand{\ZZ}{{\mathbb Z}}

\newcommand{\kk}{{\mathbbm k}}

\newcommand{\bb}{{\mathbf b}}
\newcommand{\dd}{{\mathbf d}}
\newcommand{\ee}{{\mathbf e}}
\newcommand{\gG}{{\mathbf g}}
\newcommand{\rr}{{\mathbf r}}
\newcommand{\sS}{{\mathbf s}}
\newcommand{\xx}{{\mathbf x}}
\newcommand{\yy}{{\mathbf y}}
\newcommand{\zz}{{\mathbf z}}

\newcommand{\Hbold}{{\mathbf H}}
\newcommand{\Sbold}{{\mathbf S}}
\newcommand{\sigmabold}{{\boldsymbol \sigma }}
\newcommand{\alphabold}{{\boldsymbol \alpha }}
\newcommand{\betabold}{{\boldsymbol \beta }}
\newcommand{\gammabold}{{\boldsymbol \gamma }}
\newcommand{\nubold}{{\boldsymbol \nu }}
\newcommand{\onebold}{{\boldsymbol 1}}
\newcommand{\deltabold}{{\boldsymbol \delta }}
\newcommand{\epsbold}{{\boldsymbol \varepsilon }}

\newcommand{\Acal}{{\mathcal A}}
\newcommand{\Ecal}{{\mathcal E}}
\newcommand{\Gcal}{{\mathcal G}}
\newcommand{\Hcal}{{\mathcal H}}
\newcommand{\Lcal}{{\mathcal L}}
\newcommand{\Scal}{{\mathcal S}}

\newcommand{\ctild}{\tilde{c}}
\newcommand{\Htild}{\tilde{H}}
\newcommand{\chat}{\hat{c}}
\newcommand{\sigmahat}{\widehat{\sigma }}
\newcommand{\jhat}{\widehat{\jmath}}

\DeclareFontFamily{U}{mathx}{\hyphenchar\font45}
\DeclareFontShape{U}{mathx}{m}{n}{
 <5> <6> <7> <8> <9> <10> <10.95> <12> <14.4> <17.28> <20.74> <24.88>
 mathx10}{}
\DeclareSymbolFont{mathx}{U}{mathx}{m}{n}
\DeclareFontSubstitution{U}{mathx}{m}{n}
\DeclareMathAccent{\widecheck}{0}{mathx}{"71}
\newcommand{\hGamma}{{\widehat{\Gamma}}}

\newcommand{\semis}{{\boldsymbol \delta }}

\DeclareMathOperator{\area}{area}
\DeclareMathOperator{\dinv}{dinv}
\DeclareMathOperator{\inv}{inv}
\DeclareMathOperator{\modm}{mod}
\DeclareMathOperator{\pol}{pol}
\DeclareMathOperator{\adj}{adj}
\DeclareMathOperator{\spin}{spin}
\DeclareMathOperator{\sgn}{sgn}
\DeclareMathOperator{\sshare}{ss}
\DeclareMathOperator{\Inv}{Inv}
\DeclareMathOperator{\GL}{{\it GL}}

\DeclareMathOperator{\SSYT}{SSYT}
\DeclareMathOperator{\LNDP}{LNDP}

\newcommand{\defeq}{\mathbin{\overset{\text{{\rm def}}}{=}}}

\newcommand{\fp}[1]{{\boldsymbol \{} #1 {\boldsymbol \}}}

\tikzset{cross/.style={cross out, draw=black, minimum
size=2*(#1-\pgflinewidth), inner sep=0pt, outer sep=0pt},
cross/.default={1pt}}

\begin{document}


\begin{abstract}
In a companion paper, we introduced raising operator series called
{\em Catalanimals}.  Among them are {\em Schur Catalanimals}, which
represent Schur functions inside copies $\Lambda (X^{m,n})\subset
\Ecal $ of the algebra of symmetric functions embedded in the
elliptic Hall algebra $\Ecal $ of Burban and Schiffmann.

Here we obtain a combinatorial formula for symmetric functions given
by a class of Catalanimals that includes the Schur Catalanimals.  Our
formula is expressed as a weighted sum of LLT polynomials, with terms
indexed by configurations of nested lattice paths called {\em nests},
having endpoints and bounding constraints controlled by data called a
{\em den}.

Applied to Schur Catalanimals for the alphabets $X^{m,1}$ with $n=1$,
our `nests in a den' formula proves the combinatorial formula
conjectured by Loehr and Warrington for $\nabla^m s_{\mu }$ as a
weighted sum of LLT polynomials indexed by systems of nested Dyck
paths.  When $n$ is arbitrary, our formula establishes an $(m,n)$
version of the Loehr-Warrington conjecture.

In the case where each nest consists of a single lattice path, the
nests in a den formula reduces to our previous shuffle theorem for
paths under any line.  Both this and the $(m,n)$ Loehr-Warrington
formula generalize the $(km,kn)$ shuffle theorem proven by Carlsson
and Mellit (for $n=1$) and Mellit.  Our formula here unifies these two
generalizations.
\end{abstract}

\maketitle

\section{Introduction}
\label{s:intro}

\subsection{Background}
\label{ss:background}

In this paper we prove and extend the oldest unresolved conjecture in
`$q,t$-Catalan combinatorics,' namely, the combinatorial formula for
$\nabla^m s_{\mu }$ conjectured by Loehr and Warrington
\cite{LoehWarr08}, where $s_{\mu }$ is a Schur function and $\nabla$
is the operator from \cite{BeGaHaTe99} which is important in the
theory of Macdonald polynomials.  Like other results and conjectures
in this area, beginning with the {\em shuffle theorem} conjectured by
Haglund et.\ al.\ \cite{HaHaLoReUl05} and proven by Carlsson and
Mellit \cite{CarlMell18}, the Loehr-Warrington formula is expressed as
a sum over Dyck paths (in this case, systems of nested Dyck paths) of
LLT polynomials weighted by monomials in $q$ and $t$.

Our main result, Theorem~\ref{thm:nest-identity}, is considerably more
general than the Loehr-Warrington formula.  We briefly describe some
of its further consequences.

The simplest case of the Loehr-Warrington formula, when $s_{\mu } =
e_k$ is an elementary symmetric function, reduces to the original
shuffle theorem.  The latter is the $n=1$ case of an extended {\em
$(km,kn)$ shuffle theorem} conjectured by Bergeron et.\ al.\
\cite{BeGaSeXi16} and proven by Mellit \cite{Mellit16}.  A new
consequence of our main result here, not previously formulated even as
a conjecture, is a corresponding $(m,n)$ extension of the
Loehr-Warrington formula, which reduces to the $(km,kn)$ shuffle
theorem when $s_{\mu } = e_k$.

Another generalization of the $(km,kn)$ shuffle theorem is given by
our shuffle theorem for paths under any line \cite{BHMPS-paths}.  This
too is a consequence of our main theorem here.  Thus our theorem
unifies a number of previous results and conjectures, as summarized in
the following diagram.
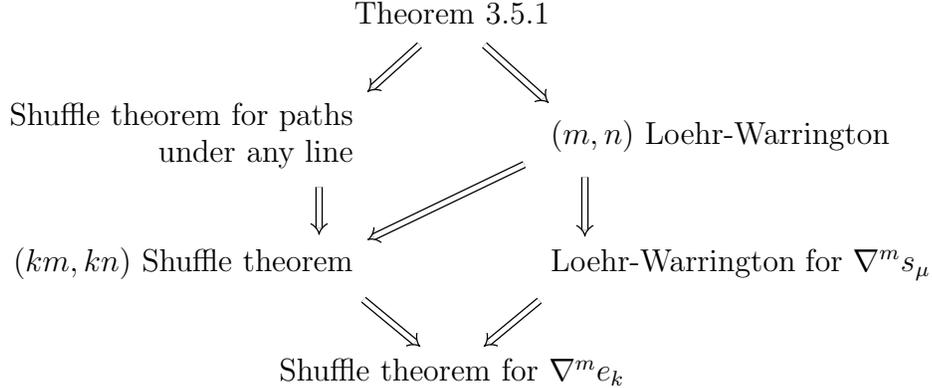
\[
\begin{tikzcd}[column sep=small]
&
\makebox[5ex]{Theorem~\ref{thm:nest-identity}}
\arrow[shorten <=1ex,ld,Rightarrow]
\arrow[shorten <=1ex,shorten >=.5ex,rd,Rightarrow]
\\
\makebox[5ex][r]{\parbox[c]{28ex}{\raggedleft Shuffle theorem for
paths under any line}}
\arrow[shorten <=.5ex,d,Rightarrow]
&&
\makebox[5ex][l]{$(m,n)$ Loehr-Warrington}
\arrow[shorten <=1ex,d,Rightarrow]
\arrow[shorten <=1ex,lld,Rightarrow]
\\
\makebox[5ex][r]{$(km,kn)$ Shuffle theorem}
\arrow[shorten <=1ex,rd,Rightarrow]
&&
\makebox[5ex][l]{Loehr-Warrington for  $\nabla ^{m} s_{\mu }$}
\arrow[shorten <=1ex,ld,Rightarrow]
\\
&
\makebox[5ex]{Shuffle theorem for $\nabla ^{m} e_{k}$}
\end{tikzcd}
\]

\subsection{Overview}
\label{ss:overview}

In \cite{BHMPS-llt}, we introduced raising operator series $H(\zz
;q,t)$ called {\em Catalanimals}; among them, we constructed examples
for which the polynomial truncation $H(\zz ;q,t)_{\pol }$ is equal, up
to an explicit factor of the form $\pm q^{r}t^{s}$, to $\omega\,
\nabla ^{m} s_{\mu }(\zz )$, where $\omega $ is the standard
involution on symmetric functions.

Our main result, Theorem~\ref{thm:nest-identity}, gives a
combinatorially defined expansion
\begin{equation}\label{e:nest-identity-preview}
H(\zz ;q,t)_{\pol } = \sum _{\pi } t^{a(\pi )} q^{\dinv _{p}(\pi )}
\Gcal _{\nubold (\pi )}(\zz ;q^{-1}),
\end{equation}
in terms of LLT polynomials $\Gcal _{\nubold }(\zz ;q)$, for a special
class of Catalanimals $H(\zz ;q,t)$ including those for which $H(\zz
;q,t)_{\pol } = \pm q^{r}t^{s}\, \omega \, \nabla ^{m} s_{\mu }(\zz
)$, as just discussed.

The terms on the right hand side of \eqref{e:nest-identity-preview}
are indexed by configurations of nested lattice paths $\pi =(\pi
_{1},\ldots,\pi _{r})$, called {\em nests}, with endpoints and
bounding constraints controlled by combinatorial data called a {\em
den}.  The statistics $a(\pi )$, $\dinv _{p}(\pi )$ generalize the
`area' and `dinv' statistics found in the shuffle theorem and its
friends.  We define these combinatorial notions in \S
\ref{s:dens-and-nests}.

In the case where the left hand side of
\eqref{e:nest-identity-preview} becomes $ \pm q^{r}t^{s}\, \omega \,
\nabla ^{m} s_{\mu }(\zz )$, formula \eqref{e:nest-identity-preview}
proves the Loehr-Warrington conjecture (see Theorem \ref{thm:mn-LW}
and \S \ref{ss:original-LW}).

Formula \eqref{e:nest-identity-preview} also applies to more general
{\em Schur Catalanimals} $H(\zz ;q,t) = H_{(\mu ^{\circ
})^{m}}^{m,n}$, which were shown in \cite{BHMPS-llt} to represent
(again up to a factor $\pm q^{r}t^{s}$) Schur functions $s_{\mu
}[-MX^{m,n}]$ in subalgebras $\Lambda (X^{m,n})\subset \Ecal $
isomorphic to the algebra of symmetric functions $\Lambda $ over $\QQ
(q,t)$, where $\Ecal $ is the elliptic Hall algebra of Burban and
Schiffmann \cite{BurbSchi12}---see \S \ref{s:Catalanimals} for
details. Under the action of $\Ecal $ on $\Lambda $ constructed by
Schiffmann and Vasserot \cite{SchiVass13}, the Schur Catalanimal $H =
H_{(\mu ^{\circ })^{m}}^{m,n}$ satisfies $H_{\pol } = \pm q^{r}t^{s}
\omega ( s_{\mu }[-M X^{m,n}]\cdot 1)(\zz )$.  For $n=1$, we have
$s_{\mu }[-M X^{m,1}]\cdot 1 = \nabla ^{m} s_{\mu }$.  In this case,
the Schur Catalanimals are the Catalanimals referred to above.

For general $n$, \eqref{e:nest-identity-preview} yields a
combinatorial formula for $s_{\mu }[-M X^{m,n}]\cdot 1$, made precise
in Theorem~\ref{thm:mn-LW}, which can be naturally understood as an
$(m,n)$ extension of the Loehr-Warrington conjecture.  For $\mu
=(1^{k})$, $s_{\mu } = e_{k}$ is the $k$-th elementary symmetric
function, and our $(m,n)$ Loehr-Warrington theorem reduces to the
$(km,kn)$ shuffle theorem of \cite{BeGaSeXi16,Mellit16}, just as the
original Loehr-Warrington conjecture for $\nabla ^{m} s_{\mu }$
reduces to the classical shuffle theorem
\cite{CarlMell18,HaHaLoReUl05} for $\nabla ^{m} e_{k}$.

Finally, for dens such that each nest consists of a single lattice
path, \eqref{e:nest-identity-preview} reduces to our shuffle theorem
for paths under any line \cite[Theorem 5.5.1]{BHMPS-paths}, which also
generalizes the $(km,kn)$ shuffle theorem.  Thus, we have the diagram
of implications in \S \ref{ss:background}, above, with formula
\eqref{e:nest-identity-preview} at the top.

As with other instances of $q,t$-Catalan combinatorics, the left hand
side of \eqref{e:nest-identity-preview} is symmetric in $q$ and $t$ by
construction; hence the right hand side shares this symmetry.  No
purely combinatorial explanation of this symmetry is yet known, even
in the case of the classical shuffle theorem for $\nabla e_{k}$.

In addition, since LLT polynomials are $q$-Schur positive
\cite{GrojHaim07}---i.e., their coefficients in terms of Schur
functions belong to $\NN [q]$---it follows that the Catalanimals to
which \eqref{e:nest-identity-preview} applies are $(q,t)$-Schur
positive.  The general question of which Catalanimals are
$(q,t)$-Schur positive seems to be a difficult one.  See
\cite[Conjecture~7.1.1]{BHMPS-paths} for one conjecture in this
direction.

\subsection{Method and outline}
\label{ss:method}

We prove our main theorem by a method parallel to the one we used to
prove the shuffle theorem for paths under any line in~\cite{BHMPS-paths}.

We obtain the combinatorial formula in Theorem~\ref{thm:nest-identity}
by taking the polynomial part of an identity between infinite series
of $\GL _{l}$ characters.  The latter identity, equation
\eqref{e:stable-main} in Theorem~\ref{thm:stable-main}, expands the
full Catalanimal $H(\zz ;q,t)$ associated with a den as an infinite
sum of {\em LLT series} $\Lcal ^{\sigma }_{\rr ,\beta /\alpha }(\zz
;q)$ weighted by powers of $t$.

Upon taking the polynomial part, all but a finite number of the terms
$t^{a}\Lcal ^{\sigma }_{\rr ,\beta /\alpha }(\zz ;q)$ in
\eqref{e:stable-main} vanish.  The surviving terms are indexed by
nests $\pi $ in the given den, and have polynomial parts $t^{a}\Lcal
^{\sigma }_{\rr ,\beta /\alpha }(\zz ;q)_{\pol } = t^{a(\pi )}
q^{\dinv _{p}(\pi )} \Gcal _{\nubold (\pi )}(\zz ;q^{-1})$, yielding
\eqref{e:nest-identity-preview}.

Given a Levi subgroup $\GL _{\rr } = \GL _{r_{1}}\times \cdots \times
\GL _{r_{k}}$ of $\GL _{l}$, the LLT series $\Lcal ^{\sigma }_{\rr
,\beta /\alpha }(\zz ;q)$ in $l$ variables $\zz =z_{1},\ldots,z_{l}$
(Definition~\ref{def:LLT-series}) encapsulates the matrix coefficients
of multiplication by arbitrary $\GL _{l}$ characters with respect to
chosen basis elements $E^{\sigma }_{\rr ,\alpha }(\zz ;q)$, $E^{\sigma
}_{\rr ,\beta }(\zz ;q)$ of the space of virtual $\GL _{\rr }$
characters.  Here $E^{\sigma }_{\rr ,\lambda }(\zz ;q)$ denotes a
(twisted) {\em semi-symmetric Hall-Littlewood polynomial}
(Definition~\ref{def:E-F-semi}).

The orthogonality of semi-symmetric Hall-Littlewood polynomials
(Proposition~\ref{prop:orthogonality}) leads to a formula for LLT
series in terms of these polynomials
(Proposition~\ref{prop:L-via-Hq}).  Using this formula, the desired
infinite series identity \eqref{e:stable-main} follows from a Cauchy
identity for semi-symmetric Hall-Littlewood polynomials,
Theorem~\ref{thm:Cauchy}, along with an auxiliary identity,
Proposition~\ref{prop:winding}, that relates semi-symmetric
Hall-Littlewood polynomials with different twists.

The steps just outlined parallel those in the proof of
\cite[Theorem~5.5.1]{BHMPS-paths}, although many of the details are more
intricate.
Readers may find the simpler argument in
\cite{BHMPS-paths}, which covers the case  $\rr = (1^{l})$, a helpful
guide to the argument here.  

Chief among the new intricacies is that the Cauchy identity for
semi-symmetric Hall-Littlewood polynomials in Theorem~\ref{thm:Cauchy}
is more subtle than the one for non-symmetric Hall-Littlewood
polynomials in \cite[Theorem~5.1.1]{BHMPS-paths}.  The new Cauchy
identity involves semi-symmetric Hall-Littlewood polynomials for two
separate Levi subgroups $\GL _{\rr }$ and $\GL _{\sS }$, along with
the choice of minimal dominant regular weights $\rho _{\rr }, \rho
_{\sS }$ for each of them.  These choices must satisfy certain
compatibilities in order for the Cauchy identity to hold.  Because of
this greater complexity, we are not able to give a short proof of
Theorem~\ref{thm:Cauchy} like we did for
\cite[Theorem~5.1.1]{BHMPS-paths}.  Instead, we devote most of \S
\ref{ss:Cauchy} to developing a series of properties of semi-symmetric
Hall-Littlewood polynomials, which we then use to prove
Theorem~\ref{thm:Cauchy}.

\section{Catalanimals and LLT polynomials}
\label{s:Catalanimals}

\subsection{Symmetric function conventions}
\label{ss:symmetric-functions}

The (French style) diagram of a partition $\lambda$ is the set of
lattice points $\{(i,j)\mid 1\leq j\leq \ell(\lambda ),\; 1\leq i \leq
\lambda _{j}\}$, where $\ell(\lambda )$ is the length of $\lambda $.
We often identify $\lambda$ and its diagram with the set of lattice
squares, or {\em boxes}, with northeast corner at a point $(i,j)\in
\lambda$.  A {\it skew diagram} is a difference $\nu = \lambda / \mu$
of partition diagrams $\mu\subseteq\lambda$, or any translate of such
a diagram.  This allows for skew diagrams $\nu = \beta /\alpha $ in
which the $x$-coordinates $\alpha _{i}, \beta _{i}$ of the left and
right ends of the rows may be negative.

The {\em content} of a box $a = (i,j)$ in row $j$, column $i$ of a
(skew) diagram is $c(a) = i-j$.

Let $\Lambda = \Lambda (X)$ be the algebra of symmetric functions in
infinitely many variables $X = x_{1},x_{2},\ldots$, with coefficients
in the field $\kk = \QQ (q,t)$.  We follow Macdonald's notation
\cite{Macdonald95} for the graded bases of $\Lambda $, the Hall inner
product $\langle -,- \rangle$ in which the Schur functions $s_{\lambda
}$ are orthonormal, and the automorphism $\omega \colon \Lambda
\rightarrow \Lambda $ such that $ \omega s_{\lambda } =
s_{\lambda^*}$, where $\lambda^*$ denotes the transpose of a partition
$\lambda$.

Given $f\in \Lambda $ and any expression $A$ involving indeterminates,
such as a polynomial, rational function, or formal series, the
plethystic evaluation $f[A]$ is defined by writing $f$ as a polynomial
in the power-sums $p_{k}$ and evaluating with $p_{k}\mapsto p_{k}[A]$,
where $p_{k}[A]$ is the result of substituting $a^{k}$ for every
indeterminate $a$ occurring in $A$.  The variables $q, t$ from our
ground field $\kk $ count as indeterminates.

By convention, the name of an alphabet $X = x_{1},x_{2},\ldots$ stands
for $x_{1}+x_{2}+\cdots $ inside a plethystic evaluation.  Then $f[X]
= f[x_{1}+x_{2}+\cdots ] = f(x_{1},x_{2},\ldots) = f(X)$.  A special
case of this convention that will arise often is the following.  We
fix
\begin{equation}\label{e:M}
M = (1-q)(1-t)
\end{equation}
here and throughout.  Then the evaluation $f[-MX]$ is the image of
$f(X)$ under the $\kk $-algebra automorphism of $\Lambda $ that sends
$p_{k}$ to $-(1-q^{k})(1-t^{k})p_{k}$.

We also allow plethystic evaluation term by term in a symmetric formal
series, provided the result makes sense formally.  In particular, the
series
\begin{equation}\label{e:Omega}
\Omega = \sum _{n=0}^{\infty } h_{n} = \exp \sum _{k=1}^{\infty }
\frac{p_{k}}{k}
\end{equation}
(with $h_{0} = 1$) has the property
\begin{equation}\label{e:Omega-pleth}
\Omega [x_{1}+x_{2}+\cdots -y_{1}-y_{2}-\cdots ] = \frac{\prod
_{i}(1 - y_{i})}{\prod _{i}(1 - x_{i})}.
\end{equation}

The linear operator $\nabla $ on $\Lambda $, introduced in
\cite{BeGaHaTe99}, is defined to act diagonally in the basis of
modified Macdonald polynomials $\Htild_{\mu
}(X;q,t)$~\cite{GarsHaim96}, with
\begin{equation}\label{e:nabla}
\nabla \Htild _{\mu } = t^{n(\mu )}q^{n(\mu ^{*})} \Htild _{\mu },
\end{equation}
where $n(\mu ) = \sum _{i} (i-1)\mu _{i}$.

\subsection{LLT polynomials}
\label{ss:LLT}

We recall the definition and basic properties of LLT polynomials
\cite{LaLeTh97}, using the `attacking inversions' formulation from
\cite{HaHaLoReUl05}.

Let $\nubold = (\nu_{(1)},\ldots,\nu_{(k)})$ be a tuple of skew
diagrams.  We consider the set of boxes in $\nubold $ to be the
disjoint union of the sets of boxes in the $\nu _{(i)}$, and define
the {\em adjusted content} of a box $a\in \nu_{(i)}$ to be $\ctild (a)
= c(a)+i\, \epsilon $, where $\epsilon $ is a fixed positive number
such that $k\, \epsilon <1$.

A {\em diagonal} in $\nubold $ is the set of boxes of a fixed adjusted
content, that is, a diagonal of fixed content in one of the $\nu
_{(i)}$.

The {\em reading order} on $\nubold $ is the total ordering $<$ on the
boxes of $\nubold $ such that $a<b \Rightarrow \ctild (a)\leq \ctild
(b)$ and boxes on each diagonal increase from southwest to northeast.
An {\em attacking pair} is an ordered pair of boxes $(a,b)$ in $
\nubold$ such that $a<b$ in reading order and $0<\ctild
(b)-\ctild(a)<1$.

A {\em semistandard tableau} on the tuple $\nubold $ is a map $T\colon
\nubold \rightarrow \ZZ _{+}$ which restricts to a semistandard Young
tableau on each component $\nu_{(i)}$.  The set of these is denoted
$\SSYT (\nubold )$.  An {\em attacking inversion} in $T$ is an
attacking pair $(a,b)$ such that $T(a)>T(b)$.  The number of attacking
inversions in $T$ is denoted $\inv (T)$.

\begin{defn}\label{def:G-nu}
The {\em LLT polynomial} indexed by a tuple of skew diagrams $\nubold
$ is the generating function, which is known to be symmetric
\cite{HaHaLoReUl05,LaLeTh97},
\begin{equation}\label{e:G-nu}
\Gcal_{\nubold  }(X;q) = \sum _{T\in \SSYT (\nubold 
)}q^{\inv (T)}\xx ^{T},
\end{equation}
where $\xx ^{T} = \prod _{a\in \nubold  } x_{T(a)}$.
\end{defn}

A similar formula expresses $\omega \, \Gcal _{\nubold }(X;q)$ as a
generating function for tableaux, as follows.  Fix an ordered alphabet
$\Acal _{-}$ of `negative' letters $\overline{1}<\overline{2} <\cdots
$ (since $\Gcal _{\nubold }(X;q)$ is symmetric, the choice of ordering
is arbitrary).

A {\em negative tableau} on $\nubold $ is a map $T\colon
\nubold\rightarrow \Acal _{-}$ that is strictly increasing on rows and
weakly increasing on columns.  Let $\SSYT _{-}(\nubold )$ be the set
of these.  An attacking inversion in a negative tableau is an
attacking pair $(a,b)$ such that $T(a)\geq T(b)$ (like for ordinary
tableaux except that equal negative entries also count as inversions).
The number of attacking inversions is again denoted $\inv (T)$.

\begin{prop}[{\cite[Corollary~4.1.3]{BHMPS-paths}}]\label{prop:omega-G-nu}
Setting $x_{\overline{i}} = x_{i}$ for indices $\overline{i}\in \Acal
_{-}$, we have
\begin{equation}\label{e:omega-G-nu}
\omega \, \Gcal_{\nubold }(X;q) = \sum _{T\in \SSYT_{-} (\nubold
)}q^{\inv (T)}\xx ^{T}.
\end{equation}
\end{prop}

As in \cite{BHMPS-paths}, formula \eqref{e:omega-G-nu} leads to the
following corollary.

\begin{cor}[{\cite[Lemma~4.1.6]{BHMPS-paths}}]\label{cor:G-nu-width}
The LLT polynomial $\Gcal _{\nubold }(X;q)$ is a linear combination of
Schur functions $s_{\lambda }$ with $\ell (\lambda )\leq l$, where
$l=\sum_i\ell(\nu_{(i)})$
is the total number of rows in the skew diagrams $\nu _{(i)}$.
\end{cor}

We also need the invariance of $\Gcal _{\nubold }(X;q)$ under shifted
rotations of $\nubold $.

\begin{prop}\label{e:G-rotation}
Let $\nubold = (\nu _{(1)},\ldots,\nu _{(k)})$ and set $\nubold ' =
(\nu _{(j+1)}^{+},\ldots,\nu _{(k)}^{+}, \nu _{(1)},\ldots,\nu
_{(j)})$ for any $1\leq j <k$, where $\nu _{(i)}^{+}$ is a translation
of $\nu _{(i)}$ such that the content of every box is increased by
$1$.  Then $\Gcal _{\nubold }(X;q) = \Gcal _{\nubold '}(X;q)$.
\end{prop}

\begin{proof}
From the construction of $\nubold '$, there is a natural bijection
between boxes of $\nubold $ and boxes of $\nubold '$ that preserves
the reading order and the set of attacking pairs.  This induces a
bijection $\SSYT (\nubold )\cong \SSYT (\nubold ')$ that preserves
$\xx ^{T}$ and $\inv (T)$.
\end{proof}

\subsection{Catalanimals}
\label{ss:Catalanimals}

Let $l$ be a positive integer and let $R_{+} = R_{+}(\GL _{l}) =
\{\alpha _{ij} = \varepsilon _{i} - \varepsilon _{j} \mid i<j \}$ be
the set of positive roots for $\GL _{l}$, where
$\varepsilon _{i}$ denotes the $i$-th unit coordinate vector in $\ZZ
^{l}$.  Given subsets $R_{q}$, $R_{t}$, $R_{qt} \subseteq R_{+}$ and a
weight $\lambda \in \ZZ ^{l}$, we define the {\em Catalanimal}
$H(R_{q},R_{t},R_{qt},\lambda )$ of length $l$ as in \cite{BHMPS-llt}
to be the symmetric rational function in $l$ variables $\zz =
z_{1},\ldots,z_{l}$ given by
\begin{equation}\label{e:Catalanimal-def}
H(R_q,R_t,R_{qt},\lambda) \defeq \sum_{w \in S_l} w\bigg(\frac{\zz
^\lambda \prod_{\alpha \in R_{qt}} (1 - q\, t\, \zz ^\alpha )}
{\prod_{\alpha \in R_+} (1 - \zz ^{-\alpha}) \prod_{\alpha \in R_q} (1
- q\, \zz ^\alpha) \prod_{\alpha \in R_t} (1 - t\, \zz
^\alpha)}\bigg),
\end{equation}
where $\zz ^\lambda$ stands for $z_{1}^{\lambda_{1}}\cdots
z_{l}^{\lambda_{l}}$.  The defining formula can also be written
\begin{equation}\label{e:Catalanimal-raising}
H(R_q,R_t,R_{qt},\lambda) = \sigmabold \bigg(\frac{\zz ^\lambda
\prod_{\alpha \in R_{qt}} (1 - q\, t\, \zz ^\alpha )} {\prod_{\alpha
\in R_q} (1 - q\, \zz ^\alpha) \prod_{\alpha \in R_t} (1 - t\, \zz
^\alpha)}\bigg),
\end{equation}
where 
\begin{equation}\label{e:Weyl-symmetrizer}
\sigmabold (f) = \sum _{w\in S_{l}} w \bigg( \frac{f}{\prod _{\alpha
\in R_{+}} (1 - \zz ^{-\alpha })} \bigg)
\end{equation}
is the Weyl symmetrization operator for $\GL _{l}$.  Recall that
$\sigmabold (\zz ^{\lambda }) = \chi _{\lambda }$ is an irreducible
$\GL _{l}$ character if $\lambda $ is a dominant weight.
For an arbitrary weight $\mu \in \ZZ ^{l}$, either $\sigmabold (\zz
^{\mu }) = \pm \chi _{\lambda }$ for a suitable dominant weight
$\lambda $, or $\sigmabold (\zz ^{\mu }) = 0$.

Expanding the factors $(1 - q\, \zz ^\alpha)^{-1} = 1 + q\, \zz
^{\alpha } + \cdots $ and $(1 - t\, \zz ^\alpha)^{-1} = 1 + t\, \zz
^{\alpha } + \cdots $ as geometric series before applying $\sigmabold
$, we can regard \eqref{e:Catalanimal-raising} as a {\em raising
operator series}, expressing $H(R_{q},R_{t},R_{qt},\lambda )$ as an
infinite formal linear combination $\sum _{\mu }a_{\mu } \chi _{\mu }$ of
irreducible $\GL _{l}$ characters with coefficients $a_{\mu }\in \ZZ [q,t]$.

The {\em polynomial characters} of $\GL _{l}$ are the irreducible
characters $\chi _{\mu }$ with $\mu \in \NN ^{l}$; thus $\mu $ is an
integer partition with at most $l$ parts and possible trailing zeroes,
and $\chi _{\mu }$ is equal to the Schur function $s_{\mu
}(z_{1},\ldots,z_{l})$.  The {\em polynomial part}
\begin{equation}\label{e:polynomial-part}
H(R_{q},R_{t},R_{qt},\lambda )_{\pol }
\end{equation}
of a Catalanimal is the truncation of its raising operator series
expansion to terms $a_{\mu }\chi _{\mu }$ for polynomial characters
$\chi _{\mu }$.  Then $H(R_{q},R_{t},R_{qt},\lambda )_{\pol }$ is a
symmetric polynomial, homogeneous of degree $|\lambda | = \sum
_{i}\lambda _{i}$, in the variables $z_{1},\ldots,z_{l}$.

We will need several results from \cite{BHMPS-llt} concerning
Catalanimals and their connection with the elliptic Hall algebra
$\Ecal $ of Burban and Schiffmann \cite{BurbSchi12} (or Schiffmann
algebra).  Before stating them, we fix notation and recall some facts
about $\Ecal $.

For each pair of coprime integers $(m,n)$, the Schiffmann algebra
$\Ecal $ contains a distinguished subalgebra $\Lambda (X^{m,n}) $
isomorphic to the algebra of symmetric functions.  By a theorem of
Schiffmann and Vasserot \cite{SchiVass13}, the `right half-plane'
subalgebra $\Ecal ^{+}\subseteq \Ecal $ generated by the $\Lambda
(X^{m,n})$ for $m>0$ is isomorphic to the shuffle algebra of Feigin et
al.\ \cite{FeHaHoShYa09}, Feigin and Tsymbauliak \cite{FeigTsym11},
and Negut \cite{Negut14}.  This shuffle algebra has many different
realizations.  Here we use the realization $\Scal _{\hGamma }$ in
\cite[\S 3.2]{BHMPS-llt}; it is a graded algebra whose degree $l$
component is a certain subspace $\Scal _{\hGamma }^{l} \subseteq \kk
(z_{1},\ldots,z_{l})^{S_{l}}$ of the space of symmetric rational
functions in $l$ variables.  The isomorphism that we use is the one
denoted
\begin{equation}\label{e:shuffle-to-Schiffmann}
\psi _{\hGamma }\colon \Scal _{\hGamma } \xrightarrow{\simeq } \Ecal
^{+}
\end{equation}
in \cite[equations (26, 28)]{BHMPS-llt}.  In \cite{SchiVass13}, Schiffmann
and Vasserot also constructed an action $\Ecal $ on $\Lambda (X)$.  We use
the version of this action given by \cite[Proposition~3.3.1]{BHMPS-paths}.

The results stated below summarize everything we need to know for the
purposes of this paper about the above algebras, isomorphism and
action.  For more details, the reader can consult
\cite{BurbSchi12,Negut14,SchiVass13}; the translation between the
notation in these papers and ours can be found in \cite[\S \S
3.2--3.3]{BHMPS-paths}, 
the defining relations of $\Ecal $ written in
our notation are in \cite[\S 3.2]{BHMPS-delta}, and the relationship
between $\Scal _{\hGamma }$ and the shuffle algebra studied by Negut
in \cite{Negut14} is explained in \cite[\S 3.6]{BHMPS-llt}.

\begin{prop}[{\cite[Proposition~4.1.3]{BHMPS-llt}}]
\label{prop:tame-Catalanimals}

Let $H = H(R_{q},R_{t},R_{qt},\lambda )$ be a {\em tame} Catalanimal
as in \cite[Definition~4.1.2]{BHMPS-llt}---that is, the root sets
satisfy $[R_{q},R_{t}] \subseteq R_{qt}$, where $[R_{q},R_{t}] =
(R_{q}+R_{t})\cap R_{+}$.  Then $H$, considered as a symmetric
rational function, is an element of $\Scal _{\hGamma }$, and as such
represents an element $\psi _{\hGamma }(H)\in \Ecal ^{+}$ of the
Schiffmann algebra.
\end{prop}

\begin{prop}\label{prop:H-pol}
If a Catalanimal $H = H(R_{q},R_{t},R_{qt},\lambda )$ of length $l$
belongs to $\Scal _{\hGamma }$, and $\zeta = \psi _{\hGamma }(H)$ is
the corresponding element of $\Ecal ^{+}$, then $\zeta $ acting on
$1\in \Lambda (X)$ satisfies
\begin{equation}\label{e:H-pol}
\omega (\zeta \cdot 1)(z_{1},\ldots,z_{l}) = H_{\pol }.
\end{equation}
In addition, $\omega (\zeta \cdot 1)$ is a linear combination of Schur
functions $s_{\mu }$ indexed by partitions $\mu $ with at most $l$
parts, so it is determined by \eqref{e:H-pol}.
\end{prop}

\begin{proof}
This follows from \cite[Proposition~3.5.2]{BHMPS-paths} in the same
way that \cite[Proposition~3.5.2]{BHMPS-llt} does.
\end{proof}

\begin{prop}[{\cite[Lemma~3.5.1]{BHMPS-llt}}]\label{prop:nabla-m-f}
For any symmetric function $f$, the element $f[-M X^{m,1}]\in \Ecal $
acting on $1\in \Lambda (X)$ is given by
\begin{equation}\label{e:nabla-m-f}
f[-M X^{m,1}]\cdot 1 = \nabla^{m} f(X).
\end{equation}
\end{prop}

\subsection{Schur Catalanimals}
\label{ss:schur-catalanimals}

Given any LLT polynomial $\Gcal _{\nubold }(X;q)$, we constructed
Catalanimals $H^{m,n}_{\nubold ^{m}}$ in \cite{BHMPS-llt} such that
$\psi _{\hGamma }(H^{m,n}_{\nubold ^{m}})$ is equal, up to a sign and
a monomial factor in $q,t$, to $\Gcal _{\nubold }[-M X^{m,n}]$. For
the proof of the Loehr-Warrington conjecture we need the special case
of this result when $\nubold $ is a single diagram and the LLT
polynomial $\Gcal _{\nubold }(X;q)$ is a Schur function.  To describe
this case, we recall some combinatorial notions from \cite[\S \S
8.1--8.2]{BHMPS-llt}.

Given $(m,n)\in \ZZ _{+}\times \ZZ $, we define the sequence of $m$
integers as in \cite[(104)]{BHMPS-llt}
\begin{equation}\label{e:bbold}
\bb (m,n)_{i} = \lceil i n/m \rceil -\lceil (i-1) n/m \rceil \qquad
(i=1,\ldots,m).
\end{equation}

The {\em $m$-stretching} of a (skew) diagram $\nu $ is the skew
diagram $\nu ^{m}$ constructed by dilating $\nu $ vertically by a
factor of $m$ in the following way: for each box $x$ of content $c$ in
$\nu $, the $m$-stretching $\nu ^{m}$ has $m$ boxes of contents $mc -
m + 1, \ldots , mc - 1, mc$ in the same column as $x$.  For example,
the $3$-stretching of the partition diagram $\nu = (3,2)$ is shown
here, with shaded bars showing the three boxes in $\nu ^{3}$ that
correspond to each box in $\nu $.
\begin{equation}\label{e:3-stretching-pic}
\nu =
\begin{tikzpicture}[scale=.3,baseline=.2cm]
\draw (2,0) -- (2,2) -- (0,2) -- (0,0) -- (3,0) -- (3,1) -- (0,1);
\draw (1,0) -- (1,2);
\end{tikzpicture}
\qquad
\nu ^{3} = 
\begin{tikzpicture}[scale=.3,baseline=.2cm]
\filldraw [black!25] (.25,3.25) rectangle (.75,5.75);
\filldraw [black!25] (.25,0.25) rectangle (.75,2.75);
\filldraw [black!25] (1.25,-1.75) rectangle (1.75,.75);
\filldraw [black!25] (1.25,1.25) rectangle (1.75,3.75);
\filldraw [black!25] (2.25,-3.75) rectangle (2.75,-1.25);
\draw (2,0) -- (0,0) -- (0,6) -- (1,6) -- (1,-2) -- (3,-2) -- (3,-1) -- (1,-1);
\draw (0,4) -- (2,4) -- (2,-4) -- (3,-4) -- (3,-2);
\draw (0,5) -- (1,5); \draw (0,3) -- (2,3); \draw (0,2) -- (2,2);
\draw (0,1) -- (2,1); \draw (2,-3) -- (3,-3);
\end{tikzpicture}
\end{equation}
We define 
\begin{equation}\label{e:gamma(nu)}
\gamma (\nu ) = (\gamma _{1},\ldots,\gamma _{h})
\end{equation}
to be the sequence of lengths of diagonals $\{x\in \nu \mid c(x) = c
\}$ in increasing order of the content $c$, and set
\begin{equation}\label{e:n'(gamma)}
n'(\gamma (\nu )) = \sum _{i=1}^{h} \binom{\gamma _{i}}{2}.
\end{equation}
The {\em magic number} $p(\nu )$ is the sum of the lengths of the
diagonals that do not contain the first box in a row of $\nu $.  Note
that the diagonals of $\nu ^{m}$ correspond to diagonals of $\nu $,
each repeated $m$ times vertically.  Using this one sees that
$n'(\gamma (\nu ^{m})) = m\, n'(\gamma (\nu ))$ and $p(\nu ^{m}) =
p(\nu )$.

A more subtle property of the magic number, which follows from
\cite[Lemma~7.2.2]{BHMPS-llt}, is that if $\nu ^{\circ }$ is the
$180^{\circ }$ rotation of $\nu $, then $p(\nu ^{\circ }) = p(\nu )$.

\begin{defn}\label{def:Schur-Catalanimal}
Given a (skew) diagram $\nu $ and coprime integers $m$, $n$ with
$m>0$, the {\em (skew) Schur Catalanimal} $H^{m,n}_{\nu ^{m}} =
H(R_{q},R_{t},R_{qt},\lambda )$ is the tame Catalanimal of length $l =
m\, |\nu | = |\nu ^{m}|$ constructed as follows, where $\nu ^{m}$ is
the $m$-stretching of $\nu $.

The root sets and weight are defined with reference to the
partition of $[l] = \{1,\ldots,l \}$
into intervals of lengths $\gamma (\nu ^{m})$.
For the root sets, we take $\alpha _{ij}\in R_{q} = R_{t}$ if $i<j$
are in distinct blocks of this partition, and $\alpha _{ij}\in
R_{qt}$ if $i<j$ are in distinct, non-adjacent blocks; equivalently,
$R_{qt} = [R_{q},R_{t}]$.

The weight $\lambda$ is defined to be constant on blocks, as follows:
for every $i \in [l]$ belonging to the $k$-th block of the partition,
we set
\begin{multline}\label{e:weight-lambda}
\lambda _{i} = \chi (\text{\rm $D_{k}$ contains the first box in
a row of $\nu ^{m}$})\\
 - \chi (\text{\rm $D_{k}$ contains the last box in a row of
$\nu ^{m}$}) + \bb(m,n)_{\modm _{m}(c)},
\end{multline}
where $D_{k}$ is the $k$-th diagonal of $\nu ^{m}$, $c$ is the content
of boxes on that diagonal, and $\modm _{m}(c)$ is the integer $j\in
[m]$ such that $j\equiv c\pmod{m}$.

When $\mu $ is a (non-skew) partition diagram, we call
$H^{m,n}_{\mu^{m}}$ a {\em Schur Catalanimal} and $H^{m,n}_{(\mu
^{\circ })^{m}}$ the {\em opposite Schur Catalanimal}, where $\mu
^{\circ }$ is the $180^{\circ }$ rotation of $\mu $.
\end{defn}

\begin{thm}[{\cite[Theorem~8.3.1]{BHMPS-llt}}]\label{thm:Schur-Catalanimal}
The Schur Catalanimal $H^{m,n}_{\mu ^{m}}$ satisfies the identity
\begin{equation}\label{e:Schur-Catalanimal}
s_{\mu }[-M X^{m,n}] = (-1)^{p(\mu )}(q\, t)^{p(\mu ) + m \, n'(\gamma
(\mu ))} \psi _{\hGamma }(H^{m,n}_{\mu ^{m}}).
\end{equation}
This identity also holds with the opposite Schur Catalanimal
$H^{m,n}_{(\mu ^{\circ})^{m}}$ in place of $H^{m,n}_{\mu ^{m}}$.
\end{thm}

\section{Nests in a den formula}
\label{s:dens-and-nests}

Our main combinatorial result, Theorem~\ref{thm:nest-identity} below,
is an identity expanding the polynomial parts of certain tame
Catalanimals as weighted sums of LLT polynomials indexed by
configurations of nested lattice paths.  In this section, we define
the required combinatorial notions and then state the theorem.  The
proof will be given in \S \ref{s:main-theorems}.

\subsection{Dens and nests}
\label{ss:dens-and-nests}

We begin by defining the data that will serve as input to
Theorem~\ref{thm:nest-identity}.

\begin{defn}\label{def:den}
A {\em den} is a tuple $(h,p,\dd ,\ee )$, where $h$ is a positive
integer, $p$ is an irrational real number, and $\dd
=(d_{0},\ldots,d_{h})$ and $\ee = (e_{0},\ldots,e_{h})$ are sequences
of integers, subject to the following conditions:
\begin{gather}
\label{e:den-d-decr}
(d_{i} - d_{j} +1)/(j-i) > p \quad \text{for} \quad 0\leq i<j\leq h-1;\\
\label{e:den-e-incr}
(e_{i} - e_{j} - 1)/(j-i) < p \quad \text{for} \quad  1\leq i<j\leq h;\\
\label{e:den-d-e}
d_{0}>e_{0},\quad d_{h} < e_{h}, \quad \text{and}\quad \sum _{i=0}^{h}
d_{i} = \sum _{i=0}^{h} e_{i}.
\end{gather}
\end{defn}

The reason for assuming $p$ irrational is to avoid having to
disambiguate equalities that might otherwise occur in comparisons such
as those in (\ref{e:den-d-decr},~\ref{e:den-e-incr}).

With any den we also define the following auxiliary notions.  The
lattice points $(i,d_{i})$ are {\em heads}, and $(i,e_{i})$ are {\em
feet}.  Points $\{(i,j)\mid e_{i}<j\leq d_{i} \}$ weakly below a head
and strictly above a foot on the same vertical line $x=i$ are {\em
sources}.  Points $\{(i,j)\mid d_{i}<j\leq e_{i} \}$ weakly below a
foot and strictly above a head are {\em sinks}.  We also set
\begin{equation}\label{e:g-vector}
\gG = (g_{1},\ldots,g_{h}),\quad \text{where}\quad g_{k} = \sum
_{i=0}^{k-1}(d_{i} - e_{i}).
\end{equation}

\begin{figure}
\[
\begin{array}{c@{\qquad\qquad}c}
\begin{tikzpicture}[scale=.5]
\draw [black!25] (0,0) grid (7,9);
\draw [fill] (0,9) circle (.2);
\draw (0,8) circle (.2);
\draw [fill] (1,7) circle (.2);
\draw (1,6) circle (.2);
\draw [fill] (2,6) circle (.2);
\draw (2,4) circle (.2);
\draw [fill, fill opacity=.3] (3,5) circle (.2);
\draw [fill] (4,4) circle (.2);
\draw (4,5) circle (.2);
\draw [fill, fill opacity=.3] (5,3) circle (.2);
\draw [fill] (6,2) circle (.2);
\draw (6,3) circle (.2);
\draw [fill] (7,0) circle (.2);
\draw (7,2) circle (.2);
\end{tikzpicture}
&
\begin{tikzpicture}[scale=.5]
\draw [black!25] (0,0) grid (7,9);
\draw [fill] (0,9) circle (.2);
\draw [fill] (1,7) circle (.2);
\draw [fill] (2,6) circle (.2);
\draw [fill] (2,5) circle (.2);
\draw (4,5) circle (.2);
\draw (6,3) circle (.2);
\draw (7,2) circle (.2);
\draw (7,1) circle (.2);
\end{tikzpicture}
\\
\begin{tikzpicture}[scale=.5,baseline=0]
\draw [fill] (0,2) circle (.2);
\node [right] at (0,2) {$\  \text{head}$};
\draw (0,1) circle (.2);
\node [right] at (0,1) {$\ \text{foot}$};
\draw [fill, fill opacity=.3] (0,0) circle (.2);
\node [right] at (0,0) {$\  \text{both}$};
\end{tikzpicture}
&
\begin{tikzpicture}[scale=.5,baseline=0]
\draw [fill] (0,2) circle (.2);
\node [right] at (0,2) {$\  \text{source}$};
\draw (0,1) circle (.2);
\node [right] at (0,1) {$\  \text{sink}$};
\end{tikzpicture}
\end{array}
\]
\caption{\label{fig:generic-den}%
Heads, feet, sources and sinks in the den with $h = 7$, $\dd =
(9,7,6,5,4,3,2,0)$, $\ee = (8,6,4,5,5,3,3,2)$.  These data define a
valid den for any $p\in (1, \frac{6}{5})$.}
\end{figure}
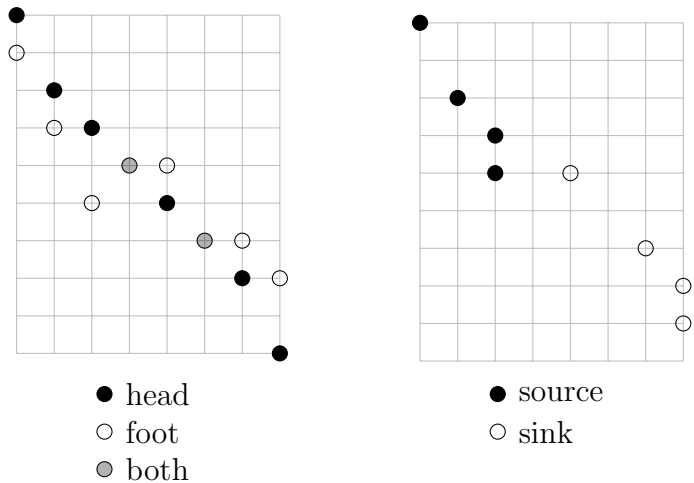

A den can be pictured by plotting the heads, feet, sources and sinks,
as shown for example in Figure~\ref{fig:generic-den}.  On each line
$x=i$ for $0\leq i\leq h$, the head, foot, and any sources or sinks
are arranged in one of the ways shown here.
\begin{equation}\label{e:head-over-heels}
\begin{array}{c@{\qquad\qquad}c@{\qquad}c}
\begin{tikzpicture}[scale=.5]
\draw [black!50] (0,-.5) -- (0,3.5);
\draw [fill] (0,3) circle (.15) node [right] {\ head};
\draw [fill] (0,2) circle (.15);
\draw [fill] (0,1) circle (.15);
\draw (0,0) circle (.15);
\draw [decorate,decoration={brace,amplitude=4pt}]
(-.4,.8) -- (-.4,3.2) node [left] at (-.75,2) {sources};
\node [right] at (.75,2) {\hphantom{sources}};
\draw (0,0) circle (.15) node [right] {\ foot};
\end{tikzpicture}
&
\begin{tikzpicture}[scale=.5]
\draw [black!50] (0,-.5) -- (0,3.5);
\draw [fill,fill opacity=.3] (0,1.5) circle (.15);
\node [right] at (.25,1.5) {head = foot};
\end{tikzpicture}
&
\begin{tikzpicture}[scale=.5]
\draw [black!50] (0,-.5) -- (0,3.5);
\draw (0,3) circle (.15) node [right] {\ foot};
\draw (0,2) circle (.15);
\draw (0,1) circle (.15);
\draw [fill] (0,0) circle (.15);
\draw [decorate,decoration={brace,amplitude=4pt}]
(-.4,.8) -- (-.4,3.2) node [left] at (-.75,2) {sinks};
\node [right] at (.75,2) {\hphantom{sinks}};
\draw (0,0) circle (.15) node [right] {\ head};
\end{tikzpicture}
\\
d_{i}> e_{i} &	d_{i} = e_{i} \hspace{12ex}
& d_{i}< e_{i}
\end{array}
\end{equation}
Condition \eqref{e:den-d-decr} means that for heads $P$ left of $Q$,
excluding the last head $(h,d_{h})$, some line of slope $-p$ passes
above $Q$ and below $P+(0,1)$.  Similarly, condition
\eqref{e:den-e-incr} means that for feet $P$ left of $Q$, excluding
the first foot $ (0,e_{0})$, some line of slope $-p$ passes above $P$
and below $Q+(0,1)$, as pictured here.
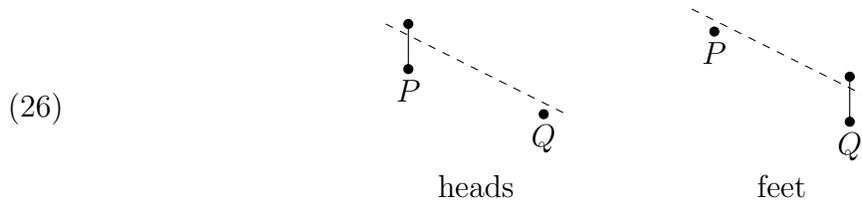
\begin{equation}\label{e:head-foot-pic}
\begin{array}{c@{\qquad\qquad}c}
\begin{tikzpicture}[scale=.6,baseline=2mm]
\draw [fill] (0,1) circle (.1) -- (0,2) circle (.1);
\node [below] at (0,1) {$P$};
\draw [fill] (3,0) circle (.1);
\node [below] at (3,0) {$Q$};
\draw [dashed] (-.5,2) -- (3.5,0);
\end{tikzpicture}
&
\begin{tikzpicture}[scale=.6,baseline=0]
\draw [fill] (0,1.5) circle (.1);
\node [below] at (0,1.5) {$P$};
\draw [fill] (3,-.5) circle (.1) -- (3,.5) circle (.1);
\node [below] at (3,-.5) {$Q$};
\draw [dashed] (-.5,2) -- (3.5,0);
\end{tikzpicture}
\\
\text{heads} & \text{feet}
\end{array}
\end{equation}

Condition \eqref{e:den-d-e} says that there is at least one source on
the line $x=0$ and at least one sink on the line $x=h$, and that the
total number of sources is equal to the total number of sinks.
Conditions \eqref{e:den-d-decr} and \eqref{e:den-e-incr} imply
\begin{equation}\label{e:den-growth-bound}
d_{j}-e_{j} \leq d_{i}-e_{i} + 1
\end{equation}
for all $0<i<j<h$.  If there is a sink with $x$ coordinate $i$ and a
source with $x$ coordinate $j$, then $d_{i}-e_{i}<0$ and
$d_{j}-e_{j}>0$.  In particular, no source and sink can be on the same
vertical line $x = i$, and \eqref{e:den-growth-bound} implies that all
sources are strictly left of all sinks.

Next we define the systems of nested lattice paths that will be
attached to a den.

\begin{defn}\label{def:nested-paths}
An {\em east end path} is a lattice path with south $(0,-1)$ and east
$(1,0)$ steps that ends with an east step.  East end paths $\pi $,
$\pi '$ are {\em nested with $\pi $ below $\pi '$} if
\begin{itemize}
\item [(i)] the interval $[a',b']$ of $x$-coordinates of points of
$\pi '$ is contained in the interval $[a,b]$ of $x$-coordinates of
points of $\pi $, and
\item [(ii)] for every integer $i\in [a',b']$ the respective intervals
$[v_{i},w_{i}]$ and $[v_{i}',w_{i}']$ of $y$-coordinates of points on
$\pi \cap (x=i)$ and $\pi '\cap (x=i)$ satisfy $v_{i}<v_{i}'$ and
$w_{i}<w_{i}'$.
\end{itemize}
\end{defn}

An example of a pair of nested east end paths is shown below.  Note
that nested paths can share south steps, but not east steps.
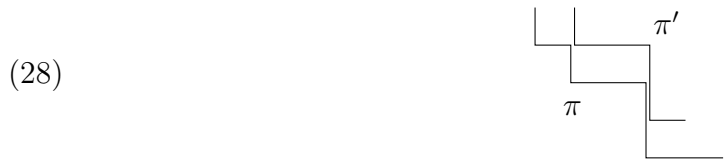
\begin{equation}\label{e:nested-pair}
\begin{tikzpicture}[scale=.5,baseline=1cm]
\draw (0,4) -- (0,3) -- (.95,3) -- (.95,2) -- (2.95,2) -- (2.95, 0) -- (5,0);
\draw (1.05,4) -- (1.05,3) -- (3.05,3) -- (3.05,1) -- (4,1);
\node [below] at (1,1.8) {$\pi $};
\node [above] at (3.5,3) {$\pi '$};
\end{tikzpicture}
\end{equation}
It is not hard to see that nesting is transitive, i.e., if $\pi $ is
nested below $\pi '$ and $\pi '$ is nested below $\pi ''$, then $\pi $
is nested below $\pi ''$.

\begin{defn}\label{def:nest}
A {\em nest} in a den $(h,p,\dd ,\ee )$ is a system of nested east end
paths $\pi = (\pi _{1},\ldots,\pi _{r})$ from the sources to the
sinks, numbered with $\pi _{k}$ nested below $\pi _{l}$ for $k<l$,
which satisfies the condition $j\leq d_{i}$ for every lattice point
$(i,j)$ other than the sink on each of the paths $\pi _{k}$.  In other
words, all non-sink lattice points in $\pi $ lie weakly below the
heads.
\end{defn}

\begin{figure}
\[
\begin{array}{c@{\qquad \qquad }c}
\begin{tikzpicture}[scale=.5]
% axes
\draw [->,gray] (-.05,-.05) -- (8,-.050);
\draw [->,gray] (-.05,-.05) -- (-.05,10);
\foreach \x in {-.05,1,2,3,4,5,6,7}
  \draw [gray] (\x ,-.050) -- (\x,-.2);
\foreach \y in {-.05,1,2,3,4,5,6,7,8,9}
  \draw [gray](-.05,\y) -- (-.2,\y);
% nest
\draw [fill] (0,9) circle (.2);
\draw (0,9) -- (0,5) -- (.95,5) -- (.95,2) -- (2,2) -- (2,1) -- (7,1);
\draw [fill] (1,7) circle (.2);
\draw [dotted, thick] (1.05,7) -- (1.05,4) -- (2,4) -- (2,3)
 -- (3,3) -- (3,2) -- (7,2);
\draw [fill] (2,6) circle (.2);
\draw [dotted, thick] (2,6) -- (2,5) -- (4,5);
\draw [fill] (2,5) circle (.2);
\draw (2,5) -- (2,4) -- (3,4) -- (3,3) -- (6,3);
\draw [fill, fill opacity=.3] (3,5) circle (.2);
\draw (4,5) circle (.2);
\draw [fill, fill opacity=.3] (5,3) circle (.2);
\draw (6,3) circle (.2);
\draw (7,2) circle (.2);
\draw (7,1) circle (.2);
% path labels on the last east steps
\node [below] at (6.5,1.1) {$\scriptstyle 1$};
\node [below] at (6.5,2.1) {$\scriptstyle 2$};
\node [below] at (5.5,3.1) {$\scriptstyle 3$};
\node [below] at (3.5,5.1) {$\scriptstyle 4$};
\end{tikzpicture}
&
\begin{tikzpicture}[scale=.5]
% axes
\draw [->,gray] (-.05,-.05) -- (8,-.050);
\draw [->,gray] (-.05,-.05) -- (-.05,10);
\foreach \x in {-.05,1,2,3,4,5,6,7}
  \draw [gray] (\x ,-.050) -- (\x,-.2);
\foreach \y in {-.05,1,2,3,4,5,6,7,8,9}
  \draw [gray](-.05,\y) -- (-.2,\y);
% nest
\draw [fill] (0,9) circle (.2);
\draw (0,9) -- (0,6) -- (.95,6) -- (.95,3) -- (2,3) -- (2,2) -- (3,2) -- (4,2)
 -- (4,1) -- (7,1);
\draw [fill] (1,7) circle (.2);
\draw [dotted, thick] (1.05,7) -- (1.05,4) -- (2,4) -- (2,3)
 -- (4,3) -- (4,2) -- (7,2);
\draw [fill] (2,6) circle (.2);
\draw [dotted, thick] (2,6) -- (2,5) -- (4,5);
\draw [fill] (2,5) circle (.2);
\draw (2,5) -- (2,4) -- (4,4) -- (4,3) -- (6,3);
\draw [fill, fill opacity=.3] (3,5) circle (.2);
\draw (4,5) circle (.2);
\draw [fill, fill opacity=.3] (5,3) circle (.2);
\draw (6,3) circle (.2);
\draw (7,2) circle (.2);
\draw (7,1) circle (.2);
% path labels on the last east steps
\node [below] at (6.5,1.1) {$\scriptstyle 1$};
\node [below] at (6.5,2.1) {$\scriptstyle 2$};
\node [below] at (5.5,3.1) {$\scriptstyle 3$};
\node [below] at (3.5,5.1) {$\scriptstyle 4$};
\end{tikzpicture}\\
\pi & \pi ^{0}\\[1ex]
\multicolumn{2}{c}{
\begin{tikzpicture}[scale=.5]
\draw [fill] (0,2) circle (.2);
\node [right] at (0,2) {$\ \text{source}$};
\draw (0,1) circle (.2);
\node [right] at (0,1) {$\ \text{sink}$};
\draw [fill, fill opacity=.3] (0,0) circle (.2);
\node [right] at (0,0) {$\ \text{head} = \text{foot}$};
\end{tikzpicture}
}
\end{array}
\]
\caption{\label{fig:nest}%
A typical nest $\pi $ in the den in Figure~\ref{fig:generic-den}, and
the unique nest $\pi ^{0}$ such that $a(\pi ^{0}) = 0$.  The sequence
$\gG = (1,2,4,4,3,3,2)$ associated with the den gives the number of
east steps at each $x$-coordinate.}
\end{figure}
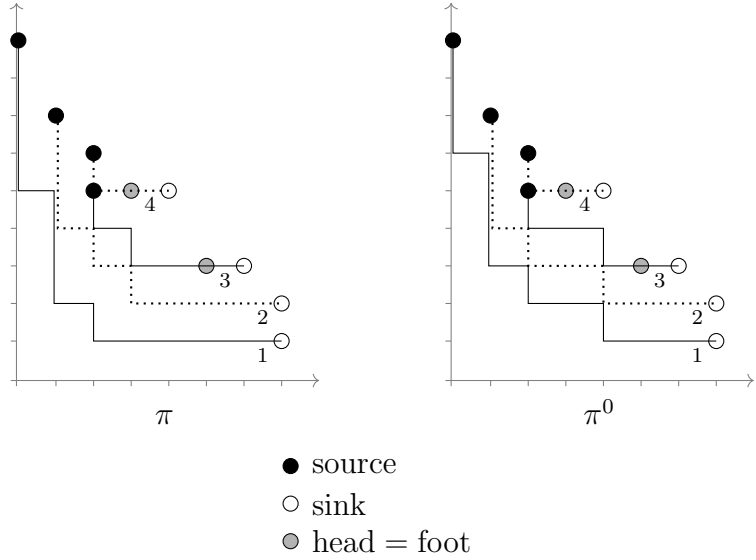

Figure~\ref{fig:nest} shows two nests belonging to the den in
Figure~\ref{fig:generic-den}.  We have marked the head ($=
\text{foot}$) on each line $x=i$ that has no source or sink in order
to make visible the condition that paths in the nest must lie weakly
below the heads.

\begin{remark}\label{rem:nest-basics}
(i) The numbering of the paths $\pi _{i}$ from nested below to nested
above implies that $\pi _{i}$ starts at the $i$-th source from left to
right, with sources on the same vertical line numbered bottom to top.
Similarly, $\pi _{i}$ ends at the $i$-th sink from right to left, again
from bottom to top on vertical lines.

(ii) Because paths in a nest are nested, any non-sink
lattice point $(i,j)$ on a path automatically lies weakly below the
head $(i,d_i)$ unless the head equals the foot $(i,e_i)$.  Only when
$d_i=e_i$, so there is no source or sink on the line $x=i$, does the
requirement $j\leq d_i$ impose an extra condition.

For example, the head and foot at $(3,5)$ in the den in
Figure~\ref{fig:generic-den}, shown as the upper left gray dot in
Figure~\ref{fig:nest}, prohibits the highest path $\pi _{4}$ from
starting with an east step to $(3,6)$, although nesting alone would
allow this.  Similarly, the head and foot at $(5,3)$ prohibits $\pi
_{3}$ from passing through $(5,4)$.

(iii) The number $g_{k}$ in \eqref{e:g-vector} counts sources minus
sinks with $x$-coordinate less than $k$.  Hence, for every nest $\pi $
in the den, $g_{k}$ is the number of paths $\pi _{j}$ in $\pi $ that
have an east step from $x = k-1$ to $x=k$, or equivalently, have a
non-sink lattice point
on the line $x = k-1$.  These are the first $g_{k}$ paths $\pi
_{1},\ldots,\pi _{g_{k}}$, since if $\pi _{j}$ has an east step from
$x = k-1$ to $x=k$, then so does every path $\pi _{i}$ nested below
$\pi _{j}$.  Since all sources are left of all sinks, and there is
assumed to be at least one source at $x = 0$ and at least one sink at
$x = h$,
the sequence $\gG $ is positive and unimodal, with maximum equal to
the number of sources (or sinks).  These properties can be seen in
Figure~\ref{fig:nest}.

(iv) It is possible to have an `abandoned' den with no nests.  The most
obvious way this can happen is if some source is lower than the
matching sink.
\end{remark}

\subsection{Parameterizing nests}
\label{ss:parametrization}

Nests in a den can be parameterized by tuples of partitions satisfying
certain inequalities, as follows.

\begin{lemma}\label{lem:parametrization}
Let $(h,p,\dd ,\ee )$ be a den, with $\gG $ as in \eqref{e:g-vector}.

\noindent (i) If $\pi = (\pi_1, \ldots, \pi_r)$ is a nest in the den, then
\begin{enumerate}
\item[(a)] for $1\leq k \leq h$, $\pi_i$ has a non-sink lattice point
on the line $x = k-1$ if and only if $1 \leq i \leq g_k$;
\item[(b)] there are unique partitions $\lambda _{(1)},\ldots,\lambda
_{(h-1)}$ of length $\ell (\lambda_{(k)})\leq \min(g_k, g_{k+1})$ such
that for $1\leq k \leq h$ and $1 \leq i \leq g_k$, the $y$-coordinates
of all non-sink lattice points of $\pi_i$ on the line $x = k-1$ form
the interval
\begin{equation}\label{e:pi(i)-interval}
I_{k,i} = [e_{k}-g_{k}+i - (\lambda _{(k)})_{i},\,
e_{k-1}-g_{k-1}+i-(\lambda_{(k-1)})_{i}]
\end{equation}
where we set $\lambda _{(0)} = \lambda _{(h)} = \varnothing $, extend
partitions with trailing zeroes if needed, and set $g_{0} = 0$.
\end{enumerate}

\noindent (ii) Set $\lambda _{(0)} = \lambda _{(h)} = \varnothing $
and $g_{0} = 0$. Let $\lambda _{(1)},\ldots,\lambda _{(h-1)}$ be
partitions of length
$\ell (\lambda_{(k)})\leq \min(g_k, g_{k+1})$
such that for $1\leq k \leq h$ and $1\leq i\leq g_{k}$,
\begin{equation}\label{e:interval-condition}
e_{k}-g_{k} - (\lambda _{(k)})_{i}\leq e_{k-1}-g_{k-1}-(\lambda
_{(k-1)})_{i}
\end{equation}
(so the intervals $I_{k,i}$ are non-empty).  Then there is a unique
nest $\pi = (\pi_1, \ldots, \pi_r)$ in the den such that the
$y$-coordinates of all non-sink lattice points of $\pi_i$ on the line
$x=k-1$ form the interval $I_{k,i}$ for all $1\leq k \leq h$ and
$1\leq i\leq g_{k}$.
\end{lemma}

\begin{proof}
Given a nest $\pi $, part (i)(a) holds by Remark~\ref{rem:nest-basics}
(iii), and the east steps from $x=k-1$ to $x=k$ in $\pi $ are on paths
$\pi _{1}$ through $\pi _{g_{k}}$.  Let $y_{k,1}<\ldots<y_{k,g_{k}}$
be the $y$-coordinates of these east steps.  The right endpoint of any
east step is weakly below the foot on the same vertical line, so
$y_{k,g_{k}}\leq e_{k}$.  Hence,
\begin{equation}\label{e:y-bound}
y_{k,i}\leq e_{k}-g_{k}+i
\end{equation}
for all $k$ and $i\leq g_{k}$.

Every east step from $x=h-1$ to $x=h$ ends at a sink, so $y_{h,i} =
e_{h}-g_{h}+i$ is fixed for all $i\leq g_{h}$, independent of the
nest.  If $k<h$ and $g_{k}>g_{k+1}$, there are $g_{k}-g_{k+1}$ sinks
on the line $x = k$.  In this case the top $g_{k}-g_{k+1}$ east steps
from $x=k-1$ to $x=k$ are fixed, with $y_{k,i} = e_{k}-g_{k}+i$ for
$g_{k+1} <i \leq g_{k}$.  This leaves the $y_{k,i}$ for $k<h$ and
$i\leq \min (g_{k},g_{k+1})$ free to vary with the nest.
To establish part (i)(b), let
\begin{equation}\label{e:lambda-parameters}
\lambda _{(k)} = (e_{k} - g_{k} + 1 - y_{k,1},\ldots,e_{k} - g_{k} + r_k
- y_{k,r_k}),
\end{equation}
for $k = 1,\ldots,h-1$,
where $r_k = \min (g_{k},g_{k+1})$.  We also set $\lambda _{(0)} =
\lambda _{(h)} = \varnothing $ and $(\lambda _{(k)})_{i} = 0$ for $i>
r_k$.

Since the $y_{k,i}$ are strictly increasing and bounded by
\eqref{e:y-bound}, $\lambda _{(k)}$ is a partition of length $\ell
(\lambda _{(k)})\leq \min (g_{k},g_{k+1})$ with possible trailing
zeroes.  Since equality holds in \eqref{e:y-bound} for the fixed east
steps that end at sinks, we have
\begin{equation}\label{e:y-from-lambda}
y_{k,i} = e_{k} - g_{k} + i -(\lambda _{(k)})_{i}
\end{equation}
for all $1 \leq k \leq h$ and $i\leq g_{k}$.

The lattice point at $x = k-1$ on $\pi _{i}$ with the smallest
$y$-coordinate is the left endpoint of the east step with
$y$-coordinate $y_{k,i}$.  If $k>1$ and $i\leq g_{k-1}$, the point
with the largest $y$-coordinate is the right endpoint of the east step
with $y$-coordinate $y_{k-1,i}$.  Otherwise, if $k=1$ or
$g_{k-1}<i\leq g_{k}$, this highest point is the $i$-th source, with
$y$-coordinate $e_{k-1} - g_{k-1} + i$, if we take $g_{0} = 0$ in the
case $k=1$.  This shows that the $y$-coordinates of the points on $\pi
_{i}$ at $x = k-1$ are given by the interval $I_{k,i}$ in all cases.
These intervals clearly determine the partitions $\lambda _{(k)}$.

For part (ii),
suppose we are given partitions $\lambda _{(1)},\ldots,\lambda
_{(h-1)}$ with $\ell (\lambda _{(k)})\leq \min (g_{k},g_{k+1})$ such
that \eqref{e:interval-condition} holds.  By
Remark~\ref{rem:nest-basics} (iii), the sequence $\gG $ is positive
and unimodal with maximum equal to the number $r$ of sources (or
sinks) in the den, and for each $i =1,\ldots, r$
the set $\{k\in [h]
\mid g_{k}\geq i \}$ is the non-empty interval $[k_{0},k_{1}]$ such
that the $i$-th source and its matching sink are at $x = k_{0}-1$ and
$x = k_{1}$.  For each $i$, we start by constructing an east end path
$\pi_{i}$ such that
the intervals $I_{k,i}$ for $k \in [k_0, k_1]$ describe the lattice points
on $\pi_{i}$.

Fix $i$ and the corresponding interval $[k_{0},k_{1}]$.  The $i$-th
source has $y$-coordinate $e_{k_{0}-1} - g_{k_{0}-1} + i$, if we set
$g_{0} = 0$ for $k_{0} = 1$.  Since $k_{0}$ is minimal with $i\leq
g_{k_{0}}$, we have $i>g_{k_{0}-1}$ and therefore $(\lambda
_{(k_{0}-1)})_{i} = 0$
since \(\ell(\lambda_{(k_0-1)}) \leq  g_{k_0-1}\).
Hence, the upper endpoint of the interval $I_{k_{0},i}$
is the $y$-coordinate of the $i$-th source.  The sink matching the
$i$-th source has $y$-coordinate $e_{k_{1}} - g_{k_{1}} + i$.  Since
$k_{1}$ is maximal with $i\leq g_{k_{1}}$, either $k_{1} = h$ or
$i>g_{k_{1}+1}$, and therefore $(\lambda _{(k_{1})})_{i} = 0$.  Hence,
the lower endpoint of $I_{k_{1},i}$
is the $y$-coordinate of the matching sink.  For the rest, if
$k_{0}\leq k<k_{1}$, the lower endpoint of
$I_{k,i}$ is equal to the upper endpoint of $I_{k+1,i}$.
Hence, there exists a unique east end path $\pi _{i}$ from the $i$-th
source to the $i$-th sink with non-sink lattice points at $x = k-1$
given by the intervals $I_{k,i}$
for this $i$ and $k_{0}\leq k\leq k_{1}$.

By construction, the paths $\pi _{i}$ defined in this way are nested
with $\pi _{i}$ below $\pi _{j}$ for $i<j$.  The upper endpoint of the
interval $I_{k,i}$
is at most $e_{k-1} - g_{k-1}+g_{k}$, which is equal to $d_{k-1}$ by
the definition of $g_{k}$, so the paths $\pi _{i}$ form a nest in the
den.
\end{proof}

\begin{example}\label{ex:parametrization}
Nests in the den in Figure~\ref{fig:generic-den} are parameterized by
partitions $\lambda _{(1)},\ldots,\lambda _{(6)}$ of lengths at most
$1,2,4,3,3,2$, subject to the inequalities in
\eqref{e:interval-condition}.  The nest $\pi ^{0}$ on the right in
Figure~\ref{fig:nest} corresponds to $\lambda _{(i)} = \varnothing $
for all $i=1,\ldots,6$.  The nest $\pi $ on the left in
Figure~\ref{fig:nest} corresponds to $\lambda _{(1)} = (1)$, $\lambda
_{(2)} = (1)$, $\lambda _{(3)} = (1)$, $\lambda _{(4)} = (1,1,1)$,
$\lambda _{(5)} = \lambda _{(6)} = \varnothing $.
\end{example}

\subsection{Combinatorial statistics associated with nests}
\label{ss:nest-stats}

We now define statistics $a(\pi )$ and $\dinv _{p}(\pi )$ for each
nest $\pi $ in a den, closely related to the $\area $ and $\dinv $
statistics seen in the Loehr-Warrington conjecture and the various
generalizations of the shuffle theorem.

\begin{defn}\label{def:a(pi)}
Let $\pi $ be a nest in a den $(h,p,\dd ,\ee )$ and let $\gG $ be the
sequence in \eqref{e:g-vector}, so there are $g_{k}$ east steps from
$x = k-1$ to $x=k$ in $\pi $.  For $k = 1,\ldots,h$, let $y_{k,1}<
\cdots < y_{k,g_{k}}$ be the $y$-coordinates of these east steps.
We define
\begin{equation}\label{e:a(pi)}
a(\pi ) = \sum _{k=1}^{h}\sum _{i=1}^{g_{k}} (e_{k}-g_{k} + i -
y_{k,i}).
\end{equation}
Equivalently, by \eqref{e:y-from-lambda}, $a(\pi ) = |\lambda | = \sum
_{k} |\lambda _{(k)}|$ in terms of the parameterization in
Lemma~\ref{lem:parametrization}.
\end{defn}

If we let $\area (\pi _{i})$ denote the number of lattice squares
above $\pi _{i}$ and below some fixed boundary---for instance, the
number of lattice squares above $\pi _{i}$ in the rectangle with
corners at its source and sink---then both $a(\pi )$ and $\sum _{i}
\area (\pi _{i})$ have the form $(\text{constant} - \sum _{k,i}
y_{k,i})$.  Hence, $a(\pi )$ differs from $\sum _{i} \area (\pi _{i})$
by a constant not depending on $\pi $.

If there is a nest $\pi ^{0}$ in the den such that $a(\pi ^{0}) = 0$,
it must correspond via Lemma~\ref{lem:parametrization} to $\lambda
_{(k)} = \varnothing $ for all $k$, or equivalently to
$y_{k,i} = e_{k} - g_{k} +i$
for all $k$ and $i$.  Since this is always an upper bound on
$y_{k,i}$, the path $\pi _{i}$ in any nest $\pi $ lies weakly
below the path $\pi ^{0}_{i}$, and $a(\pi )$ is equal to the sum of
the areas $\area (\pi _{i})$ between $\pi _{i}$ and $\pi ^{0}_{i}$.
We can make this more precise as follows.

\begin{cor}\label{cor:a(pi)=0}
A den $(h,p,\dd ,\ee )$ has a nest $\pi ^{0}$ such that $a(\pi ^{0}) =
0$ if and only if $e_{k}\leq d_{k-1}$ for all $k=1,\ldots,h$.  Such a
nest $\pi ^{0}$ is unique.  If it exists, then each path $\pi _{i}$ in
any nest $\pi $ lies weakly below $\pi ^{0}_{i}$, and $a(\pi )$ is
equal to the sum of the areas $\area (\pi _{i}) = |\pi ^{0}_{i}/\pi
_{i}|$ enclosed between the paths $\pi _{i}$ and $\pi ^{0}_{i}$.

If $\dd $ is weakly decreasing, which is always the case if $p>0$ by
condition \eqref{e:den-d-decr}, then the den is either abandoned (has
no nests), or it has a nest $\pi ^{0}$ as above.
\end{cor}

\begin{proof}
By the definition of $\gG $ in \eqref{e:g-vector}, we have $e_{k}\leq
d_{k-1}$ if and only if $e_{k} - g_{k}\leq e_{k-1} - g_{k-1}$, taking
$g_{0} = 0$.  This is equivalent to
\eqref{e:interval-condition} for all $\lambda _{(k)} = \varnothing $.
The first paragraph then follows from the preceding observations and
Lemma~\ref{lem:parametrization}.

For the last part, assume that $\dd $ is weakly decreasing.
If $e_{k} > d_{k-1}$ for some $k$, it follows that $e_{k} > d_{k}$, so
$(k,e_{k})$ is a sink.  The east step ending at this sink in any nest
would start at $(k-1,e_{k})$, but this is not allowed, since $e_{k} >
d_{k-1}$.  Hence, the den is abandoned unless $e_{k}\leq d_{k-1}$ for
all $k = 1,\ldots,h$.
\end{proof}

\begin{example}\label{ex:pi-0}
The den in Figure~\ref{fig:generic-den} has a unique nest $\pi ^{0}$
with $a(\pi ^{0}) = 0$, shown on the right in Figure~\ref{fig:nest}.
For the nest $\pi $ on the left in Figure~\ref{fig:nest}, the areas
between corresponding paths $\pi ^{0}_{i}$ and $\pi _{i}$ add up to
$a(\pi ) = 4+1+1+0 = 6$.
\end{example}

\begin{defn}\label{def:dinv}
Let $\pi $ be a nest in a den $(h,p,\dd ,\ee )$.  We define $\dinv
_{p}(\pi )$ to be the number of tuples $(P, i, S, j)$, where $P$ is a
non-sink lattice point on $\pi _{i}$, $S$ is a south step on $\pi
_{j}$, $P$ is strictly left of $S$, and the line of slope $-p$ through
$P$ passes through $S$ (necessarily through the interior of $S$, since
we assume $p$ is irrational).
\end{defn}

In effect, $\dinv_{p}(\pi )$ counts pairs $P$, $S$ in $\pi $ whose
relative position is as indicated, with multiplicities if $P$ or $S$
lies on more than one path in the nest.
\begin{equation}\label{e:dinv-pair}
\begin{array}{c}
\begin{tikzpicture}[scale=.5,]
\draw [thick] (3,-.2) -- (3,.8);
\node [right] at (3,.8) {$S$};
\draw [dashed] (-.5,2) -- (3.5,0);
\draw [fill] (0,1.75) circle (.1);
\node [below left] at (0,1.75) {$P$};
\end{tikzpicture}
\end{array}
\end{equation}

\begin{figure}
\[
\begin{array}{c@{\qquad \qquad }c}
\begin{tikzpicture}[scale=.5]
\draw [fill, gray] (0,9) circle (.1);
\draw [fill, gray] (0,8) circle (.1);
\draw [fill, gray] (0,7) circle (.1);
\draw [fill, gray] (0,6) circle (.1);
\draw [fill, gray] (0,5) circle (.1);
\draw [fill, gray] (1,7) circle (.1);
\draw [fill, gray] (1,6) circle (.1);
\draw [fill, gray] (1,5) circle (.1);
\draw [fill, gray] (1,4) circle (.1);
\draw [fill, gray] (1,3) circle (.1);
\draw [fill, gray] (1,2) circle (.1);
\draw [fill, gray] (2,6) circle (.1);
\draw [fill, gray] (2,5) circle (.1);
\draw [fill, gray] (2,4) circle (.1);
\draw [fill, gray] (2,3) circle (.1);
\draw [fill, gray] (2,2) circle (.1);
\draw [fill, gray] (2,1) circle (.1);
\draw [fill, gray] (3,5) circle (.1);
\draw [fill, gray] (3,4) circle (.1);
\draw [fill, gray] (3,3) circle (.1);
\draw [fill, gray] (3,2) circle (.1);
\draw [fill, gray] (3,1) circle (.1);
\draw (4,5) circle (.1);
\draw [fill, gray] (4,3) circle (.1);
\draw [fill, gray] (4,2) circle (.1);
\draw [fill, gray] (4,1) circle (.1);
\draw [fill, gray] (3,1) circle (.1);
\draw [fill, gray] (5,3) circle (.1);
\draw [fill, gray] (5,2) circle (.1);
\draw [fill, gray] (5,1) circle (.1);
\draw (6,3) circle (.1);
\draw [fill, gray] (6,2) circle (.1);
\draw [fill, gray] (6,1) circle (.1);
\draw (7,2) circle (.1);
\draw (7,1) circle (.1);
\draw [gray] (-.2,10) -- (7.5,2);
\draw (0,9) -- (0,5) -- (.95,5) -- (.95,2) -- (2,2) -- (2,1) -- (7,1);
\draw [dotted, thick] (1.05,7) -- (1.05,4) -- (2,4) -- (2,3)
 -- (3,3) -- (3,2) -- (7,2);
\draw [dotted, thick] (2,6) -- (2,5) -- (4,5);
\draw (2,5) -- (2,4) -- (3,4) -- (3,3) -- (6,3);
\end{tikzpicture}
&
\begin{tikzpicture}[scale=.5]
% (44/33)
\draw (3,0) rectangle (4,2);
\draw (3,1) -- (4,1);
% (6555/5444) translated by (1,1)
\draw (6,1) rectangle (7,2);
\draw (5,2) rectangle (6,5);
\draw (5,3) -- (6,3);
\draw (5,4) -- (6,4);
% (65/32) translated by (5,5)
\draw (8,7) -- (8,5) -- (11,5) -- (11,6) -- (7,6) -- (7,7) -- (10,7)
 -- (10,5);
\draw (9,7) -- (9,5);
% (4/0) translated by (8,8)
\draw (8,8) rectangle (12,9);
\draw (9,8) -- (9,9);
\draw (10,8) -- (10,9);
\draw (11,8) -- (11,9);
% content diagonals
\draw [dashed,gray] (2.5,2.5) -- (8,8); \draw [dashed,gray] (9,9) -- (9.5,9.5);
\draw [dashed,gray] (2.5,1.5) -- (5,4); \draw [dashed,gray] (6,5) -- (7,6); 
\draw [dashed,gray] (8,7) -- (9,8); \draw [dashed,gray] (10,9) -- (10.5,9.5); 
\draw [dashed,gray] (2.5,.5) -- (3,1); \draw [dashed,gray] (4,2) -- (5,3); 
\draw [dashed,gray] (6,4) -- (8,6); \draw [dashed,gray] (9,7) -- (10,8); 
\draw [dashed,gray] (11,9) -- (11.5,9.5); 
\draw [dashed,gray] (2.5,-.5) -- (3,0); \draw [dashed,gray] (4,1) -- (5,2); 
\draw [dashed,gray] (6,3) -- (8,5); \draw [dashed,gray] (10,7) -- (11,8); 
\draw [dashed,gray] (12,9) -- (12.5,9.5); 
\draw [dashed,gray] (3.5,-.5) -- (9,5);
\draw [dashed,gray] (10,6) -- (12.5,8.5);
\draw [dashed,gray] (4.5,-.5) -- (6,1); \draw [dashed,gray] (7,2) -- (10,5);
\draw [dashed,gray] (11,6) -- (12.5,7.5);
\end{tikzpicture}\\
\pi & \nubold (\pi )
\end{array}
\]
\caption{\label{fig:nest-dinv-nu}%
(i) The nest $\pi $ drawn below a line $y+p\, x = s$ with $p \approx
1.04$ and $s \approx 9.80$. (ii) The tuple of skew diagrams $\nubold
(\pi )$, arranged southwest to northeast with dashed lines showing
boxes of equal content.}
\end{figure}
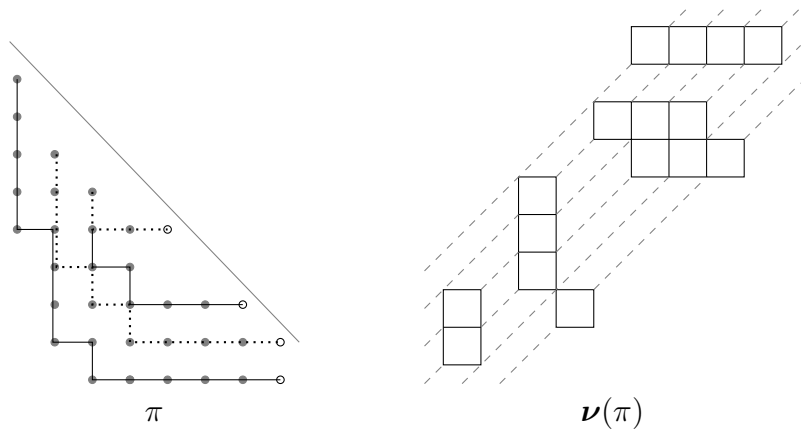

\begin{example}\label{ex:dinv}
In the den in Figure~\ref{fig:generic-den}, take $p = 1+\epsilon $ for
a small $\epsilon > 0$.  Then $P$ and $S$ contribute to $\dinv
_{p}(\pi )$ if $P$ is northwest of the upper end of $S$ on the same
diagonal of slope $-1$.  In Figure~\ref{fig:nest-dinv-nu} we have
redrawn the nest $\pi $ from Figure~\ref{fig:nest}, displaying lattice
points, as an aid to checking that $\dinv _{p}(\pi ) = 22$.
\end{example}

\subsection{LLT polynomial associated with a nest}
\label{ss:nest-LLT}

Next we define a tuple of skew diagrams $\nubold (\pi )$ attached to
each nest $\pi $ in a den.

Given a den $(h, p, \dd ,\ee )$, fix a real number $s$ such that the
line $y+px=s$ passes weakly above all the heads and feet.  For $i =
1,\ldots,h$, let $c_{i} = \fp{s - p\, (i-1)}$ be the height of the gap
between the line $y+p\, x = s$ and the highest lattice point weakly
below it at $x = i-1$, where $\fp{a} = a - \lfloor a \rfloor$ denotes
the fractional part of a real number $a$.

Let $\sigma \in S_{h}$ be the permutation such that $\sigma
(1),\ldots,\sigma (h)$ are in the same relative order as
$c_{1},\ldots,c_{h}$, i.e., such that $\sigma (c_{1},\ldots,c_{h})$ is
increasing.  Note that the $c_{i}$ are distinct, since $p$ is
irrational.

\begin{defn}\label{def:nu(pi)}
Given a den $(h,p,\dd ,\ee )$, fix $s$ and define $c_{1},\ldots,c_{h}$
and $\sigma $ as above.  For each nest $\pi $ in the den, $\nubold
(\pi )$ is the tuple of skew diagrams $(\nu _{(1)},\ldots,\nu _{(h)})$
with $\nu _{(j)}$ defined as follows.

Let $k = \sigma ^{-1}(j)$.  By Remark~\ref{rem:nest-basics} (iii), the
paths in $\pi $ which have a non-sink lattice point on the line $x =
k-1$ are $\pi _{1},\ldots,\pi _{g_{k}}$, for $g_{k}$ defined in
\eqref{e:g-vector}.  For $i = 1,\ldots,g_{k}$, let $y_{i}$ and $w_{i}$
be the maximum and minimum $y$-coordinates of lattice points on $\pi
_{i}$ at $x = k-1$.  Then we set
\begin{equation}\label{e:nu(pi)(j)}
\nu _{(j)} = \beta /\alpha = (\beta _{1},\ldots,\beta
_{g_{k}})/(\alpha _{1},\ldots,\alpha _{g_{k}}),
\end{equation}
where $\alpha _{i} = \lfloor s - p \,(k - 1) \rfloor - y_{i} + i$ and
$\beta _{i} = \lfloor s - p \, (k - 1) \rfloor - w _{i} + i$.  Note
that $(y_{1},\ldots,y_{g_{k}})$ and $(w_{1},\ldots,w_{g_{k}})$ are
strictly increasing, with $y_{r}\geq w_{r}$, so $\beta /\alpha $ makes
sense as a skew diagram.
\end{defn}

Less formally, for each $k = 1,\ldots,h$, if $O_{k} = (k-1,\lfloor s -
p\, (k-1) \rfloor )$ is the highest lattice point at $x = k-1$ below
the line $y+p\, x = s$, we construct a skew diagram by turning runs of
south steps in paths in $\pi$ at $x = k-1$ into rows of a skew
diagram, placed so that the content of the box corresponding to a
south step $S$ is the distance between $O_{k}$ and the south endpoint
of $S$.  Then $\nubold (\pi )$ is the list of these skew diagrams in
increasing order of the gaps $c_{k}$.  Note that this makes the
reading order on boxes of $\nubold (\pi )$ correspond to the
ordering of south steps in $\pi $ by increasing distance below the
line $y+p\, x = s$, with occurrences of the same south step $S$ on two
paths $\pi _{i}$, $\pi _{j}$ ordered by $i<j$.

\begin{example}\label{ex:nu(pi)}
In Figure~\ref{fig:nest-dinv-nu}, we have re-drawn the nest $\pi $
from Example~\ref{ex:dinv} and Figure~\ref{fig:nest} below a line $y +
p\, x = s$ with $p = 1 +\epsilon $ for $\epsilon \approx .04$, and
$s\approx 9.8$ chosen so the line passes a little above the source at
$(0,9)$ and the sinks at $(4,5)$, $(6,3)$ and $(7,2)$.
The gaps $c_{k}$ between this line and the highest lattice points
below it increase in the order $c_{7}<\cdots <c_{1}$, giving $\sigma =
w_{0}$, the longest permutation in $S_{7}$.  Accordingly, the skew
diagrams $\nu _{(1)},\ldots,\nu _{(7)}$ in $\nubold (\pi )$ are
associated to south runs on the lines $x = k-1$ in the order $k = 7,
6, \ldots, 1$.

The first three diagrams in $\nubold (\pi )$ are empty.  The last four
are plotted in Figure~\ref{fig:nest-dinv-nu}, arranged from southwest
to northeast and positioned so that boxes of equal content are on the
same diagonal line.
\end{example}

\begin{remark}\label{rem:nu(pi)-vs-lambda}
In terms of the parameterization in Lemma~\ref{lem:parametrization},
the skew diagram $\nu _{(j)} = \beta /\alpha $ in \eqref{e:nu(pi)(j)}
is given by $\alpha = (a^{g_{k}}) + \lambda _{(k-1)}$, $\beta
=(b^{g_{k}}) + \lambda _{(k)}$, where $k = \sigma ^{-1}(j)$, $a =
{\lfloor s - p\, (k-1) \rfloor} - e_{k-1} + g_{k-1}$, $b = {\lfloor s
- p\, (k-1) \rfloor} - e_{k} + g_{k}$.
\end{remark}

Although the definition of $\nubold (\pi )$ involves an auxiliary
choice of the line $y+p\, x = s$, one can check that when $s$ varies,
$\nubold (\pi )$ changes by rotations of the kind in
Proposition~\ref{e:G-rotation}.  The LLT polynomial $\Gcal _{\nubold
(\pi )}(X;q)$ therefore depends only on $\pi $ and the slope parameter
$p$.  We can make this more explicit as follows.

Let $\Sbold (\pi )$ be the set of pairs $(S,i)$ such that $S$ is a
south step on $\pi _{i}$, and let $\chat (S)$ denote the vertical
distance $s - (l +p\, k)$ between the line $y+p\, x = s$ and the south
endpoint $(k,l)$ of $S$.  We say that
an ordered pair of elements $(S,i),\,(S',j)\in \Sbold (\pi )$ is
an {\em attacking pair}
if $0< \chat (S') - \chat (S)< 1$.  This means that $S$ and $S'$ are
distinct, some line of slope $-p$ passes through them both, and they
are ordered with $\chat (S)<\chat (S')$.  The differences $\chat (S')
- \chat (S)$ and the set of attacking pairs do not depend on $s$.  Via
the natural correspondence between $\Sbold (\pi )$ and the set of
boxes in $\nubold (\pi )$, attacking pairs in $\Sbold (\pi )$
correspond to attacking pairs in $\nubold (\pi )$.

\begin{defn}\label{def:labelings}
A {\em negative labeling} of a nest $\pi $ is a map $N\colon \Sbold
(\pi )\rightarrow \ZZ _{+}$ that satisfies the conditions
\begin{itemize}
\item [(i)] $N$ is strictly increasing from north to south along each
run of south steps in each $\pi _{i}$;
\item [(ii)] if $(S,i)$ and $(S',i+1)$ are on the same vertical line
with $S'$ immediately above $S$, then $N(S,i)\leq N(S',i+1)$.
\end{itemize}
We define $\inv (N)$ to be the number of attacking pairs $(S,i)$,
$(S',j)$ in $\Sbold (\pi )$ such that $N(S,i)\geq N(S',j)$.

A {\em positive labeling} $P\colon \Sbold (\pi )\rightarrow \ZZ _{+}$
is defined similarly, but with `weakly increasing' in place of
`strictly increasing' in (i), with $P(S,i)< P(S',i+1)$ in place of
$N(S,i)\leq N(S',i+1)$ in (ii), and with $P(S,i) > P(S',j)$ in place of
$N(S,i)\geq N(S',j)$ when defining $\inv (P)$.
\end{defn}

By construction, if we transfer labels from steps in $\Sbold (\pi)$ to
the corresponding boxes in $\nubold (\pi )$, positive or negative
labelings $P$ or $N$ correspond to tableaux $T\in \SSYT _{\pm
}(\nubold (\pi ))$ with positive or negative letters.  The definition
of $\Gcal _{\nubold }$ and Proposition~\ref{prop:omega-G-nu} therefore
yield
\begin{gather}\label{e:pi-labeling-pos}
\Gcal _{\nubold (\pi )}(X;q) = \sum _{P} q^{\inv (P)} \xx ^{P},\\
\label{e:pi-labeling-neg}
\omega \, \Gcal _{\nubold (\pi )}(X;q) = \sum _{N} q^{\inv (N)} \xx ^{N},
\end{gather}
where the sums are over positive and negative labelings $P$ and $N$,
respectively, and $\xx ^{P} = \prod _{(S,i)\in \Sbold (\pi )}
x_{P(S,i)}$, with $\xx ^{N}$ defined similarly.

\subsection{Main theorem}
\label{ss:main-theorem}

We have now defined the ingredients needed to state our main
combinatorial result.

\begin{thm}\label{thm:nest-identity}
Given a den $(h,p,\dd ,\ee )$, with $\gG $ as in \eqref{e:g-vector},
define a Catalanimal
\begin{equation}\label{e:H-for-den}
H = H\big(R_{q}, R_{t}, R_{qt}, ((d_{0}-e_{1})^{g_{1}}, \ldots,
(d_{h-1}-e_{h})^{g_{h}})\big)
\end{equation}
of length $l = |\gG | = \sum _{k} g_{k}$, taking $R_{q} = R_{t}$ to be
the set of positive roots $\alpha _{ij}$ ($i<j$) such that $i$, $j$
are in distinct blocks of the partition of $[l]$ into intervals of
lengths $g_{k}$, and $R_{qt}$ to be the subset of these roots with
$i$, $j$ in non-adjacent blocks.

Then the polynomial part of $H$ is given by
\begin{equation}\label{e:nest-identity}
H_{\pol }(\zz ) = \sum _{\pi } t^{a(\pi )} q^{\dinv _{p}(\pi )} \Gcal
_{\nubold (\pi )}(z_{1},\ldots,z_{l};q^{-1}),
\end{equation}
where the sum is over all nests $\pi $ in the given den, and $a(\pi
)$, $\dinv _{p}(\pi )$ and $\Gcal _{\nubold (\pi )}(X;q)$ are as
defined in \S \S \ref{ss:nest-stats}--\ref{ss:nest-LLT}.
\end{thm}

There are several alternative ways to formulate the conclusion of
Theorem~\ref{thm:nest-identity}.  Using Proposition~\ref{prop:H-pol},
we can connect it with the Schiffmann algebra, as follows.  Note that
the root sets in \eqref{e:H-for-den} satisfy $R_{qt} = [R_{q},R_{t}]$,
so the Catalanimal $H$ is tame.

\begin{cor}\label{cor:nest-identity-Schiffmann}
Given a den $(h,p,\dd ,\ee )$, let $\zeta = \psi _{\hGamma }(H) \in
\Ecal ^{+}$ be the Schiffmann algebra element represented by the tame
Catalanimal $H$ in \eqref{e:H-for-den}.  Then
\begin{equation}\label{e:nest-identity-Schiffmann}
\zeta \cdot 1 = \sum _{\pi } t^{a(\pi )} q^{\dinv _{p}(\pi )} \omega\,
\Gcal _{\nubold (\pi )}(X;q^{-1}).
\end{equation}
\end{cor}

We can also reformulate \eqref{e:nest-identity} and
\eqref{e:nest-identity-Schiffmann} in terms of labeled nests using
(\ref{e:pi-labeling-pos}--\ref{e:pi-labeling-neg}).
The resulting identities are
\begin{gather}\label{e:nest-identity-bis}
H_{\pol }(\zz ) = \sum _{\pi , P} t^{a(\pi )} q^{\dinv _{p}(\pi ) - \inv
(P)} \zz ^{P},\\
\label{e:nest-identity-Schiffmann-bis}
\zeta \cdot 1 = \sum _{\pi , N} t^{a(\pi )} q^{\dinv _{p}(\pi ) - \inv
(N)} \xx ^{N},
\end{gather}
where the sums are over nests $\pi $ and positive or negative
labelings $P$ or $N$ of $\pi $ (Definition~\ref{def:labelings}),
with the positive labelings $P$ in \eqref{e:nest-identity-bis} having
labels between $1$ and $l$.

\begin{cor}\label{cor:q-t-symmetry-and-positivity} 
(i) The right hand sides of \eqref{e:nest-identity} through
\eqref{e:nest-identity-Schiffmann-bis} are symmetric in $q$ and $t$.

\noindent (ii) The left hand sides of \eqref{e:nest-identity} through
\eqref{e:nest-identity-Schiffmann-bis} are $q,t$ Schur positive, i.e.,
they are linear combinations of Schur functions with coefficients in
$\NN [q,t]$.
\end{cor}

\begin{proof}
(i) The Catalanimal $H$ in the theorem is symmetric in $q$ and $t$ by
construction.

(ii) A priori, the coefficients are in $\ZZ [q^{\pm 1},t]$, but (i)
implies that they are in $\ZZ [q,t]$.  It was shown in
\cite{GrojHaim07} that LLT polynomials $\Gcal _{\nubold }(X;q)$ are $q$
Schur positive.  Hence, the coefficients are in $\NN [q,t]$.
\end{proof}

\begin{remark}\label{rem:abandonded-den}
If the den $(h,p,\dd ,\ee )$ has no nests,
Theorem~\ref{thm:nest-identity} implies that the left hand sides of
\eqref{e:nest-identity} through \eqref{e:nest-identity-Schiffmann-bis}
are zero.
\end{remark}

\begin{example}\label{ex:nest-identity}
To illustrate Theorem~\ref{thm:nest-identity}, we write everything out
for the den defined by
\begin{equation}\label{e:ex-nest-identity-den}
h = 4, \quad p = \frac{1}{2}+\epsilon ,\quad \dd = (3,2,2,1,-1),\quad
\ee =(1,2,2,1,1).
\end{equation}
This den has sources at $(0,2)$, $(0,3)$ and sinks at $(4,0)$,
$(4,1)$.  Its nests are pairs $(\pi _{1},\pi _{2})$ of nested
generalized Dyck paths, with $\pi _{1}$ from $(0,2)$ to $(4,0)$ and
$\pi _{2}$ from $(0,3)$ to $(4,1)$, each path staying weakly below the
line of slope $-1/2$ connecting its endpoints.  For this den we have
$\gG =(2,2,2,2)$, since each nest consists of two paths from $x=0$ to
$x=4$.

The Catalanimal $H = H(R_{q},R_{t},R_{qt},\lambda )$ on the left hand
side of \eqref{e:nest-identity} has length $8$, with root sets and
weight displayed below.  Matrix position $(i,j)$ in the diagram
represents the root $\alpha _{ij}$.  The weight $\lambda $ is written
along the diagonal.
\begin{equation}\label{e:ex-nest-identity-cat}
\begin{tikzpicture}[scale=.4,baseline=1.6cm]
\draw [fill,black!25] (2,8) -- (2,6) -- (4,6) --(4,4) -- (6,4) --
(6,2) -- (8,2) -- (8,8) -- (2,8);
\draw [fill] (4.5,7.5) circle (.2);
\draw [fill] (5.5,7.5) circle (.2);
\draw [fill] (6.5,7.5) circle (.2);
\draw [fill] (7.5,7.5) circle (.2);
\draw [fill] (4.5,6.5) circle (.2);
\draw [fill] (5.5,6.5) circle (.2);
\draw [fill] (6.5,6.5) circle (.2);
\draw [fill] (7.5,6.5) circle (.2);
\draw [fill] (6.5,5.5) circle (.2);
\draw [fill] (7.5,5.5) circle (.2);
\draw [fill] (6.5,4.5) circle (.2);
\draw [fill] (7.5,4.5) circle (.2);
\node at (0.5,7.5) {\small $1$};
\node at (1.5,6.5) {\small $1$};
\node at (2.5,5.5) {\small $0$};
\node at (3.5,4.5) {\small $0$};
\node at (4.5,3.5) {\small $1$};
\node at (5.5,2.5) {\small $1$};
\node at (6.5,1.5) {\small $0$};
\node at (7.5,0.5) {\small $0$};
\draw (0,8) -- (8,8) -- (8,0);
\draw (1,8) -- (1,7) -- (8,7);
\draw (2,8) -- (2,6) -- (8,6);
\draw (3,8) -- (3,5) -- (8,5);
\draw (4,8) -- (4,4) -- (8,4);
\draw (5,8) -- (5,3) -- (8,3);
\draw (6,8) -- (6,2) -- (8,2);
\draw (7,8) -- (7,1) -- (8,1);
\end{tikzpicture}
\qquad \qquad 
\begin{tikzpicture}[scale=.4,baseline=1.6cm]
\draw [fill=black!25] (0,4.25) rectangle (1,5.25) node [right] at
(1.5,4.6) {$R_{q} = R_{t}$};
\draw [fill] (.5,3.25) circle (.2);
\draw  (0,2.75) rectangle (1,3.75) node [right] at (1.5,3.1) {$R_{q t}$};
\end{tikzpicture}
\end{equation}
One can verify by expanding the raising operator series that the
polynomial part of this Catalanimal is given by
\begin{multline}\label{e:ex-nest-identity-total}
H_{\pol }(\zz ) = (q^{3} t + q^{2} t^{2} + q\, t^{3})\, s_{31}(\zz )
+ (q^{4} + q^{3} t +2\, q^{2}t^{2} + q\, t^{3} + t^{4})\, s_{22}(\zz )
\\
+ (q^{3} + 2 \, q^{2}t+2\, q\, t^{2} + t^{3})\, s_{211}(\zz ) + (q^{2}
+ q\, t + t^{2})\, s_{1111}(\zz ).
\end{multline}

The following table displays each nest $\pi $ in the den along with
the corresponding term $t^{a(\pi )} q^{\dinv _{p}(\pi )} \Gcal
_{\nubold (\pi )}(\zz ;q^{-1})$ on the right hand side of
\eqref{e:nest-identity}, with $\Gcal _{\nubold (\pi )}(\zz ;q^{-1})$
expanded as a linear combination of Schur functions $s_{\lambda }(\zz
)$.  The reader can verify that \eqref{e:ex-nest-identity-total} is
the sum of these terms.
\begin{equation}\label{e:ex-nest-identity-nests}
\begin{aligned}
\begin{tikzpicture}[scale=.3,baseline=.3cm]
\draw [fill] (0,3) circle (.2); \draw [fill] (0,2) circle (.2);
\draw (4,1) circle (.2); \draw  (4,0) circle (.2);
\draw [dotted,thick,xshift=2pt] (0,3) -- (0,2) -- (2,2) -- (2,1) -- (4,1);
\draw [xshift=-2pt] (0,2) -- (0,1) -- (2,1) -- (2,0) -- (4,0);
\end{tikzpicture}
&\qquad t^{0}\,  q^{5}\, \big( q^{-1} s_{22}(\zz ) + q^{-2}
s_{211}(\zz )+q^{-3} s_{1111}(\zz ) \big)
\\[1ex]
\allowbreak 
\begin{tikzpicture}[scale=.3,baseline=.3cm]
\draw [fill] (0,3) circle (.2); \draw [fill] (0,2) circle (.2);
\draw (4,1) circle (.2); \draw  (4,0) circle (.2);
\draw [dotted,thick,xshift=2pt] (0,3) -- (0,2) -- (2,2) -- (2,1) -- (4,1);
\draw [xshift=-2pt] (0,2) -- (0,1) -- (1,1) -- (1,0) -- (4,0);
\end{tikzpicture}
&\qquad t^{1}\, q^{3}\, \big( s_{31}(\zz ) + s_{22}(\zz ) + 2\, q^{-1}
s_{211}(\zz ) + q^{-2} s_{1111}(\zz ) \big)
\\[1ex]
\allowbreak 
\begin{tikzpicture}[scale=.3,baseline=.3cm]
\draw [fill] (0,3) circle (.2); \draw [fill] (0,2) circle (.2);
\draw (4,1) circle (.2); \draw  (4,0) circle (.2);
\draw [dotted,thick,xshift=2pt] (0,3) -- (0,2) -- (2,2) -- (2,1) -- (4,1);
\draw [xshift=-2pt] (0,2) -- (0,0) --(4,0);
\end{tikzpicture}
&\qquad t^{2}\, q^{2}\, \big( s_{31}(\zz ) + s_{22}(\zz ) + q^{-1}
s_{211}(\zz ) \big)
\\[1ex]
\begin{tikzpicture}[scale=.3,baseline=.3cm]
\draw [fill] (0,3) circle (.2); \draw [fill] (0,2) circle (.2);
\draw (4,1) circle (.2); \draw  (4,0) circle (.2);
\draw [dotted,thick,xshift=2pt] (0,3) -- (0,2) -- (1,2) -- (1,1) -- (4,1);
\draw [xshift=-2pt] (0,2) -- (0,1) -- (1,1) -- (1,0) -- (4,0);
\end{tikzpicture}
&\qquad t^{2}\, q^{3}\, \big( q^{-1} s_{22}(\zz ) + q^{-2} s_{211}(\zz
) + q^{-3} s_{1111}(\zz ) \big)
\\[1ex]
\allowbreak 
\begin{tikzpicture}[scale=.3,baseline=.3cm]
\draw [fill] (0,3) circle (.2); \draw [fill] (0,2) circle (.2);
\draw (4,1) circle (.2); \draw  (4,0) circle (.2);
\draw [dotted,thick,xshift=2pt] (0,3) -- (0,2) -- (1,2) -- (1,1) -- (4,1);
\draw [xshift=-2pt] (0,2) -- (0,0) -- (4,0);
\end{tikzpicture}
&\qquad t^{3}\, q^{2}\, \big( q^{-1} s_{31}(\zz ) + q^{-1} s_{22}(\zz
) + q^{-2} s_{211}(\zz ) \big)
\\[1ex]
\allowbreak 
\begin{tikzpicture}[scale=.3,baseline=.3cm]
\draw [fill] (0,3) circle (.2); \draw [fill] (0,2) circle (.2);
\draw (4,1) circle (.2); \draw  (4,0) circle (.2);
\draw [dotted,thick,xshift=2pt] (0,3) -- (0,1) -- (4,1);
\draw [xshift=-2pt] (0,2) -- (0,0) -- (4,0);
\end{tikzpicture}
&\qquad t^{4}\,  q^{0}\, s_{22}(\zz )
\\
\end{aligned}
\end{equation}
\end{example}

\subsection{The single path case}
\label{ss:single-path}

As a further class of examples, we examine more closely the instances
of Theorem~\ref{thm:nest-identity} and
Corollary~\ref{cor:nest-identity-Schiffmann} for dens $(h,p,\dd ,\ee
)$ with just one source and one sink, that is,
\begin{equation}\label{e:single-path}
\dd - \ee = (1,0,\ldots,0,-1).
\end{equation}
We assume that $h>1$, as $h=1$ gives trivial dens with at most one
nest.  We also note that if $p<0$, then $d_{i} = e_{i}\leq e_{h}$ for
$0<i<h$.  This again implies that the den has at most one nest, so we
assume that $p\geq 0$.

Let $y+p\, x = s$ be the highest line of slope $-p$ that passes
through one of the heads ($= \text{feet}$) $(i,d_{i})$ for $0<i<h$.
Conditions (\ref{e:den-d-decr}--\ref{e:den-e-incr}) hold if and only
if all the heads for $0<i<h$ lie in the band $s-1<y+p\, x\leq s$, and
the source and sink at $x=0$ and $x = h$ are above the lower boundary
$y+p\, x = s-1$ of this band.  In other words, for $0<i<h$, the head
$(i,d_{i})$ is the highest lattice point below $y+p\, x = s$ on the
line $x = i$, and the source at $x = 0$ and sink at $x = h$ are weakly
above the highest lattice points on these respective lines.

Translating the picture vertically, we can assume that $e_{h} = 0$,
i.e., the sink is on the $x$-axis.  A nest $\pi $ in this den is then
a lattice path from $(0,d_{0})$ to $(h,0)$ that stays weakly below the
line $y+p\, x = s$ except possibly for an initial south run along the
$y$-axis and the final east step along the $x$-axis.

Let $r = s/p$ be the $x$-intercept of the line $y+p\, x = s$.  If
$r\geq h$, the sink $(h,0)$ is weakly below the line, $(h,1)$ is above
the line, and every lattice path weakly below the line that ends at
$(h,0)$ ends with an east step.  If $r<h$, the sink $(h,0)$ is above
the line, the point $(h-1,0)$ is below the line, so $h= \lfloor r
\rfloor + 1$, and deleting the final east step in each nest $\pi $
gives a path ending at $(\lfloor r \rfloor,0)$.  In either case, nests
$\pi $ correspond one-to-one with lattice paths from $(0,d_{0})$ to
$(\min (\lfloor r \rfloor,h), 0)$ that stay weakly below $y +p \, x =
s$ except for a possible south run along the $y$-axis.

In this picture, $a(\pi )$ is the area between the path $\pi $ and the
highest such path $\pi ^{0}$.  In \cite[Definition~5.4.1]{BHMPS-paths},
$\dinv _{p}(\pi )$ was defined to be the number of $p$-balanced hooks
whose arm and leg end on the path $\pi $; but this was also shown in
the proof of \cite[Proposition~5.4.4]{BHMPS-paths} to coincide with
$\dinv _{p}(\pi )$ as defined here.

Since there is only one path, we have $g_{i}=1$ for all $i$, so the
root sets for the Catalanimal $H$ in Theorem~\ref{thm:nest-identity}
are $R_{q} = R_{t} = R_{+}$ and $R_{qt} = [R_{+},R_{+}]$, and its
weight is given by $\lambda _{i} = d_{i-1}-e_{i}$, which is the number
of south steps on $x = i-1$ in the highest path $\pi ^{0}$ under the
line $y+p\, x = s$.  If $r<h$, so $h = \lfloor r \rfloor + 1$, then
$\lambda = \bb $, where $\bb = (b_{1},\ldots,b_{h})$ is as in
\cite[Theorem~5.5.1]{BHMPS-paths}, $H$ is equal to the function $\Hcal
_{\bb }(\zz )$ in \cite[Definition~3.7.1]{BHMPS-paths} (which is a
Catalanimal), and we have $\psi _{\hGamma }(H) = D_{\bb }$ in the
notation of \cite[\S 3.6]{BHMPS-paths} ($D_{\bb }$ is a {\em Negut
element} in $\Ecal ^{+}$).  In general, $h\leq \lfloor r \rfloor + 1$,
and $H$ is a possibly shorter Catalanimal than $\Hcal _{\bb }$, with
$\psi _{\hGamma } (H) = D_{\lambda }$, where $\lambda $ is obtained by
dropping some trailing zeroes from $\bb $.  By
\cite[Lemma~3.6.2]{BHMPS-paths}, the trailing zeroes do not matter and
we have $D_{\bb }\cdot 1 = D_{\lambda } \cdot 1$.

In the single path case, Corollary~\ref{cor:nest-identity-Schiffmann}
now reduces to the generalized shuffle theorem
\cite[Theorem~5.5.1]{BHMPS-paths} for paths under the line $y+p\, x=
s$, including the more general version in
\cite[Remark~5.5.2]{BHMPS-paths} for paths extended along the
$y$-axis.

\section{The Loehr-Warrington conjecture and its
\texorpdfstring{$(m,n)$}{(m,n)} extension}
\label{s:LW}

In this section we construct a den such that the associated
Catalanimal $H$ in Theorem~\ref{thm:nest-identity} represents the
element $s_{\mu }[-M X^{m,n}]$ in the Schiffmann algebra, for any
partition $\mu $ and coprime integers $m,n>0$.  In the case $n = 1$,
Corollary~\ref{cor:nest-identity-Schiffmann} then yields a
combinatorial formula for $\nabla ^{m} s_{\mu }$, which we will show
agrees with the one conjectured by Loehr and Warrington in
\cite{LoehWarr08}.

\subsection{LW dens}
\label{ss:LW-dens}

Given a partition $\mu $ of length $\ell (\mu )$, we define the
following for use in constructing its associated dens. 

\begin{align}
\label{e:eta(mu)}
h(\mu ) & = \mu _{1} + \ell (\mu ) - 1
              = \text{largest hook length in $\mu $},\\
\label{e:delta(mu)}
\delta _{i}(\mu ) & = \chi (\text{$\mu _{1}-1-i$ is the content of the
last box in some row of $\mu $}),\\ 
\label{e:epsilon(mu)}
\varepsilon _{i}(\mu ) & = \chi (i \geq \mu _{1})
\end{align}
for $0 \leq i \leq h(\mu )$, where $\chi (P) = 1$ if $P$ is true, $0$
if $P$ is false.  Note that the contents of boxes in $\mu $ range from
$\mu _{1} - h(\mu )$ to $\mu _{1}-1$.  Since $\mu $ is a partition,
for $i\leq h(\mu )$ we have $\varepsilon _{i}(\mu) = 1$ if and only if
$\mu _{1} - i$ is the content of the first box in some row of $\mu $.

\begin{defn}\label{def:LW-den}
The {\em LW den} associated to a partition $\mu $ and a pair of
coprime positive integers $m,n$ is the den $(h,p,\dd ,\ee )$ defined
as follows:
\begin{equation}\label{e:LW-h-p}
h = m\, h(\mu ); \quad p = n/m -\epsilon,
\end{equation}
where $\epsilon >0$ is small; 
\begin{equation}\label{e:LW-d-e}
\begin{gathered}
d_{i} = e_{i} = \lfloor n\,  h(\mu ) - i\, n/m \rfloor \  \text{for $i$
 not a multiple of $m$},\\
d_{jm} =  n\, h(\mu ) - n\, j + \delta _{j}(\mu ) - 1, \\
e_{jm} = n\, h(\mu ) - n\, j + \varepsilon _{j}(\mu ) - 1,
\end{gathered}
\end{equation}
for $0\leq i\leq m\, h(\mu )$, $0\leq j\leq h(\mu )$.
\end{defn}

We will see below (Proposition~\ref{prop:LW-den}) that the LW den is
indeed a den.

\begin{remark}\label{rem:LW-den-and-nests}
(i) The $p$ dependence in all constructions involving dens and nests
comes from comparisons between $p$ and finitely many rational numbers
$r$. By saying that $\epsilon > 0$ is `small' in \eqref{e:LW-h-p}, we
mean that such comparisons give $p<r$ if $n/m\leq r$, and $p>r$ if
$n/m > r$.

(ii) Since $d_{0} = n\, h(\mu ) = (n/m) h$, the line $y+(n/m)x=d_{0}$
meets the coordinate axes at $(h,0)$ and $(0, d_{0})$.  The quantity
$\lfloor n\, h(\mu ) - (n/m)i \rfloor$ is the $y$-coordinate of the
highest lattice point on the line $x = i$ weakly below
$y+(n/m)x=d_{0}$.

If $i$ is not a multiple of $m$, the highest lattice point is strictly
below $y+(n/m)x=d_{0}$, and $(i,d_{i}) = (i,e_{i})$ is this point.
For $i = jm$, the highest lattice point is on the line
$y+(n/m)x=d_{0}$.  In this case, $(i,d_{i})$ and $(i,e_{i})$ are each
either on the line or one unit below, depending on the values of
$\delta _{j}(\mu )$ and $\varepsilon _{j}(\mu )$.  It follows that the
sources and sinks in the LW den all lie on the bounding line $y+(n/m)x
= d_{0}$, and that paths in every nest in this den stay weakly below
the bounding line.

For $i = jm$, if $\delta _{j}(\mu ) = \varepsilon _{j}(\mu ) = 0$, the
head and foot $(i,d_{i}) = (i,e_{i})$ are both one unit below the
bounding line.  In this case there is no source or sink at $x = jm$,
and the condition that paths in a nest lie weakly below the heads
forbids the paths from touching the bounding line at $x = jm$.
The other heads impose no further conditions.
Every system of nested east end paths from the sources to the sinks,
which stay weakly below the bounding line and do not touch it at the
forbidden points, is therefore a nest in the LW den. 
\end{remark}

\begin{figure}
\[
\begin{array}{c@{\qquad \qquad }c}
\begin{tikzpicture}[scale=.5]
% axes
\draw [->,gray] (-.05,-.05) -- (9,-.050);
\draw [->,gray] (-.05,-.05) -- (-.05,9);
\foreach \x in {-.05,1,2,3,4,5,6,7,8}
  \draw [gray] (\x ,-.050) -- (\x,-.2);
\foreach \y in {-.05,1,2,3,4,5,6,7,8}
  \draw [gray](-.05,\y) -- (-.2,\y);
% bounding line
\draw [black!30] (0,8) -- (8,0);
% heads and feet
\draw [fill] (0,8) circle (.2);
\draw (0,7) circle (.2);
\draw [fill, fill opacity=.3] (1,6) circle (.2);
\draw [fill] (2,6) circle (.2);
\draw (2,5) circle (.2);
\draw [fill] (3,5) circle (.2);
\draw (3,4) circle (.2);
\draw [fill, fill opacity=.3] (4,4) circle (.2);
\draw [fill] (5,2) circle (.2);
\draw (5,3) circle (.2);
\draw [fill, fill opacity=.3] (6,2) circle (.2);
\draw (6,2) circle (.2);
\draw [fill] (7,0) circle (.2);
\draw (7,1) circle (.2);
\draw [fill] (8,-1) circle (.2);
\draw (8,0) circle (.2);
%forbidden point
\draw (1,7) node[very thick, cross=4pt] {};
% highest nest
\draw (3,5) -- (3,4) -- (4,4) -- (4,3) -- (5,3);
\draw [thick,dotted] (2.05,6) -- (2.05,4) -- (3,4) -- (3,3) -- (4,3) -- (4,2)
 --(6,2) -- (6,1) -- (7,1);
\draw  (.05,8) -- (.05,6) -- (1,6) -- (1,5) -- (1.95,5) -- (1.95,3) -- (3,3)
 --(3,2) -- (4,2) -- (4,1) -- (6,1) -- (6,.05) -- (8,.05);
\end{tikzpicture}
&
\begin{tikzpicture}[scale=.5]
% axes
\draw [->,gray] (-.05,-.05) -- (17,-.050);
\draw [->,gray] (-.05,-.05) -- (-.05,9);
\foreach \x in {-.05,1,2,3,4,5,6,7,8,9,10,11,12,13,14,15,16}
  \draw [gray] (\x ,-.050) -- (\x,-.2);
\foreach \y in {-.05,1,2,3,4,5,6,7,8}
  \draw [gray](-.05,\y) -- (-.2,\y);
% bounding line
\draw [black!30] (0,8) -- (16,0);
% heads and feet
\draw [fill] (0,8) circle (.2);
\draw (0,7) circle (.2);
\draw [fill, fill opacity=.3] (1,7) circle (.2);
\draw [fill, fill opacity=.3] (2,6) circle (.2);
\draw [fill, fill opacity=.3] (3,6) circle (.2);
\draw [fill] (4,6) circle (.2);
\draw (4,5) circle (.2);
\draw [fill, fill opacity=.3] (5,5) circle (.2);
\draw [fill] (6,5) circle (.2);
\draw (6,4) circle (.2);
\draw [fill, fill opacity=.3] (7,4) circle (.2);
\draw [fill, fill opacity=.3] (8,4) circle (.2);
\draw [fill, fill opacity=.3] (9,3) circle (.2);
\draw [fill] (10,2) circle (.2);
\draw (10,3) circle (.2);
\draw [fill, fill opacity=.3] (11,2) circle (.2);
\draw [fill, fill opacity=.3] (12,2) circle (.2);
\draw (12,2) circle (.2);
\draw [fill, fill opacity=.3] (13,1) circle (.2);
\draw [fill] (14,0) circle (.2);
\draw (14,1) circle (.2);
\draw [fill, fill opacity=.3] (15,0) circle (.2);
\draw [fill] (16,-1) circle (.2);
\draw (16,0) circle (.2);
%forbidden point
\draw (2,7) node[very thick, cross=4pt] {};
% highest nest
\draw (6,5) -- (6,4) -- (8,4) -- (8,3) -- (10,3);
\draw [thick,dotted] (4,6) -- (4,5) -- (5,5) -- (5,4) -- (6,4) -- (6,3)
 -- (8,3) -- (8,2) -- (12,2) -- (12,1) -- (14,1);
\draw  (.05,8) -- (.05,7) -- (1,7) -- (1,6) -- (3,6) -- (3,5)
 -- (4,5) -- (4,4) -- (5,4) -- (5,3) -- (6,3)
 --(6,2) -- (8,2) -- (8,1) -- (12,1) -- (12,.05) -- (16,.05);
\end{tikzpicture}
\\
\multicolumn{2}{c}{
\begin{tikzpicture}[scale=.5]
\draw [fill] (0,2) circle (.2);
\node [right] at (0,2) {$\ \text{head}$};
\draw (0,1) circle (.2);
\node [right] at (0,1) {$\ \text{foot}$};
\draw [fill, fill opacity=.3] (0,0) circle (.2);
\node [right] at (0,0) {$\ \text{both}$};
\draw (0,-1) node[very thick, cross=4pt] {};
\node [right] at (0,-1) {$\ \text{forbidden point}$};
\end{tikzpicture}
}
\end{array}
\]
\caption{\label{fig:LW-den-examples}%
LW dens for $\mu =(4,3,3,3,2)$ and $m = n =1$ (left), $m=2$, $n=1$
(right), with the highest nest $\pi ^{0}$ in each den.  Sources and
sinks are the points on the bounding line that are heads or feet but
not both.}
\end{figure}
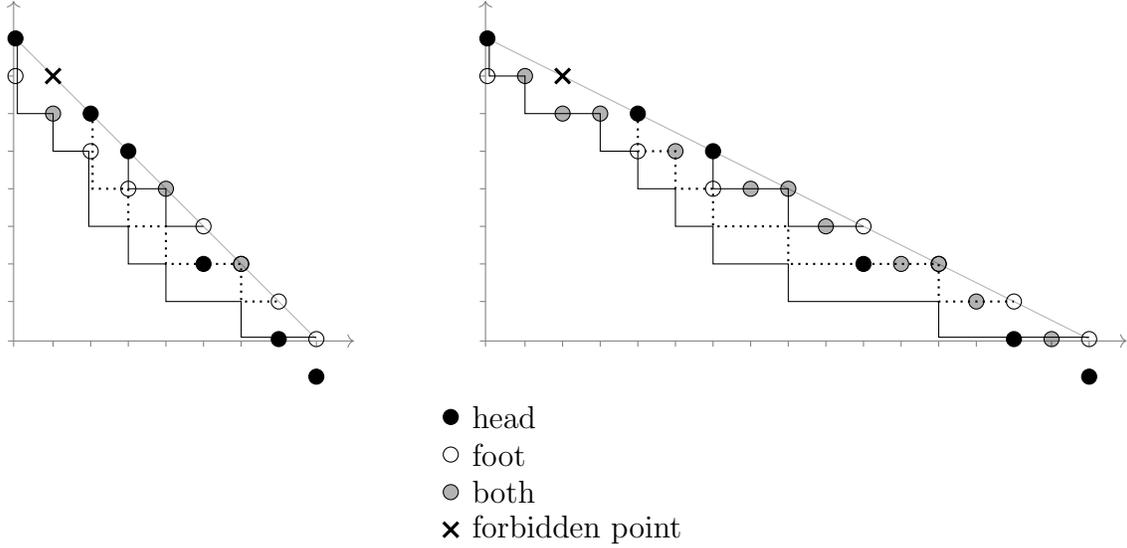

\begin{example}\label{ex:LW-den}
For the partition $\mu = (4,3,3,3,2)$ we have $h(\mu ) = 8$,
$\deltabold = (\delta _{0}(\mu ),\ldots,\delta _{8}(\mu )) =
(1,0,1,1,1,0,1,0,0)$, and $\epsbold =(\varepsilon _{0}(\mu
),\ldots,\varepsilon _{8}(\mu )) = (0,0,0,0,1,1,1,1,1)$.  In the
simplest LW den, for $m = n = 1$, we get $\dd $, $\ee $ by adding
$\deltabold - \onebold$, $\epsbold -\onebold$ to the vector
$(8,7,\ldots,1,0)$, giving
\begin{equation}\label{e:LW-den-1,1-example}
(h,p,\dd ,\ee ) = \big( 8, \, 1-\epsilon, \, (8,6,6,5,4,2,2,0,-1), \, 
(7,6,5,4,4,3,2,1,0) \big).
\end{equation}
For $(m,n) = (2,1)$ we interleave the above $\dd $ and $\ee $ with the
sequence $(7,6,5,4,3,2,1,0)$ of numbers $\lfloor 8-i/2 \rfloor$ for
$i$ odd, to get
\begin{multline}\label{e:LW-den-2,1-example}
(h,p,\dd ,\ee ) = \big( 16, \, \frac{1}{2}-\epsilon, \,
(8,7,6,6,6,5,5,4,4,3,2,2,2,1,0,0,-1),\\
(7,7,6,6,5,5,4,4,4,3,3,2,2,1,1,0,0) \big).
\end{multline}
These two dens are plotted in Figure~\ref{fig:LW-den-examples}.  In
each den, we also display the nest $\pi ^{0}$ such that $a(\pi ^{0}) =
0$, for later reference.
\end{example}

\begin{prop}\label{prop:LW-den}
The data $(h,p,\dd ,\ee )$ in Definition~\ref{def:LW-den} define a
valid den.
\end{prop}

\begin{proof}
For any $\mu $ we have $\delta _{0}(\mu )=1$, $\varepsilon _{0}(\mu )
= 0$, $\delta _{h(\mu )}(\mu ) = 0$ (because $\mu _{1}-1-h(\mu )$ is
less than the content of any box in $\mu $) and $\varepsilon _{h(\mu
)}(\mu ) = 1$.  We also have $\sum _{i=0}^{h(\mu )}\delta _{i}(\mu ) =
\sum _{i=0}^{h(\mu )}\varepsilon _{i}(\mu ) = \ell (\mu )$.  This
implies $d_{0}-e_{0} = 1 = e_{h} - d_{h}$ and $\sum _{i=0}^{h} d_{i} =
\sum _{i=0}^{h} e_{i}$,
verifying condition \eqref{e:den-d-e}.

All heads and feet are between or on the lines $y+(n/m)x = d_{0}$ and
$y+(n/m)x = d_{0}-1$, which implies $(d_{i}-d_{j}+1)/(j-i)\geq n/m
\geq (e_{i}-e_{j}-1)/(j-i)$ for all $i<j$.  Since all feet on the
lower line are left of all feet on the upper line, the second
inequality is strict.  This implies conditions
(\ref{e:den-d-decr}--\ref{e:den-e-incr}) for $p = n/m - \epsilon $.
\end{proof}

\subsection{An \texorpdfstring{$(m,n)$}{(m,n)} Loehr-Warrington formula}
\label{ss:mn-LW}

Our next theorem is a combinatorial formula for the symmetric function
$s_{\mu }[- X^{m,n}]\cdot 1$.  As we will see, this generalizes both
the $(km,kn)$-shuffle theorem \cite{BeGaSeXi16,Mellit16}, when $s_{\mu
} = e_{k}$, and the Loehr-Warrington conjecture \cite{LoehWarr08},
when $n=1$.

\begin{lemma}\label{lem:LW-H=schur-H} Let $(h,p,\dd ,\ee )$ be the LW
den for $\mu $ and $m,n$.

(i) With notation as in
Definition~\ref{def:Schur-Catalanimal}, the sequence $\gG $ in
\eqref{e:g-vector} is the same as $\gamma ((\mu ^{\circ })^{m})$, or
equivalently $\gamma (\mu ^{\circ })$ with each entry repeated $m$
times.

(ii) The Catalanimal $H$ in Theorem~\ref{thm:nest-identity} is the
same as the opposite Schur Catalanimal $H^{m,n}_{(\mu ^{\circ })^{m}}$
in Definition~\ref{def:Schur-Catalanimal} and
Theorem~\ref{thm:Schur-Catalanimal}.
\end{lemma}

\begin{proof}
The sources in the LW den have $x$-coordinate $i\, m$ for $\delta
_{i}(\mu ) = 1$, $\varepsilon _{i}(\mu ) = 0$, and the sinks have
$x$-coordinate $i\, m$ for $\delta _{i}(\mu ) = 0$, $\varepsilon
_{i}(\mu ) = 1$.  By Remark~\ref{rem:nest-basics} (iii), the
associated sequence $\gG $ is therefore obtained by repeating each
entry of the sequence $\gamma = (\gamma _{1},\ldots,\gamma _{h(\mu
)})$ $m$ times, where $\gamma _{i+1}-\gamma _{i} = \delta _{i}(\mu
)-\varepsilon _{i}(\mu )$ for $i=0,\ldots,h(\mu )-1$, and we set
$\gamma _{0} = 1$.  To prove (i), we need to show that $\gamma =\gamma
(\mu ^{\circ })$, which is the reverse of $\gamma (\mu )$.

Consider the decomposition of $\mu $ into hooks with corner on the
main diagonal, as shown here for the partition $\mu = (4,3,3,3,2)$ in
Example~\ref{ex:LW-den}.
\begin{equation}\label{e:main-hooks}
\begin{array}{c}
\begin{tikzpicture}[scale=.25]
\draw (0,0) -- (4,0);
\draw (0,1) -- (4,1);
\draw (0,2) -- (3,2);
\draw (0,3) -- (3,3);
\draw (0,4) -- (3,4);
\draw (0,5) -- (2,5);
\draw (0,0) -- (0,5);
\draw (1,0) -- (1,5);
\draw (2,0) -- (2,5);
\draw (3,0) -- (3,4);
\draw (4,0) -- (4,1);
\draw [very thick, gray] (.5,4.5) -- (.5,.5) -- (3.5,.5);
\draw [very thick, gray] (1.5,4.5) -- (1.5,1.5) -- (2.5,1.5);
\draw [very thick, gray] (2.5,3.5) -- (2.5,2.5);
\end{tikzpicture}
\end{array}
\end{equation}
From the definitions we see that $\delta _{i}(\mu ) = 1$ and
$\varepsilon _{i}(\mu ) = 0$ if and only if $\mu _{1} - i - 1$ is the
content of the rightmost box in the arm of a main diagonal hook.  This
is also the condition on $i\in [0, h(\mu )]$ to have $\gamma
_{i+1}(\mu ^{\circ }) - \gamma _{i}(\mu ^{\circ }) = 1$, if we define
$\gamma _{i}(\mu ^{\circ }) = 0$ for $i<1$ or $i>h(\mu )$.

Noting that a diagonal of content $c\leq 0$ contains the top of the
leg in a main diagonal hook if and only if the diagonal of content
$c-1$ does not contain the box at the end of a row of $\mu $, we see
that $\delta _{i}(\mu ) = 0$ and $\varepsilon _{i}(\mu ) = 1$ if and
only if $\mu _{1} - i$ is the content of the highest box in a leg.
This is the condition on $i$ to have $\gamma _{i+1}(\mu ^{\circ }) -
\gamma _{i}(\mu ^{\circ }) = -1$.

When neither of these two conditions hold, we have $\gamma _{i+1}(\mu
^{\circ } ) = \gamma _{i}(\mu ^{\circ })$.  Hence, $\gamma _{i+1}(\mu
^{\circ })-\gamma _{i}(\mu ^{\circ }) = \delta _{i}(\mu )-\varepsilon
_{i}(\mu )$ for all $i$, and therefore $\gamma = \gamma (\mu ^{\circ
})$, proving (i).

Part (i) implies that the Catalanimals $H$ and $H^{m,n}_{(\mu ^{\circ
})^{m}}$ have the same root sets.  The weight $\lambda $ for
$H^{m,n}_{(\mu ^{\circ })^{m}}$, given by \eqref{e:weight-lambda}, is
constant on blocks of lengths $g_{1},\ldots,g_{h}$.  We need to verify
that its value on the $i$-th block is $d_{i-1} - e_{i}$.

We can write \eqref{e:LW-d-e} in the form
\begin{align}\label{e:LW-d}
d_{i} & =  n\, h(\mu )- \lceil  i\, n/m \rceil +
 \chi (m|i)(\delta _{i/m}(\mu ) - 1)\\
\label{e:LW-e}
e_{i} & = n\, h(\mu )- \lceil i\, n/m \rceil +
 \chi (m|i)(\varepsilon _{i/m}(\mu ) - 1).
\end{align}
 Combining
(\ref{e:LW-d}--\ref{e:LW-e}) with the definition of $\bb (m,n)$ in
\eqref{e:bbold} gives
\begin{equation}
\label{e:d-minus-e} d_{i-1} - e_{i} = \bb (m, n)_{\modm _{m}(i)} +
\chi (m|i-1)(\delta _{(i-1)/m}(\mu )-1) + \chi (m|i)(1-\varepsilon
_{i/m}(\mu )).
\end{equation}
To compare this with \eqref{e:weight-lambda}, note that the content
$c$ on the $i$-th diagonal from northwest to southeast in $(\mu
^{\circ })^{m}$ has $c\equiv i \pmod{m}$ by construction, so the term
$\bb (m, n)_{\modm _{m}(i)}$ here agrees with $\bb (m, n)_{\modm
_{m}(c)}$ in \eqref{e:weight-lambda}.

The $i$-th diagonal in $(\mu ^{\circ })^{m}$ always contains
the first box in a row of $(\mu ^{\circ })^{m}$ if $i\not \equiv 1
\pmod{m}$.  If $i = j\, m+1$, the $i$-th diagonal contains the first
box in a row of $(\mu ^{\circ })^{m}$ if and only if the $(j+1)$-st
diagonal in $\mu ^{\circ }$ contains the first box in a row of $\mu
^{\circ }$, that is, if and only if $\delta _{j}(\mu ) = 1$.  Hence,
\begin{multline}\label{e:mu-op-row-starts}
\chi (\text{$i$-th diagonal in
$(\mu ^{\circ })^{m}$ contains the first box in a row})\\
= 1 + \chi (m|i-1)(\delta _{(i-1)/m}(\mu )-1).
\end{multline}
Similarly, the $i$-th diagonal in $(\mu ^{\circ })^{m}$ always
contains the last box in a row of $(\mu ^{\circ })^{m}$ if $i \not
\equiv 0 \pmod{m}$. If $i = j\, m$, the $i$-th diagonal contains the
last box in a row of $(\mu ^{\circ })^{m}$ if and only if the $j$-th
diagonal in $\mu ^{\circ }$ contains the last box in a row of $\mu
^{\circ }$, that is, if and only if $j\geq \mu _{1}$.  Hence,
\begin{multline}\label{e:mu-op-row-ends}
\chi (\text{$i$-th diagonal in
$(\mu ^{\circ })^{m}$ contains the last box in a row})\\
= 1 + \chi (m|i)(\varepsilon _{i/m}(\mu )-1).
\end{multline}
Using \eqref{e:mu-op-row-starts} and \eqref{e:mu-op-row-ends} we see
that \eqref{e:d-minus-e} agrees with \eqref{e:weight-lambda}.
\end{proof}

Our $(m,n)$ Loehr-Warrington formula is now given by the following
theorem, which follows from Proposition~\ref{prop:nabla-m-f},
Theorem~\ref{thm:Schur-Catalanimal},
Corollary~\ref{cor:nest-identity-Schiffmann}, and
Lemma~\ref{lem:LW-H=schur-H}.

\begin{thm}\label{thm:mn-LW} 
For every partition $\mu $ and coprime positive integers $m,n$, we
have the identity
\begin{equation}\label{e:mn-LW}
s_{\mu }[- M X^{m,n}]\cdot 1 = (-1)^{p(\mu )}(q\, t)^{p(\mu ) + m \,
n'(\gamma (\mu ))} \sum _{\pi } t^{a(\pi )} q^{\dinv _{p}(\pi )}
\omega\, \Gcal _{\nubold (\pi )}(X;q^{-1}),
\end{equation}
where $p(\mu )$ is the magic number (not to be confused with $p =
n/m-\epsilon $), $n'(\gamma (\mu ))$ is given by \eqref{e:n'(gamma)},
and the sum is over nests $\pi $ in the LW den for $\mu $ and $m,n$.

For $n = 1$, the expression in \eqref{e:mn-LW} is equal to $\nabla
^{m} s_{\mu }(X)$.
\end{thm}

\begin{remark}\label{rem:(km,kn)-shuffle}
In the case $\mu = (1^{k})$, giving $s_{\mu } = e_{k}$, we have $p(\mu
) = n'(\gamma (\mu )) = 0$.  As explained in \cite[\S 6]{BHMPS-paths},
Theorem~\ref{thm:mn-LW} is then equivalent to the $(km,kn)$-shuffle
theorem conjectured by Bergeron et\ al.~\cite{BeGaSeXi16} and proven
by Mellit \cite{Mellit16}.  This case is also included in the path
case discussed in \S \ref{ss:single-path}, above.
\end{remark}

\subsection{Comparison with the original Loehr-Warrington formula}
\label{ss:original-LW}

Loehr and Warrington \cite[Conjecture~2.4]{LoehWarr08} conjectured a
combinatorial formula for $\nabla ^{m} s_{\mu }$, which in their
notation reads
\begin{equation}\label{e:LW-conjecture}
\nabla ^{m }s_{\mu }(X) =
\sgn (\mu ) \sum _{(G,R)\in \LNDP ^{m}_{\mu }} t^{\area (G,R)}
q^{\dinv (G,R)} \xx _{R}.
\end{equation}
We prove their conjecture (after correcting a mistake in the
definition of $\dinv (G,R)$---see below) by verifying that this
formula agrees with Theorem~\ref{thm:mn-LW} for $n=1$, when we express
the latter in terms of labeled nests as in
\eqref{e:nest-identity-Schiffmann-bis}, namely
\begin{equation}\label{e:m-1-LW}
\nabla ^{m }s_{\mu }(X) = (-1)^{p(\mu )}(q\, t)^{p(\mu ) + m \,
n'(\gamma (\mu ))} \sum _{\pi ,N} t^{a(\pi )} q^{\dinv _{p}(\pi ) -
\inv (N)} \xx ^{N}.
\end{equation}
We now work out some of the details needed to see that
\eqref{e:LW-conjecture} and \eqref{e:m-1-LW} agree.

The notation $(G,R)\in \LNDP ^{m}_{\mu }$ stands for a system of
labeled nested $m$-Dyck paths which, after reflecting about a
horizontal line, corresponds to a pair ($\pi , N$), where $\pi $ is a
nest in the LW den for $\mu $ and $m,n$ with $n=1$, and $N$ is a negative
labeling of $\pi $ as in \eqref{e:pi-labeling-neg}.

In more detail, the proof of Lemma~\ref{lem:LW-H=schur-H} shows that
the sources and sinks in the LW den are located on the bounding line
at positions $x = jm$, where $\mu _{1}-1-j$ is the content of the box
at either end of a main diagonal hook in $\mu $.  Loehr and Warrington
start and end their Dyck paths at corresponding positions constructed
from a dissection of $\mu $ into border strips
(\cite[Fig.~1]{LoehWarr08}).  To see that these agree, observe that
the heights, widths and content ranges of the border strips for any
$\mu $ match those of the main diagonal hooks, as in this example:
\begin{equation}\label{e:main-hooks=border-strips}
\begin{array}{c@{\qquad \qquad }c}
\begin{tikzpicture}[scale=.25]
\draw (0,0) -- (4,0);
\draw (0,1) -- (4,1);
\draw (0,2) -- (3,2);
\draw (0,3) -- (3,3);
\draw (0,4) -- (3,4);
\draw (0,5) -- (2,5);
\draw (0,0) -- (0,5);
\draw (1,0) -- (1,5);
\draw (2,0) -- (2,5);
\draw (3,0) -- (3,4);
\draw (4,0) -- (4,1);
\draw [very thick, gray] (.5,4.5) -- (.5,.5) -- (3.5,.5);
\draw [very thick, gray] (1.5,4.5) -- (1.5,1.5) -- (2.5,1.5);
\draw [very thick, gray] (2.5,3.5) -- (2.5,2.5);
\end{tikzpicture}
&
\begin{tikzpicture}[scale=.25]
\draw (0,0) -- (4,0);
\draw (0,1) -- (4,1);
\draw (0,2) -- (3,2);
\draw (0,3) -- (3,3);
\draw (0,4) -- (3,4);
\draw (0,5) -- (2,5);
\draw (0,0) -- (0,5);
\draw (1,0) -- (1,5);
\draw (2,0) -- (2,5);
\draw (3,0) -- (3,4);
\draw (4,0) -- (4,1);
\draw [very thick, gray] (.5,4.5) -- (1.5,4.5) -- (1.5,3.5) -- (2.5,3.5)
 -- (2.5,.5) -- (3.5,.5);
\draw [very thick, gray] (.5,3.5) -- (.5,2.5) -- (1.5,2.5) -- (1.5,.5);
\draw [very thick, gray] (.5,1.5) -- (.5,.5);
\end{tikzpicture}
\end{array}
\end{equation}
This also shows that the statistic $\spin (\mu ^{*}) = \adj (\mu )$
in \cite[\S 2.2]{LoehWarr08}, defined as the sum $\sum
_{B}({\text{width(B)}-1})$ over border strips $B$ in $\mu $, is equal to
the magic number $p(\mu )$, since the latter is
the number of boxes strictly below the main diagonal, and thus the sum of the
arm lengths
of the main diagonal hooks.

Loehr and Warrington introduce additional Dyck paths of length zero
which serve to forbid the other paths from touching the bounding line
at certain points.  These correspond to the points where nests in the
LW den cannot touch the bounding line, as in
Remark~\ref{rem:LW-den-and-nests} (ii).  The rest of
Remark~\ref{rem:LW-den-and-nests} (ii) then shows that nests $\pi $ in
the LW den correspond to systems of $m$-Dyck paths $G$ for $(G,R)\in
\LNDP ^{m}_{\mu }$.

Deciphering the notation in \cite{LoehWarr08} further, the indices
$(a,u)$ of entries $g_{a}^{(u)}$, $r_{a}^{(u)}$ in $(G,R)$ correspond
to pairs $(S,i)$, where $S$ is a south step in
a path $\pi _{i}$ in the nest $\pi $.
When $(a,u)$ corresponds to $(S,i)$, the
index $u$ is a strictly increasing function of $i$.
The integer $g_{a}^{(u)}$ is equal to $m$ times the vertical
distance from the north endpoint of $S$ to the bounding line $y+(1/m)x
= d_{0}$ for the LW den.
The labels $r_{a}^{(u)}$ are subject to the same
conditions as our labels $N(S,i)$.

To finish reconciling \eqref{e:LW-conjecture} with \eqref{e:m-1-LW},
we need to show the following.

\begin{lemma}\label{lem:LW-stats}
For $(G,R)$ corresponding to $(\pi ,N)$, the combinatorial statistics
in \cite{LoehWarr08} are related to ours by
\begin{align}
\label{e:LW-sgn}
\sgn (\mu ) & \defeq (-1)^{\spin (\mu ^{*})} = (-1)^{p(\mu )}\\
\label{e:LW-area}
\area (G,R) & = p(\mu ) + m\, n'(\gamma (\mu )) + a(\pi )\\
\label{e:LW-dinv}
\dinv (G,R) & = p(\mu ) + m\, n'(\gamma (\mu )) + \dinv _{p}(\pi ) -
\inv (N),
\end{align}
where $p = 1/m-\epsilon $.
\end{lemma}

\begin{proof}
For \eqref{e:LW-sgn}, we already observed that $\spin (\mu ^{*}) =
p(\mu )$.

For \eqref{e:LW-area}, the definition of $\area (G,R)$ corresponds to
the sum of the areas $|\rho /\pi _{i}|$, where $\rho $ is the highest
lattice path from $(0,d_{0}) = (0, h(\mu ))$ to $(h,0) = (m\, h(\mu
),0)$ that stays weakly below the bounding line $y + (1/m)x = d_{0}$.
To verify \eqref{e:LW-area} it therefore suffices to show that for the
nest $\pi ^{0}$ with $a(\pi ^{0}) = 0$, we have $\sum _{i} |\rho /\pi
^{0}_{i}| = p(\mu ) + m n'(\gamma (\mu ))$.

Now $\sum _{i} |\rho /\pi ^{0}_{i}| = \sum _{E} a_{E}$, where the sum
is over east steps $E$ in $\pi ^{0}$, and $a_{E}$ is the vertical
distance between $E$ and the east step weakly above it in $\rho $.
Recall that $g_{k}$ as defined in \eqref{e:g-vector} is the
number of east steps from $x = k-1$ to $x = k$ in any nest.  For $k =
jm$ with $1\leq j< \mu _{1}$,
the point at
$x = k$
on the bounding line for the LW den is either a source or a forbidden
point, and thus is not the right endpoint of an east step in $\pi
^{0}$.
For these
values of $k$, the numbers $a_{E}$ for east steps $E$ from $x = k-1$
to $x = k$ are $1,\ldots,g_{k}$.
For other values of $k$ the $a_{E}$ are $0,1,\ldots,g_{k}-1$.
It may be instructive to verify this with the examples
in Figure~\ref{fig:LW-den-examples}.

It follows that
\begin{equation}\label{e:aE-sum}
\sum _{i} |\rho /\pi ^{0}_{i}| = \sum _{E} a_{E} = \sum _{k=1}^{m
h(\mu )}\binom{g_{k}}{2} + \sum _{j=1}^{\mu _{1} - 1} g_{jm}.
\end{equation}
The first sum on the right is $m\, n'(\gamma (\mu ))$
by Lemma~\ref{lem:LW-H=schur-H}(i).
Since $g_{jm}$ is the number of boxes on the diagonal of content $\mu
_{1} - j$ in $\mu $, the second sum is $p(\mu )$.

For \eqref{e:LW-dinv}, the statistic $\dinv (G,R)$ is defined by
\cite[(11)]{LoehWarr08}, except that
the expression $\chi (a\leq b)$ in the middle sum there should read
$\chi ((a < b) \vee ((a=b) \wedge (u<v)))$, as in the last sum.  After
correcting this mistake and exchanging indices $(a,u)$ and $(b,v)$ in
the first two sums, we can rewrite \cite[(11)]{LoehWarr08} in the form
\begin{multline}\label{e:LW-dinv-def}
\dinv (G,R) =\\
\begin{aligned}
\adj (\mu )
& + \sum \chi (0 < g_{b}^{(v)} - g_{a}^{(u)} \leq m)\, \chi (a\geq b)\, 
 \chi (r_{a}^{(u)} < r_{b}^{(v)})\\
& + \sum \chi (0 \leq  g_{b}^{(v)} - g_{a}^{(u)} < m)\, \chi (a < b)\, 
 \chi (r_{a}^{(u)} < r_{b}^{(v)})\\
& + \sum \chi (0 \leq  g_{b}^{(v)} - g_{a}^{(u)} < m)\,  \chi (a = b
 \wedge u < v)\, \chi (r_{a}^{(u)} < r_{b}^{(v)})\\
& + \sum |(g_{b}^{(v)} - g_{a}^{(u)} + [0,m-1])\cap [1,m-1]|\, \chi
 ((a> b) \vee (a = b \wedge u > v)),
\end{aligned}
\end{multline}
where the sums are over all pairs of valid indices $(a,u)$, $(b,v)$.
Let $(S,i)$, $(S',j)$ be the south steps on paths in $\pi $
corresponding to $(a,u)$ and $(b,v)$.  The condition 
\begin{equation}\label{e:attacking-Gs}
\big( (0 < g_{b}^{(v)} - g_{a}^{(u)} \leq m )\wedge (a \geq b) \big)
\vee \big( (0 \leq g_{b}^{(v)} - g_{a}^{(u)} < m )\wedge (a < b) \big)
\end{equation}
holds if and only if $(S,i)$ and $(S',j)$ form an attacking pair in
$\Sbold (\pi )$, as defined in \S \ref{ss:nest-LLT}, for $p = 1/m -
\epsilon $. We leave it as an exercise for the reader to verify this,
with the hint that if $g_{b}^{(v)} = g_{a}^{(u)}$, then $S$ is
strictly to the left of $S'$ if and only if $a<b$, while if
$g_{b}^{(v)} - g_{a}^{(u)} = m$, then $S$ is strictly to the right of
$S'$ if and only if $a\geq b$.  The first two sums in
\eqref{e:LW-dinv-def} therefore count attacking pairs that do not
contribute to $\inv (N)$, that is, they add up to
\begin{equation}\label{e:non-attacking}
\Acal (\pi ) - \inv (N),
\end{equation}
where $\Acal (\pi )$ is the number of attacking pairs in $\Sbold (\pi
)$, or equivalently in $\nubold (\pi )$.

Turning to the third sum in \eqref{e:LW-dinv-def}, If $a = b$ and
$g_{a}^{(u)}\leq g_{b}^{(v)}$, then $S$ and $S'$ have the same
$y$-coordinates, with $S'$ weakly to the left of $S$.  If $S'$ is
strictly to the left of $S$, then nesting implies $v<u$.  Hence, if
$u<v$, then $S = S'$ is a shared south step on paths $\pi _{i}$ and
$\pi _{j}$ with $i<j$.  In this case, the conditions on the labeling
imply $r_{a}^{(u)}<r_{b}^{(v)}$.  The third sum in
\eqref{e:LW-dinv-def} therefore reduces to the number of pairs
$\{(S,i), (S,j)\}$ of shared south steps in $\pi $.  We denote this
number by $\sshare (\pi )$.

Now we consider the last sum in \eqref{e:LW-dinv-def}.  The terms with
$a = b$, $g_{a}^{(u)} = g_{b}^{(v)}$ and $u>v$ contribute $m-1$ for
each pair of shared south steps, giving $(m-1)\sshare (\pi )$.

The remaining terms are zero if $|g_{b}^{(v)} - g_{a}^{(u)}|\geq m$.
Otherwise, they correspond to pairs $(S,i)$, $(S',j)$ with $S\not =S'$
and some line of slope $-1/m$ passing through the interiors of both
$S$ and $S'$.  If $S'$ is to the right of $S$, this implies $b\geq a$,
with $u<v$ if $a = b$, by nesting.  Hence, $S'$ is to the left of $S$.
To describe the contribution from such a term geometrically, let
$B_{S}$ be the region bounded on the right by $S$ and above and below
by lines of slope $-1/m$ through the endpoints of $S$, with open
boundaries above and on the right, and a closed boundary below.  Then
$S'$ has an endpoint in $B_{S}$, and the contribution from the
corresponding term is given by
\begin{equation}\label{e:BS-term}
r-1,\ \text{for}\quad  
\begin{tikzpicture}[scale=.7,baseline=.7cm]
\fill [fill=black!10] (0,1) -- (3,0) -- (3,1) -- (0,2);
\draw (0,1) -- (3,0);
\draw [dotted] (3,1) -- (0,2);
\draw [very thick, dashed] (3,0) -- (3,1);
\node [right] at (3,.5) {$S$};
\draw [fill] (1.5,1) circle (.1) -- (1.5,2);
\node [right] at (1.5,2) {$S'$};
\end{tikzpicture}
\qquad \qquad  r,\ \text{for}\quad 
\begin{tikzpicture}[scale=.7,baseline=.55cm]
\fill [fill=black!10] (0,1) -- (3,0) -- (3,1) -- (0,2);
\draw (0,1) -- (3,0);draw [dotted] (3,1) -- (0,2);
\draw [very thick, dashed] (3,0) -- (3,1);
\node [right] at (3,.5) {$S$};
\draw [fill] (1.5,1) circle (.1) -- (1.5,0);
\node [left] at (1.5,0) {$S'$};
\end{tikzpicture}\ ,
\end{equation}
where 
$r = m - |g_{b}^{(v)} - g_{a}^{(u)}|$
is the integer such that the line segment $S'\cap B_{S}$ has
length $r/m$.
Note that we require the upper endpoint of $S'$ in the second
picture to be in the interior of $B_{S}$ (although if it were on the
boundary, it would contribute zero anyway).

Let $\delta (\pi )$ denote the sum of the contributions
in \eqref{e:BS-term} for all pairs $(S,i)$, $(S',j)$ with $S$ and $S'$
positioned as shown.  Then the last sum in \eqref{e:LW-dinv-def} is
equal to $\delta (\pi ) + (m-1)\sshare (\pi )$.

Combining $\adj (\mu ) = p(\mu )$ with the above, we obtain
\begin{equation}\label{e:LW-dinv-first-formula}
\dinv (G,R) =p(\mu ) + m \sshare (\pi ) + \delta (\pi ) + \Acal (\pi ) -
\inv (N).
\end{equation}
The following lemma now completes the verification of \eqref{e:LW-dinv}.
\end{proof}

\begin{lemma}\label{lem:dinv-p=1/m-eps}
For $p = 1/m - \epsilon $ and $\pi $ a nest in the LW den for $\mu $ and
$m, n$ with $n=1$,
we have
\begin{equation}\label{e:dinv-p=1/m-eps}
\dinv _{p}(\pi ) = m\, \sshare (\pi ) - m\, n'(\gamma (\mu )) + \Acal (\pi
) + \delta (\pi ).
\end{equation}
\end{lemma}

\begin{proof}
We evaluate $\dinv _{p}(\pi ) - \Acal (\pi ) - \delta (\pi )$.

For $p = 1/m - \epsilon $, each unordered pair $\{ (S,i), (S',j) \}$
with $S,S'$ positioned as in \eqref{e:BS-term} forms an attacking 
pair in $\Sbold (\pi )$ when ordered with the smaller of $\hat{c}(S)$,
$\hat{c}(S')$ first.  Every attacking pair has this form, so $\Acal
(\pi )$ is the number of such unordered pairs.

Given a pair $(S,i)$ in $\Sbold (\pi )$ and a path $\pi _{j}$,
consider the connected components of $\pi _{j}\cap B_{S}$. These are
of four possible types,
depicted in Figure~\ref{fig:BS-components},
depending on whether they enter and exit $B_{S}$ along the upper or
the lower boundaries.

\begin{figure}
\[
\begin{tikzpicture}[scale=.7]
\fill [fill=black!10] (-1,1) -- (15.5,-4.5) -- (15.5,-3.5) -- (-1,2);
\draw (15.5,-4.5) -- (-1,1);
\draw [dotted] (15.5,-3.5) -- (-1,2);
\draw [very thick, dashed] (15.5,-4.5) -- (15.5,-3.5);
\node [right] at (15.5,-4) {$S$};
\draw (0,2) -- (0,1) -- (2,1); \draw [fill] (0,1) circle (.1);
\draw [fill] (1,1) circle (.1);
\draw (3,1) -- (3,0) -- (4,0) -- (4,-1);
\draw [fill] (3,0) circle (.1);
\draw [fill] (4,0) circle (.1);
\draw (5.5,-1.5) -- (8.5,-1.5) -- (8.5,-2.5);
\draw [fill] (6.5,-1.5) circle (.1);
\draw [fill] (7.5,-1.5) circle (.1);
\draw [fill] (8.5,-1.5) circle (.1);
\draw (10,-3) -- (14,-3); \draw [fill] (11,-3) circle (.1);
\draw [fill] (12,-3) circle (.1); \draw [fill] (13,-3) circle (.1);
\end{tikzpicture}
\]
\caption{\label{fig:BS-components}%
The possibilities for components of the intersection of a path $\pi
_{j}$ with $B_{S}$, illustrated with $m = 3$, $p = 1/3-\epsilon $.}
\end{figure}
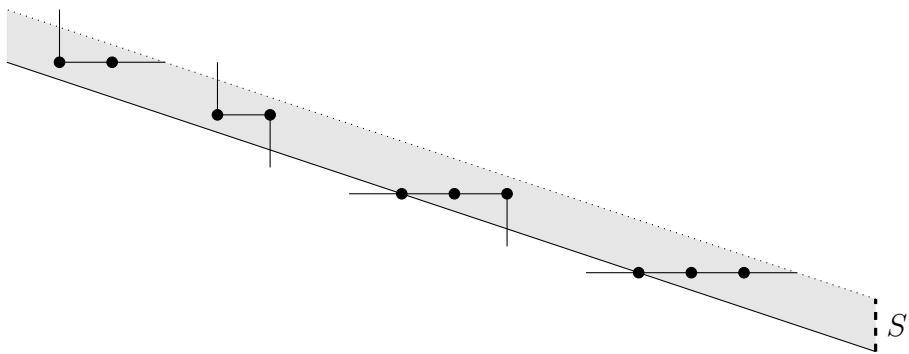

For $p = 1/m-\epsilon $, $\dinv _{p}(\pi )$ counts pairs $(S,i)$,
$(P,j)$, where $P$ is a point on $\pi _{j}$ that lies in $B_{S}$.
Using this, the description of $\Acal (\pi )$ above, and the
definition of $\delta (\pi )$, one can check that each component of
$\pi _{j}\cap B_{S}$ contributes $-m$ to $\dinv _{p}(\pi ) - \Acal
(\pi ) - \delta (\pi )$ for components that cross $B_{S}$ from top to
bottom, $m$ for components that cross from bottom to top, and zero for
components of the other two types.

Since all sources are on the bounding line, the leftmost component (if
any) of $\pi _{j}\cap B_{S}$ enters from above.
From left to right, the components that cross $B_{S}$ alternate
between the second and fourth types shown in
Figure~\ref{fig:BS-components}, possibly with components of the other
types in between.

If $i\leq j$, so $\pi _{j}$ is equal to or nested above $\pi _{i}$, or
if $i>j$ and the paths $\pi _{i}$ and $\pi _{j}$ share the south step $S$,
then $\pi _{j}\cap B_{S}$ is either empty, or its last component exits
$B_{S}$ on the upper boundary.  In this case, crossings from top to
bottom cancel those from bottom to top, giving a net contribution of zero
to $\dinv _{p}(\pi ) - \Acal (\pi ) - \delta (\pi )$.

Otherwise, if $i>j$ and $S$ is not on $\pi _{j}$, there is one more
crossing from top to bottom than from bottom to top, for a net
contribution of $-m$.  This shows that $\dinv _{p}(\pi ) - \Acal (\pi
) - \delta (\pi )$ is equal to $-m$ times the number of tuples
$(S,i,j)$ with $S$ on $\pi _{i}$ and $i>j$, plus $m \sshare (\pi )$.  The
number of such tuples $(S,i,j)$ is the sum over all $i>j$ of the
number of south steps on $\pi _{i}$.

For the LW den, the number of south steps on $\pi _{i}$ is the length
of the $i$-th main diagonal hook in $\mu $, from southwest to northeast,
and there are $i-1$ indices $j<i$.  If we write $i-1$ in each box on
the $i$-th main diagonal hook, the sum of these numbers is therefore
the number of tuples $(S,i,j)$.  But these numbers sum to
$\binom{\gamma }{2}$ on a content diagonal of length $\gamma $, so the
sum of them all is $n'(\gamma (\mu ))$.  This gives
\begin{equation}\label{e:dinv-A-delta}
\dinv _{p}(\pi ) - \Acal (\pi ) - \delta (\pi ) = m\,\sshare (\pi )- m\,
n'(\gamma (\mu )),
\end{equation}
as desired.
\end{proof}

\section{LLT series and semi-symmetric Hall-Littlewood polynomials}
\label{s:LLT}

Using a strategy similar to that in \cite{BHMPS-paths}, we will prove
Theorem~\ref{thm:nest-identity} by taking the polynomial part of an
infinite series identity between the full Catalanimal on the left hand
side and a sum of LLT series on the right.  LLT series associated to a
reductive group $G$ and Levi subgroup $L\subseteq G$ were defined in
\cite{GrojHaim07}; for $G = \GL _{l}$, they are series versions of the
LLT polynomials $\Gcal _{\nubold }(X;q)$.  Since \cite{GrojHaim07} is
unpublished, in \cite[\S 4]{BHMPS-paths} we gave a self-contained
treatment (with some improvements) of the case when $G = \GL _{l}$ and
$L = T$ is the torus, which corresponds to $\Gcal _{\nubold }(X;q)$
when $\nubold $ is a tuple of single-row skew diagrams.  Building on
\cite{BHMPS-paths}, we extend this here to the case of any Levi
subgroup $L = \GL _{r_{1}}\times \cdots \times \GL _{r_{k}} \subseteq
\GL _{l}$ and any $\Gcal _{\nubold }(X;q)$.

\subsection{Hecke algebra and root system preliminaries}
\label{ss:Hecke}

We set $\kk = \QQ (q,t)$ as in \S \ref{s:Catalanimals}.  The algebra
of Laurent polynomials $\kk [z_{1}^{\pm 1},\ldots,z_{l}^{\pm }]$ is
the group algebra of the weight lattice of $\GL _{l}$, with monomials
$\zz ^{\lambda } = z_{1}^{\lambda _{1}}\cdots z_{l}^{\lambda _{l}}$
corresponding to weights $\lambda \in \ZZ ^{l}$.  As in \S
\ref{s:Catalanimals}, we denote the roots by $\alpha _{ij} =
\varepsilon _{i} - \varepsilon _{j}$. For simple roots we abbreviate
this to $\alpha _{i} = \alpha _{i,i+1}$.

The Weyl group $S_{l}$ acts by permuting the variables, with Coxeter
generators (simple reflections) given by the transpositions $s_{i} =
(i \leftrightarrow i+1)$.  Given $w\in S_{l}$, we let $\ell (w)$
denote the length of a reduced factorization $w = s_{i_{1}}\cdots
s_{i_{\ell}}$; this is also the number of inversions in $w$.  The
longest element of $S_{l}$ or any finite Coxeter group is denoted
$w_{0}$.  Usually it will be clear from the context what group $w_{0}$
belongs to; otherwise we indicate it with a superscript such as
$w_{0}^{l}$.

The Demazure-Lusztig operators
\begin{equation}\label{e:Ti}
T_{i} = q\, s_{i} + (1 - q) \frac{1}{1- \zz ^{-\alpha _{i}}} (s_{i} - 1)
\end{equation}
generate an action of the Hecke algebra $\Hcal (S_{l})$ on $\kk
[z_{1}^{\pm 1},\ldots,z_{l}^{\pm }]$.  We have normalized them so that
$(T_{i} - q)(T_{i} + 1) = 0$.  As usual, for $w\in S_{l}$, we set
$T_{w} = T_{i_{1}}\cdots T_{i_{\ell }}$, where $w = s_{i_{1}}\cdots
s_{i_{\ell }}$ is a reduced factorization.

We use an overbar $\overline{\, \vphantom{t} \cdot \, }$ to signify
inverting the variables $q$, $t$, $z_{i}$; thus
\begin{equation}\label{e:T-bar-def}
\overline{T_{i}} = q^{-1}\, s_{i} + (1 - q^{-1}) \frac{1}{1-
\zz ^{\alpha _{i}}} (s_{i} - 1).
\end{equation}
One can then check that
\begin{equation}\label{e:T-bar}
\overline{T_{i}} = T_{i}^{-1}; \quad \text{hence} \quad
\overline{T_{w}} = T_{w^{-1}}^{-1}.
\end{equation}

Given a composition $\rr = (r_{1},\ldots,r_{k})$ of $l$, we denote the
corresponding Levi subgroup of $\GL _{l}$ and its Weyl group (which is
a Young subgroup of $S_{l}$) by
\begin{gather}
\label{e:GL_r}
\GL _{\rr } = \GL _{r_{1}} \times \cdots \times \GL _{r_{k}} \subseteq
 \GL _{l},\\ 
\label{e:S_r}
S_{\rr } = S_{r_{1}} \times \cdots \times S_{r_{k}}\subseteq S_{l}.
\end{gather}
Note that $R_{+}(\GL _{\rr })$ is the set of positive roots $\alpha
_{ij} \in R_{+}(\GL _{l})$ such that $i$, $j$ are in the same block of
the partition of $[l]$ into intervals of lengths $r_{1},\ldots,
r_{k}$.  

Here we have implicitly taken $\rr$ to be a strict composition with
all entries $r_{i}>0$.  If $\rr$ is a {\em weak composition} with
entries $r_{i} = 0$ allowed, we define $\GL _{\rr} = \GL _{\sS }$
where $\sS = (r_{i_{1}},\ldots,r_{i_{j}})$ is the subsequence of
non-zero entries in $\rr $.

The Levi subgroup $\GL _{\rr }$ has the same weight lattice $\ZZ ^{l}$
as $\GL _{l}$.  A weight $\lambda $ is dominant (resp.\ dominant and
regular) for $\GL _{\rr }$ iff $\lambda _{i}\geq \lambda _{i+1}$
(resp.\ $\lambda _{i} > \lambda _{i+1}$) for all $i$ such that $\alpha
_{i}\in R_{+}(\GL _{\rr})$, or equivalently such that $s_{i}\in S_{\rr
}$.  We denote the set of dominant weights by $X^{+}(\GL _{\rr })$ and
the set of regular dominant weights by $X^{++}(\GL _{\rr })$.

We write $\rho _{\rr }$ (or just $\rho $ if $\rr = (l)$) for a weight
such that
\begin{equation}\label{e:rho-r}
\langle \alpha _{i}^{\vee }, \rho _{\rr } \rangle \defeq (\rho _{\rr })_{i}
- (\rho _{\rr })_{i+1} = 1 \quad \text{for every simple root
$\alpha _{i}\in R_{+}(\GL _{\rr })$}.
\end{equation}
Such a weight is unique up to adding an $S_{\rr }$ invariant
weight.  When we use this notation, the choice of $\rho _{\rr }$ will
either be fixed or make no difference.

We define the semi-symmetrization operator for $\GL _{\rr }$ by 
means of the following lemma.

\begin{lemma}\label{lem:delta_r}
For any composition $\rr $ of $l$, there is an operator $\semis _{\rr
}$ on $\kk [z_{1}^{\pm 1},\ldots,z_{l}^{\pm 1}]$ given by either of
two equivalent formulas
\begin{align}\label{e:delta_r}
\semis _{\rr } & = \frac{1}{\prod _{\alpha \in R_{+}(\GL _{\rr })} (1
- \zz ^{-\alpha })} \sum _{w\in S_{\rr }} (-1)^{\ell (w)} w
\\
\label{e:delta_r-q-form}
& = \frac{q^{\ell (w_{0}^{\rr })}}{\prod _{\alpha \in R_{+}(\GL _{\rr
})} (1 - q\, \zz ^{-\alpha })} \sum _{w\in S_{\rr }} (-q)^{-\ell (w)}
T_{w}.
\end{align}
\end{lemma}

\begin{proof}
Fixing a choice of $\rho _{\rr }$, formula \eqref{e:delta_r} can also
be written
\begin{equation}\label{e:delta_r-vs-sigma_r}
\semis _{\rr } = \zz ^{\rho _{\rr }} \, \sigmabold _{\rr }\, \zz ^{-\rho
_{\rr }},
\end{equation}
where $\sigmabold _{\rr }$ is the Weyl symmetrization operator for
$\GL _{\rr }$ in \eqref{e:Weyl-symmetrizer}.  In particular,
the operator $\semis _{\rr }$ defined
by \eqref{e:delta_r} acts on $\kk [z_{1}^{\pm 1},\ldots,z_{l}^{\pm
1}]$.  Let $\semis _{\rr }^{q}$ denote the operator defined by
\eqref{e:delta_r-q-form}.  We are to prove that $\semis _{\rr }^{q} =
\semis _{\rr }$.

Let $V = \kk [z_{1}^{\pm 1},\ldots,z_{l}^{\pm 1}]$.  For each $i$ such
that $s_{i}\in S_{\rr }$, let $W_{i} = V^{s_{i}} = \{f\in V\mid s_{i}
f = f \}$ be the subspace of $s_{i}$ invariant functions.  Let $A_{\rr
} = \sum _{w\in S_{\rr }} (-1)^{\ell (w)} w$ denote the
antisymmetrization operator.  It is a general property of Coxeter group
representations that $\sum_{i}W_{i}$ is an invariant subspace and that
$V/\sum _{i}W_{i}$ carries the sign representation, thus
every $w\in S_{\rr }$ acts on $V/\sum
_{i}W_{i}$ as $(-1)^{\ell (w)}$, and $A_{\rr }$ acts as the scalar
$|S_{\rr }|$.  This implies, first, that $A_{\rr }$ is surjective on
$V/\sum _{i}W_{i}$, that is, the space of all antisymmetric functions
spans $V/\sum _{i}W_{i}$, and, second, that if $V'\subseteq V$ is
any subspace such that $A_{\rr }V'$ is the space of all antisymmetric
functions, then $V'$ spans $V/\sum _{i}W_{i}$, that is, $V = V'+\sum
_{i} W_{i}$.  In particular, this holds with $V' = \zz ^{\rho _{\rr }}
V^{S_{\rr }}$, since every antisymmetric $f\in V$ has the form $f =
A_{\rr }\, \zz ^{\rho _{\rr }} g$, where $g\in V$ is $S_{\rr }$
invariant.

Since $T_{i}f =q\, f$ for $f\in W_{i}$, both operators $\semis _{\rr
}$ and $\semis _{\rr }^{q}$ kill the subspaces $W_{i}$.  Both
operators also commute with multiplication by any $S_{\rr }$ invariant
function $g$.  Hence, to prove $\semis _{\rr }^{q} = \semis _{\rr } $
it suffices to show that $\semis _{\rr }^{q} \zz ^{\rho _{\rr }}=
\semis _{\rr } \zz ^{\rho _{\rr}}$.  By \eqref{e:delta_r-vs-sigma_r},
we have $\semis _{\rr }\zz ^{\rho _{\rr }} = \zz ^{\rho _{\rr }}$.
Meanwhile, $\semis _{\rr }^{q}\zz ^{\rho _{\rr }} = \zz ^{\rho _{\rr
}}$ is equivalent to the well-known identity \cite{BakDunFor00}
\begin{equation}\label{e:rho-Hecke-anti}
\sum _{w\in S_{\rr }} (-q)^{-\ell (w)} T_{w} \, \zz ^{\rho _{\rr }} =
q^{-\ell (w_{0}^{\rr })}\, \zz ^{\rho _{\rr }} \prod _{\alpha \in
R_{+}(\GL _{\rr })} (1 - q\, \zz ^{-\alpha }).
\qedhere
\end{equation}
\end{proof}

\begin{remark}\label{rem:Hecke-antisymmetry}
For $i$ such that $s_{i}\in S_{\rr }$, the Hecke algebra
antisymmetrization operator $A_{\rr }^{q} = \sum _{w\in S_{\rr }}
(-q)^{-\ell (w)} T_{w}$ in \eqref{e:rho-Hecke-anti} can be factored in
the form $B\cdot (T_{i}-q)$ for an element $B$ of the Hecke algebra,
and therefore satisfies $A_{\rr }^{q}\, T_{i} = -A_{\rr }^{q}$.  More
generally, this implies $A_{\rr }^{q}\, T_{w}^{\pm 1} = (-1)^{\ell
(w)} A_{\rr }^{q}$, for all $w\in S_{\rr }$, and consequently also
$\semis _{\rr }\, T_{w}^{\pm 1} = (-1)^{\ell (w)} \semis _{\rr }$.
\end{remark}

\subsection{Semi-symmetric Hall-Littlewood polynomials}
\label{ss:semi-symmetric}

As in \cite[(72)]{BHMPS-paths}, we define non-symmetric
Hall-Littlewood polynomials for $\GL _{l}$ by
\begin{equation}\label{e:E_lambda}
E_{\lambda }(\zz ;q) = q^{-\ell(w)} T_{w} \, \zz ^{\lambda _{+}},\\ 
\end{equation}
where $\lambda = w(\lambda _{+})$ with $w\in S_{l}$ and $\lambda _{+}$
dominant.  If $\lambda _{+}$ has non-trivial stabilizer, $w$ is not
unique, but the formula does not depend on the choice.  For $\sigma \in
S_{l}$, we also define twisted versions
\begin{align}
\label{e:E-twist} E^{\sigma }_{\lambda }(\zz ;q) & = q^{|{\Inv (\sigma
^{-1})\cap \Inv (\lambda +\epsilon \rho )}|}\, \overline{T_{\sigma
}}\, E_{\sigma
^{-1}(\lambda )}(\zz ;q)\\
\label{e:F-twist}
F^{\sigma }_{\lambda }(\zz ;q) & = \overline{E^{\sigma
w_{0}}_{-\lambda }(\zz ;q)},
\end{align}
where $\Inv ((a_{1},\ldots,a_{l})) = \{(i<j) \mid a_{i}>a_{j} \} $,
and $\epsilon $ is small, so $\Inv (\sigma ^{-1}) = \{(i<j)\mid \sigma
^{-1}(i) > \sigma ^{-1}(j) \}$ and $\Inv (\lambda +\epsilon \rho ) =
\{(i<j) \mid \lambda _{i}\geq \lambda _{j} \}$.
Note that for $\sigma =1$, $E^{\sigma }_{\lambda }(\zz ;q)$ reduces to
the untwisted $E_{\lambda }(\zz ;q)$.

From the definition, one can verify the recurrence
\cite[(76)]{BHMPS-paths}
\begin{equation}\label{e:E-recurrence}
E^{\sigma }_{\lambda } =
\begin{cases}
q^{-\chi (\lambda _{i}\leq \lambda _{i+1})}\, T_{i}\, E^{s_{i}\sigma
}_{s_{i}\lambda },
 &	  s_{i}\sigma >\sigma, \\
q^{\chi (\lambda _{i}\geq \lambda _{i+1})}\, T_{i}^{-1}\, E^{s_{i}\sigma
}_{s_{i}\lambda }, & s_{i}\sigma <\sigma.
\end{cases}
\end{equation}
The $E^{\sigma }_{\lambda }$ are determined by this recurrence and the
initial condition $E^{\sigma }_{\lambda } = \zz ^{\lambda }$ for all
$\sigma $ if $\lambda$ is a dominant weight.

\begin{defn}\label{def:E-F-semi}
(i) Given a composition $\rr = (r_{1},\ldots,r_{k})$ and a permutation
$\sigma \in S_{k}$, let $l = |\rr | = r_{1}+\cdots +r_{k}$ and define
$\sigmahat \in S_{l}$ to be the permutation that carries intervals of
lengths $\sigma ^{-1}(\rr ) = (r_{\sigma (1)},\ldots,r_{\sigma (k)})$
to intervals of lengths $\rr $ in the order given by $\sigma $.  More
precisely,
\begin{equation}\label{e:sigma-hat}
\sigmahat (r_{\sigma (1)} + \dots + r_{\sigma (i-1)} + j) = 
r_{1} + \dots + r_{\sigma (i)-1} + j
\end{equation}
for $i=1,\ldots,k$ and $j = 1,\ldots,r_{\sigma (i)}$.

(ii) Given $\rr $ and $\sigma $ as in (i), and $\mu \in X^{++}(\GL
_{\rr })$ a regular dominant weight for $\GL _{\rr }$, we define {\em
semi-symmetric Hall-Littlewood polynomials}
\begin{equation}\label{e:E-F-semi}
E^{\sigma }_{\rr ,\mu }(\zz ;q) = \semis _{\rr }\, E^{\sigmahat }_{\mu
}(\zz ;q),\qquad F^{\sigma }_{\rr ,\mu }(\zz ;q) = \semis _{\rr }\,
F^{\sigmahat }_{\mu }(\zz ;q).
\end{equation}
\end{defn}

\begin{remark}\label{rem:sigma-hat-and-weak-r}
(i) For simplicity, we have suppressed $\rr $ from the notation for
$\sigmahat $.  Although this should not usually cause confusion, one
should note that $\sigmahat\, ^{-1}$ is not given by $\widehat{\sigma
^{-1}}$ for the same $\rr $, but rather by $\widehat{\sigma ^{-1}}$
defined relative to the composition $\sigma ^{-1}(\rr )$.

For example, if $\sigma =(2,3,1)$ in one-line notation, and $\rr
=(1,4,3)$, then $\sigma ^{-1} = (3,1,2)$ and $\sigma ^{-1}(\rr ) =
(4,3,1)$.  Partitioning the set $\{1,\ldots,8\}$ into intervals $I_{1}
= \{1 \}$, $I_{2} = \{2,3,4,5 \}$, $I_{3} = \{6,7,8 \}$ of lengths
$\rr $ and intervals $J_{1} = \{1, 2, 3 ,4 \}$, $J_{2} = \{5,6,7 \}$,
$J_{3} = \{8 \}$ of lengths $\sigma ^{-1}(\rr )$, the permutation
$\sigmahat = (2,3,4,5,6,7,8,1)$ carries $J_{1}$, $J_{2}$, $J_{3}$ to
$I_{2} =I_{\sigma (1)}$, $I_{3} = I_{\sigma (2)}$, $I_{1} = I_{\sigma
(3)}$.  The inverse permutation $\sigmahat ^{-1}$ that carries
$I_{1}$, $I_{2}$, $I_{3}$ back to $J_{3}$, $J_{1}$, $J_{2}$ is
$\widehat{\sigma ^{-1}}$ defined relative to the composition $\sigma
^{-1}(\rr ) = (4,3,1)$ that gives the intervals $J_{i}$, rather than
the original composition $\rr = (1,4,3)$ that gave the intervals
$I_{i}$.

(ii) If $\rr $ is a weak composition, so $\GL _{\rr} = \GL _{\sS }$ 
where $\sS = (r_{i_1},\ldots, r_{i_j})$ 
is the subsequence of non-zero entries
in $\rr $, then the definitions of $E^{\sigma }_{\rr ,\mu }(\zz ;q)$ and
$F^{\sigma }_{\rr ,\mu }(\zz ;q)$ reduce to $E^{\sigma }_{\rr ,\mu
}(\zz ;q) = E^{\tau }_{\sS ,\mu }(\zz ;q)$, $F^{\sigma }_{\rr ,\mu
}(\zz ;q) = F^{\tau }_{\sS ,\mu }(\zz ;q)$ where $\tau \in S_{j}$ is
the permutation such that $\tau ^{-1}(1),\ldots,\tau ^{-1}(j)$ are in
the same relative order as $\sigma ^{-1}(i_{1}),\ldots,\sigma
^{-1}(i_{j})$.
\end{remark}

\begin{example}\label{ex:special-semi-syms}
(i) At $q = 1$, we have $E^{\sigma }_{\mu }(\zz ;1) = \zz ^{\mu }$ for
any $\sigma $.  Using \eqref{e:delta_r-vs-sigma_r}, it follows that
\begin{equation}\label{e:q=1}
\zz ^{-\rho _{\rr }} E^{\sigma }_{\rr ,\lambda +\rho _{\rr }}(\zz ;1)
= \chi _{\lambda }(\GL _{\rr })
\end{equation}
if $\lambda \in X^{+}(\GL _{\rr })$
is a dominant weight for $\GL_{\rr}$,
where $\chi _{\lambda }(\GL _{\rr })$ is the irreducible $\GL _{\rr }$
character with highest weight $\lambda $.

(ii) For $\rr = (l)$, we must have $\sigma =1\in S_{1}$.  If $\lambda
$ is a dominant weight for $\GL _{l}$,
then $E_{\lambda +\rho}(\zz ;q)
= \zz ^{\lambda +\rho }$, which implies $\zz ^{-\rho } E^{1}_{(l),
\lambda +\rho }(\zz ;q) = \chi _{\lambda }$, independent of $q$. Note
that this is quite different from the usual symmetric Hall-Littlewood
polynomial $P_{\lambda }(\zz; q)$.
\end{example}

We develop some initial properties of these polynomials for later use,
beginning with expressions for the semi-symmetric polynomials
$F^{\sigma }_{\rr }(\zz ;q)$ in terms of $E^{\sigma }_{\rr }(\zz ;q)$.

\begin{lemma}\label{lem:F-vs-E}
Given $\rr = (r_{1},\ldots,r_{k})$, $\sigma \in S_{k}$ and $\mu \in
X^{++}(\GL _{\rr })$, we have
\begin{align}\label{e:F-vs-E-bar}
F^{\sigma }_{\rr ,\mu }(\zz ,q) & = \zz ^{\rho _{\rr } - w_{0}^{\rr
}(\rho _{\rr })} \overline{E^{\sigma w_{0}^{k}}_{\rr ,-w_{0}^{\rr
}(\mu )}(\zz ;q)}
\\
\label{e:F-vs-E-w0} & = \zz ^{\rho _{\rr } - w_{0}^{\rr }(\rho _{\rr
})} w_{0}\big( E^{w_{0}^{k}\sigma }_{w_{0}^{k}(\rr ),w_{0}w_{0}^{\rr
}(\mu )}(\zz ;q^{-1}) \big).
\end{align}
Note that $\rho _{\rr } - w_{0}^{\rr }(\rho _{\rr }) = \sum _{\alpha
\in R_{+}(\GL _{\rr })}\alpha $ does not depend on the choice of $\rho
_{\rr }$.
\end{lemma}

\begin{proof}
From \eqref{e:delta_r} we find $\overline{\semis _{\rr }} = (-1)^{\ell
(w_{0}^{\rr })} \zz ^{w_{0}^{\rr }(\rho _{\rr }) - \rho _{\rr }}
\semis _{\rr }$.  Applying $\semis _{\rr }$ on both sides of the
definition $F^{\sigmahat }_{\mu } = \overline{E^{\sigmahat
w_{0}}_{-\mu }}$ therefore gives
\begin{equation}\label{e:F-vs-E-bar-step}
F^{\sigma }_{\rr ,\mu } = (-1)^{\ell (w_{0}^{\rr })} \zz ^{\rho _{\rr
} - w_{0}^{\rr }(\rho _{\rr })} \overline{\semis _{\rr} E^{\sigmahat
w_{0}}_{-\mu }}.
\end{equation}
Since $\sigmahat w_{0}$ is maximal in its coset $S_{\rr } \sigmahat
w_{0}$, and $-\mu $ is regular and anti-dominant for $\GL _{\rr }$, it
follows by repeated use of the recurrence \eqref{e:E-recurrence} that
$E^{\sigmahat w_{0}}_{-\mu } = \overline{T_{w_{0}^{\rr }}}
E^{w_{0}^{\rr }\sigmahat w_{0}}_{-w_{0}^{\rr }(\mu )}$.  Now
$\widehat{\sigma w_{0}^{k}} = w_{0}^{\rr }\sigmahat w_{0}$, and
$\semis _{\rr}T_{i}^{\pm 1} = - \semis _{\rr }$ for $s_{i}\in S_{\rr
}$ by Remark~\ref{rem:Hecke-antisymmetry}, hence $\semis _{\rr
}\overline{T_{w_{0}^{\rr }}} = (-1)^{\ell (w_{0}^{\rr })} \semis _{\rr
}$.  Combining these, the right hand side of \eqref{e:F-vs-E-bar-step}
becomes $\zz ^{\rho _{\rr } - w_{0}^{\rr }(\rho _{\rr })}
\overline{\semis _{\rr} E^{\widehat{\sigma w_{0}^k}}_{-w_{0}^{\rr
}(\mu )}} = \zz ^{\rho _{\rr } - w_{0}^{\rr }(\rho _{\rr })}
\overline{E^{\sigma w_{0}^k}_{\rr ,-w_{0}^{\rr }(\mu )}}$, giving
\eqref{e:F-vs-E-bar}.

To prove \eqref{e:F-vs-E-w0}, we use the identity
\begin{equation}\label{e:w0-identity}
\overline{E^{\tau }_{-\lambda }(\zz ;q)} = w_{0}\big( E^{w_{0}\tau
w_{0}}_{w_{0}(\lambda )}(\zz ;q^{-1}) \big).
\end{equation}
Like the equivalent identity \cite[(109)]{BHMPS-delta}, one can prove
\eqref{e:w0-identity} by verifying that after applying $w_{0}$, both
sides are characterized by the recurrence \eqref{e:E-recurrence} with
the variables reversed and inverted.  We also observe that
\eqref{e:delta_r} implies $\overline{\semis _{\rr }}w_{0} =
w_{0}\semis _{\rr '}$, where $\rr ' = w_{0}^{k}(\rr ) =
(r_{k},\ldots,r_{1})$, and that $w_{0}(\widehat{\sigma
w_{0}^{k}})w_{0} = \widehat{w_{0}^{k}\sigma }$, where
$\widehat{w_{0}^{k}\sigma }$ is defined with respect to $\rr '$ rather
than $\rr $.  Then \eqref{e:F-vs-E-w0} follows from
\eqref{e:F-vs-E-bar} and
\begin{equation}\label{e:F-vs-E}
\overline{E^{\sigma w_{0}^{k}}_{\rr ,-w_{0}^{\rr }(\mu )}(\zz ;q)} =
\overline{\semis _{\rr } E^{\widehat{\sigma w_{0}^{k}}}_{-w_{0}^{\rr
}(\mu )}(\zz ;q)} = w_{0} \big( \semis _{\rr
'}E^{\widehat{w_{0}^{k}\sigma }}_{w_{0}w_{0}^{\rr }(\mu )}(\zz
;q^{-1}) \big) = w_{0} \big(E^{w_{0}^{k}\sigma }_{\rr
',w_{0}w_{0}^{\rr }(\mu )}(\zz ;q^{-1}) \big) \,.
\end{equation}
\end{proof}

By \cite[Corollary~4.3.1]{BHMPS-paths}, $E^{\sigma }_{\lambda }$
has the monic and triangular form
\begin{equation}\label{e:ns-HL-triangular}
E^{\sigma }_{\lambda }(\zz ;q) = \zz ^{\lambda } + \sum _{\mu <\lambda
} c_{\lambda ,\mu }(q) \, \zz ^{\mu }
\end{equation}
with respect to a suitable partial ordering $<$ on the weight lattice
$\ZZ ^{l}$.  If $\mu $ is a regular weight for $\GL _{\rr }$, then
$\zz ^{-\rho _{\rr }} \semis _{\rr }(\zz ^{\mu }) = \pm \chi
_{\nu}(\GL _{\rr })$, where $\nu +\rho _{\rr }$ is the unique dominant
(and regular) weight in the orbit $S_{\rr }\cdot \mu $.  The ordering
$<$ has the property that the dominant weight for $\GL _{\rr }$ in any
$S_{\rr }$ orbit is the unique minimal element in that orbit; hence
$\nu +\rho _{\rr }\leq \mu $ in this case.  If $\mu $ is not regular
for $\GL _{\rr }$, then $\semis _{\rr }(\zz^{\mu} ) = 0$.

For every $\lambda \in X^{+}(\GL _{\rr })$ it now follows from
\eqref{e:ns-HL-triangular} that $\zz ^{-\rho _{\rr }} E^{\sigma }_{\rr
,\lambda + \rho _{\rr } }$ has the form
\begin{equation}\label{e:semi-sym-triangular}
\zz ^{-\rho _{\rr }} E^{\sigma }_{\rr ,\lambda + \rho _{\rr } }(\zz
;q) = \chi _{\lambda
}(\GL _{\rr }) + \sum _{\substack{\nu  \in X^{+}(\GL _{\rr })\\
\nu +\rho _{\rr } <\lambda +\rho _{\rr } }} a_{\lambda ,\nu }(q) \,
\chi _{\nu }(\GL _{\rr }).
\end{equation}
More precisely, given the choice of $\rho _{\rr }$,
\eqref{e:semi-sym-triangular} holds for all $\lambda \in X^{+}(\GL
_{\rr })$, although the coefficients $a_{\lambda ,\nu }(q)$ and the
set of weights $\nu $ that occur depend on $\rho _{\rr }$.  In
particular, for any fixed choice of $\rho _{\rr }$, it follows that
\begin{equation}\label{e:E-basis}
\{\zz ^{-\rho _{\rr }} E^{\sigma
}_{\rr ,\mu } \mid \mu \in X^{++}(\GL _{\rr }) \}
\end{equation}
is a basis of $\kk [z_{1}^{\pm 1},\ldots,z_{l}^{\pm }]^{S_{\rr }}$.
Then \eqref{e:F-vs-E-bar} implies that
\begin{equation}\label{e:F-basis}
\{\zz ^{-\rho _{\rr }} F^{\sigma
}_{\rr ,\mu } \mid \mu \in X^{++}(\GL _{\rr }) \}
\end{equation}
is also a basis.  Note that $\kk [z_{1}^{\pm 1},\ldots,z_{l}^{\pm
1}]^{S_{\rr }}$ is the algebra of virtual $\GL _{\rr }$ characters
with coefficients in $\kk $.

\begin{remark}\label{rem:coefficient-rings}
The coefficients $c_{\lambda ,\mu }(q)$ in \eqref{e:ns-HL-triangular},
and therefore also $a_{\lambda ,\nu }(q)$ in
\eqref{e:semi-sym-triangular}, are in $\ZZ [q^{-1}]$.  Hence,
$E^{\sigma }_{\rr ,\mu }(\zz ;q)$ and $F^{\sigma }_{\rr, \mu }(\zz
;q)$ have coefficients in $\ZZ [q^{-1}]$ and $\ZZ [q]$, respectively,
and (\ref{e:E-basis}--\ref{e:F-basis}) are free module bases over
these coefficient rings in place of $\kk $.
\end{remark}

\begin{lemma}\label{lem:Tw-on-EF-sigmahat}
Given $\rr $, $\sigma $ and $\sigmahat $ as in
Definition~\ref{def:E-F-semi}, $\mu \in X^{++}(\GL _{\rr })$, and
$w\in S_{\rr }$, we have
\begin{align}\label{e:Tw-on-E-sigmahat}
E^{\sigmahat }_{w(\mu )} & = q^{-\ell (w)} T_{w} E^{\sigmahat }_{\mu },\\
\label{e:Tw-on-F-sigmahat}
F^{\sigmahat }_{w(\mu )} & = \overline{T_{w}} F^{\sigmahat }_{\mu }.
\end{align}
\end{lemma}

\begin{proof}
The case $\sigma =1$ of \eqref{e:Tw-on-E-sigmahat} follows from the
definition \eqref{e:E_lambda} of $E_{\mu }$,
because $\mu \in X^{++}(\GL _{\rr })$ implies that if $v\in
S_{l}$ is such that $\mu =v(\mu _{+})$, and $w\in S_{\rr }$, then
$w\cdot v$ is a reduced factorization of $wv$.
For general $\sigma $, let $v = \sigmahat ^{-1} w \sigmahat $, and
note that $v\in S_{\sigma ^{-1}(\rr )}$.  Since $\sigmahat ^{-1}$ is
minimal in both of its cosets $S_{\sigma ^{-1}(\rr )} \, \sigmahat
^{-1}$ and $\sigmahat ^{-1} \, S_{\rr }$, each side of $\sigmahat
^{-1}\cdot w = v\cdot \sigmahat ^{-1}$ is a reduced factorization,
giving $T_{\sigmahat ^{-1}} T_{w} = T_{v} T_{\sigmahat ^{-1}}$, or
equivalently $T_{w} \overline{T_{\sigmahat }} = \overline{T_{\sigmahat
}} T_{v}$.  We also have $\ell (v) = \ell (w)$.  Then
\begin{multline}\label{e:Tw-E}
q^{-\ell (w)}\, T_{w}\, E^{\sigmahat }_{\mu } = q^{e - \ell (w)}\,
T_{w} \overline{T_{\sigmahat }}\, E_{\sigmahat ^{-1}(\mu )} = q^{e}\,
\overline{T_{\sigmahat }}\, q^{-\ell (v)}\, T_{v}\, E_{\sigmahat
^{-1}(\mu )}
\\
= q^{e}\, \overline{T_{\sigmahat }}\, E_{v \sigmahat ^{-1} (\mu )} =
q^{e}\, \overline{T_{\sigmahat }}\, E_{\sigmahat ^{-1} w (\mu )},
\end{multline}
where $e = |{\Inv (\sigmahat ^{-1})\cap \Inv (\mu +\epsilon \rho )}|$.
For the third equality, we used $\sigmahat ^{-1}(\mu )\in X^{++}(\GL
_{\sigma ^{-1}(\rr )})$ and the $\sigma = 1$ case.  Since $\sigmahat
^{-1}$ is increasing on intervals of lengths $r_{1},\ldots,r_{k}$, we
have $e = |{\Inv (\sigmahat ^{-1})\cap \Inv (w(\mu) +\epsilon \rho
)}|$ for any $w\in S_{\rr }$.  The last formula in \eqref{e:Tw-E}
therefore reduces to $E^{\sigmahat }_{w(\mu )}$.

For \eqref{e:Tw-on-F-sigmahat}, let $u = w_{0} \sigmahat ^{-1} w
\sigmahat w_{0}$ (here $w_{0} = w_{0}^{l}\in S_{l}$), and note that
$u\in S_{w_{0}^{k}\sigma ^{-1}(\rr )}$.  Since $w_{0}\sigmahat ^{-1}$
is maximal in both of its cosets $w_{0}\sigmahat ^{-1} S_{\rr }$ and
$S_{w_{0}^{k}\sigma ^{-1}(\rr )} w_{0}\sigmahat ^{-1}$, the
factorizations $(w_{0}\sigmahat ^{-1}w^{-1})\cdot w$ and $u\cdot
(u^{-1} w_{0}\sigmahat ^{-1})$ are reduced, giving
$T_{w}T_{w_{0}\sigmahat ^{-1}}^{-1} = T_{w_{0}\sigmahat
^{-1}w^{-1}}^{-1} = T_{u^{-1} w_{0}\sigmahat ^{-1}}^{-1} =
T_{w_{0}\sigmahat ^{-1}}^{-1} T_{u}$, or equivalently $T_{w}
\overline{ T_{\sigmahat w_{0}}} = \overline{ T_{\sigmahat w_{0}}}
T_{u}$.  Then
\begin{multline}\label{e:Tw-F}
T_{w}\, E^{\sigmahat w_{0}}_{-\mu } = q^{d} \, T_{w} \, \overline{
T_{\sigmahat w_{0}}} \, E_{- w_{0}\sigmahat ^{-1}(\mu )} = q^{d} \,
\overline{ T_{\sigmahat w_{0}}} \, T_{u} \, E_{- w_{0}\sigmahat
^{-1}(\mu )}
\\
= q^{d + \ell (u)} \, \overline{ T_{\sigmahat w_{0}}} \, E_{- u
w_{0}\sigmahat ^{-1}(\mu )} = q^{d + \ell (u)} \, \overline{
T_{\sigmahat w_{0}}} \, E_{- w_{0}\sigmahat ^{-1} w(\mu )},
\end{multline}
where $d = |{\Inv (w_{0}\sigmahat ^{-1})\cap \Inv (-\mu +\epsilon \rho
)}|$. For the third equality, we used the $\sigma = 1$ case of
\eqref{e:Tw-on-E-sigmahat} with $-w_{0}\sigmahat ^{-1}(\mu )\in
X^{++}(\GL _{w_{0}^{k}\sigma ^{-1}(\rr )})$.  Since $w_{0}\sigmahat
^{-1}$ is decreasing on intervals of lengths $r_{1},\ldots,r_{k}$,
changing $-\mu $ to $-w(\mu )$ in the formula for $d$ creates $\ell
(w) = \ell (u)$ new inversions, giving $d+\ell (u) = |{\Inv
(w_{0}\sigmahat ^{-1})\cap \Inv (-w(\mu )+\epsilon \rho )}|$.  The
last formula in \eqref{e:Tw-F} now reduces to $E^{\sigmahat
w_{0}}_{-w(\mu )}$, showing that $E^{\sigmahat w_{0}}_{-w(\mu )} =
T_{w}\, E^{\sigmahat w_{0}}_{-\mu }$.  Taking $\overline{\,
\vphantom{t} \cdot \, }$ on both sides gives
\eqref{e:Tw-on-F-sigmahat}.
\end{proof}

Finally, since $\semis _{\rr }$ and  $T_{i}$ commute with
multiplication by $z_{1}\cdots z_{l}$, we have the identities
\begin{equation}\label{e:shift}
(z_{1}\cdots z_{l})^{m}\, E^{\sigma }_{\rr ,\mu }(\zz ;q) = E^{\sigma
}_{\rr ,(m^{l})+\mu }(\zz ;q),\qquad (z_{1}\cdots z_{l})^{m}\,
F^{\sigma }_{\rr ,\mu }(\zz ;q) = F^{\sigma }_{\rr ,(m^{l})+\mu }(\zz
;q).
\end{equation}

\subsection{Orthogonality}
\label{ss:orthogonality}

For $f\in \kk [z_{1}^{\pm 1},\ldots,z_{l}^{\pm 1}]^{S_{\rr }}$, let
$\langle 1_{\GL _{\rr }} \rangle\, f$ denote the coefficient of the
trivial character when $f$ is expanded in terms of irreducible $\GL
_{\rr }$ characters.  The formula
\begin{equation}\label{e:1-G-in-Weyl}
\langle 1_{\GL _{\rr }} \rangle\, \sigmabold _{\rr } (f) = \langle \zz
^{0} \rangle\, f\cdot \prod \nolimits_{\alpha \in R_{+}(\GL _{\rr })}
(1-\zz ^{\alpha })
\end{equation}
holds for any $f\in \kk [z_{1}^{\pm 1},\ldots,z_{l}^{\pm 1}]$, as can
be verified by reducing to the case $f = \zz ^{\lambda }$, for which
both sides become $(-1)^{\ell (w)}$ if $\lambda +\rho _{\rr } = w(\rho
_{\rr })$ for $w\in S_{\rr }$, or zero otherwise.
If $f$ is $S_{\rr }$ invariant, then $\sigmabold _{\rr }(f) = f
\sigmabold _{\rr }(1) = f$, and we obtain
\begin{equation}\label{e:1-GL_r-coef}
\langle 1_{\GL _{\rr }} \rangle \, f = \langle \zz ^{0} \rangle\, f
\cdot \prod \nolimits_{\alpha \in R_{+}(\GL _{\rr })} (1-\zz ^{\alpha
}).
\end{equation}

We define a symmetric inner product on $\kk [z_{1}^{\pm
1},\ldots,z_{l}^{\pm 1}]^{S_{\rr }}$ by
\begin{equation}\label{e:r-inner-prod}
\langle f,g \rangle^{\rr }_{q} \, \defeq \langle 1_{\GL _{\rr }}
\rangle\, f\, g\, \prod \nolimits_{\alpha \in R_{+}\setminus R_{+}(\GL
_{\rr }) } \frac{1 - \zz ^{\alpha }}{1 - q^{-1} \zz ^{\alpha }}.
\end{equation}
For $\rr = (1^{l})$ this reduces to the inner product
\begin{equation}\label{e:HL-inner-prod}
\langle f,g \rangle _{q} = \langle \zz ^{0} \rangle\, f\, g\, \prod
\nolimits_{\alpha \in R_{+}} \frac{1 - \zz ^{\alpha }}{1 - q^{-1} \zz
^{\alpha }}
\end{equation}
in \cite[Proposition~4.3.2]{BHMPS-paths}.  For general $\rr $, 
\eqref{e:1-GL_r-coef} implies that $\langle f,g \rangle^{\rr
}_{q}$ and $\langle f,g \rangle_{q} $ are related by
\begin{equation}\label{e:r-vs-HL-inner}
\langle f,g \rangle^{\rr }_{q} = \langle\, f,\,  g \cdot \prod
\nolimits_{\alpha \in R_{+}(\GL _{\rr })} (1-q ^{-1}\zz ^{\alpha })\,
\rangle _{q}.
\end{equation}

We remark that \eqref{e:r-inner-prod} and \eqref{e:HL-inner-prod} are
to be interpreted by expanding the factors $(1 - q^{-1} \zz ^{\alpha
})^{-1} = 1 + q^{-1} \zz ^{\alpha }+\cdots $ as geometric series,
yielding a power series in $q^{-1}$ over $\ZZ [z_{1}^{\pm
},\ldots,z_{l}^{\pm 1}]$, which is $S_{\rr }$ invariant in the case of
\eqref{e:r-inner-prod}.  Upon taking the coefficient $\langle 1_{\GL
_{\rr }} \rangle$ or $\langle \zz ^{0} \rangle$, only finitely many
terms in the series survive.  In the case of \eqref{e:HL-inner-prod},
this is clear, and for \eqref{e:r-inner-prod} it then follows from
\eqref{e:r-vs-HL-inner}.

\begin{prop}\label{prop:orthogonality}
Given $\rr $ and $\sigma $ as in Definition~\ref{def:E-F-semi}, and
any choice of $\rho _{\rr }$ as in \eqref{e:rho-r}, we have dual bases
of $\kk [z_{1}^{\pm 1},\ldots,z_{l}^{\pm 1}]^{S_{\rr }}$
\begin{equation}\label{e:orthogonality}
\big\langle \, \zz ^{-\rho _{\rr }} E^{\sigma }_{\rr ,\lambda },\,
\overline{\zz ^{-\rho _{\rr }} F^{\sigma }_{\rr ,\mu }} \,
\big\rangle ^{\rr }_{q} = \delta _{\lambda ,\mu }\qquad (\lambda ,\mu \in
X^{++}(\GL _{\rr })).
\end{equation}
\end{prop}

\begin{proof} We have already seen that the two sets 
$\{ \zz ^{-\rho _{\rr }} E^{\sigma }_{\rr ,\lambda } \mid \lambda \in
X^{++}(\GL _{\rr }) \}$ and $\{ \overline{\zz ^{-\rho _{\rr }}
F^{\sigma }_{\rr ,\mu }} \mid \mu \in X^{++}(\GL _{\rr }) \}$ are
bases of $\kk [z_{1}^{\pm 1},\ldots,z_{l}^{\pm 1}]^{S_{\rr }}$.  We
also note that $\langle \zz ^{-\rho _{\rr }} E^{\sigma }_{\rr ,\lambda
},\, \overline{\zz ^{-\rho _{\rr }} F^{\sigma }_{\rr ,\mu }}\rangle
^{\rr }_{q}$ is independent of the choice of $\rho _{\rr }$, since
$\langle f,g \rangle ^{\rr }_{q}$ is a function of $fg$.  The case
$\rr = (1^{l})$ is \cite[Proposition~4.3.2]{BHMPS-paths}.  We will use
this result to prove the general case.

For $\mu \in X^{++}(\GL _{\rr })$, define
\begin{align}\label{e:E-minus}
E^{\sigma ,-}_{\rr ,\mu } & =  \big( q^{-\ell (w_{0}^{\rr })} \prod
\nolimits _{\alpha \in R_{+}(\GL _{\rr })} (1- q \zz ^{-\alpha })\big)
E^{\sigma }_{\rr ,\mu } = \sum _{w\in S_{\rr }} (-q)^{-\ell (w)} T_{w}
E^{\sigmahat }_{\mu },
\\
\label{e:F-minus}
F^{\sigma ,-}_{\rr ,\mu } & = \big( q^{-\ell (w_{0}^{\rr })} \prod
\nolimits _{\alpha \in R_{+}(\GL _{\rr })} (1- q \zz ^{-\alpha })\big)
F^{\sigma }_{\rr ,\mu } = \sum _{w\in S_{\rr }} (-q)^{-\ell (w)} T_{w}
F^{\sigmahat }_{\mu }.
\end{align}
Using Lemma~\ref{lem:Tw-on-EF-sigmahat}, we can also write
\begin{equation}\label{e:E-minus-alt}
E^{\sigma ,-}_{\rr ,\mu } = \sum _{w\in S_{\rr }} (-1)^{\ell (w)}
E^{\sigmahat }_{w(\mu )}.
\end{equation}
Defining $A_{\rr }^{q} = \sum _{w\in S_{\rr }} (-q)^{-\ell (w)} T_{w}$
as in Remark~\ref{rem:Hecke-antisymmetry}, we have
\begin{equation}\label{e:antisymmetrizer-trick}
\sum _{w\in S_{\rr }} (-q)^{-\ell (w)} T_{w} = A_{\rr }^{q} =A_{\rr
}^{q} \cdot (-1)^{\ell (w_{0}^{\rr })} T_{w_{0}^{\rr }}^{-1} =
q^{-\ell (w_{0}^{\rr })} \sum _{w\in S_{\rr }} (-q)^{\ell (w)}
\overline{T_{w}},
\end{equation}
and Lemma~\ref{lem:Tw-on-EF-sigmahat} then implies 
\begin{equation}\label{e:F-minus-alt}
F^{\sigma ,-}_{\rr ,\mu } = q^{- \ell (w_{0}^{\rr })}\sum _{w\in S_{\rr
}} (-q)^{\ell (w)} F^{\sigmahat }_{w(\mu )}.
\end{equation}

Macdonald's identity \cite[Theorem~(2.8)]{Macdonald72} for $\GL _{\rr
}$ gives
\begin{equation}\label{e:Macdonald-identity}
W_{\rr }(q) \defeq \sum _{w\in S_{\rr }} q^{\ell (w)} = \sum _{w\in
S_{\rr }} w\left(\prod \nolimits _{\alpha \in R_{+}(\GL _{\rr })}
\frac{1-q\, \zz ^{-\alpha }}{1 - \zz ^{- \alpha }} \right) =
\sigmabold _{\rr }\big( \prod \nolimits _{\alpha \in R_{+}(\GL _{\rr
})} (1-q\, \zz ^{-\alpha }) \big).
\end{equation}
Combining this with \eqref{e:1-G-in-Weyl},
if $f$ is $S_{\rr }$
invariant, we find
\begin{multline}\label{e:Mac-ident-var}
W_{\rr }(q)\, \langle 1_{\GL _{\rr }} \rangle\, f = \langle 1_{\GL
_{\rr }} \rangle\, \sigmabold _{\rr } \big(f \prod \nolimits _{\alpha
\in R_{+}(\GL _{\rr })} (1-q\, \zz ^{-\alpha }) \big)
\\
= \langle \zz ^{0} \rangle\, f \cdot \prod \nolimits _{\alpha \in
R_{+}(\GL _{\rr })} \big( (1 - \zz ^{\alpha })(1-q\, \zz ^{-\alpha })
\big).
\end{multline}
Now we calculate
\begin{equation}\label{e:inner-prod-proof}
\begin{aligned}
W_{\rr }(q) \big\langle\, & \zz ^{-\rho _{\rr }} E^{\sigma }_{\rr
,\lambda },\, \overline{\zz ^{-\rho _{\rr }} F^{\sigma }_{\rr ,\mu}}\,
\big\rangle ^{\rr }_{q}
\\
& = W_{\rr }(q)\, \langle 1_{\GL _{\rr }} \rangle\, E^{\sigma }_{\rr
,\lambda }\, \overline{F^{\sigma }_{\rr ,\mu }}\, \prod
\nolimits_{\alpha \in R_{+}\setminus R_{+}(\GL _{\rr }) } \frac{1 -
\zz ^{\alpha }}{1 - q^{-1} \zz ^{\alpha }}
\\
& = \langle \zz ^{0} \rangle \, E^{\sigma }_{\rr ,\lambda }\,
\overline{F^{\sigma }_{\rr ,\mu }} \, \prod \nolimits _{\alpha \in
R_{+}(\GL _{\rr })} \big( (1 - q^{-1} \zz ^{\alpha })(1-q\, \zz
^{-\alpha }) \big) \prod \nolimits _{\alpha \in R_{+}} \frac{1- \zz
^{\alpha }}{1 - q^{-1}\, \zz ^{\alpha }}
\\
& = \big\langle \, E^{\sigma }_{\rr ,\lambda }\, \prod \nolimits
_{\alpha \in R_{+}(\GL _{\rr })} (1-q\, \zz ^{-\alpha }),\,
\overline{F^{\sigma }_{\rr ,\mu }\, \prod \nolimits _{\alpha \in
R_{+}(\GL _{\rr })} (1-q\, \zz ^{-\alpha })}\, \big\rangle _{q}
\\
& = \langle \, q^{\ell (w_{0}^{\rr })} E^{\sigma ,-}_{\rr ,\lambda
},\, \overline{q^{\ell (w_{0}^{\rr })} F^{\sigma ,-}_{\rr ,\mu }}\,
\rangle _{q}
\\
& = \big\langle\, \sum _{w\in S_{\rr }} (-1)^{\ell (w)} E^{\sigmahat
}_{w(\lambda )},\, q^{\ell (w_{0}^{\rr })}\sum _{w\in S_{\rr }} (-q)^{
-\ell (w)} \overline{F^{\sigmahat }_{w(\mu )}} \, \big\rangle _{q},
\end{aligned}
\end{equation}
where we canceled $q^{\pm \ell (w_{0}^{\rr })}$ in the penultimate
line and used \eqref{e:E-minus-alt} and \eqref{e:F-minus-alt} to get
the last line.  By \cite[Proposition~4.3.2]{BHMPS-paths}, the
functions $E^{\sigmahat }_{\nu }$ and $\overline{F^{\sigmahat }_{\nu
}}$ are dual bases for $\langle -,- \rangle_{q}$.  Hence, the last line
in \eqref{e:inner-prod-proof} simplifies to $W_{\rr }(q)\, \delta
_{\lambda ,\mu }$, and the result follows.
\end{proof}

\subsection{LLT series}
\label{ss:LLT-series}

Generalizing \cite[Definition~4.4.1]{BHMPS-paths}, we now define LLT
series associated to $\GL _{l}$ and any Levi subgroup $\GL _{\rr }$.

\begin{defn}\label{def:LLT-series}
Given a composition $\rr = (r_{1},\ldots,r_{k})$ of $l$, a permutation
$\sigma \in S_{k}$, and weights $\alpha ,\beta \in X^{++}(\GL _{\rr
})$, the {\em LLT series} $\Lcal ^{\sigma }_{\rr ,\beta /\alpha }(\zz
;q)$ is the infinite formal linear combination of irreducible $\GL
_{l}$ characters with coefficients defined by
\begin{equation}\label{e:LLT-series}
\langle \chi _{\lambda } \rangle\, \Lcal ^{\sigma ^{-1}}_{\rr ,\beta
/\alpha }(\zz ;q^{-1}) = \langle E^{\sigma }_{\rr , \beta }(\zz ;q)
\rangle\, \chi _{\lambda }\cdot E^{\sigma }_{\rr , \alpha }(\zz ;q)
\end{equation}
in terms of the basis $\{E^{\sigma }_{\rr ,\mu } \mid \mu \in
X^{++}(\GL _{\rr }) \}$ of the space $\zz ^{\rho _{\rr }} \kk
[z_{1}^{\pm 1},\ldots,z_{l}^{\pm 1}]^{S_{\rr }}$ (note that this space
is independent of the choice of $\rho _{\rr }$ and closed under
multiplication by $\GL _{l}$ characters).
\end{defn}

\begin{remark}\label{rem:L-alternative-etc}
(i) The elements $E^{\sigma ,-}_{\rr ,\mu }(\zz ;q)$ in
\eqref{e:E-minus} form a basis of the space of Laurent polynomials
antisymmetric with respect to the action of the Hecke algebra $\Hcal
(S_{\rr })$.  In terms of this basis, an alternative formulation
equivalent to \eqref{e:LLT-series} is
\begin{equation}\label{e:LLT-series-alt}
\langle \chi _{\lambda } \rangle\, \Lcal ^{\sigma ^{-1}}_{\rr ,\beta
/\alpha }(\zz ;q^{-1}) = \langle E^{\sigma , -}_{\rr , \beta }(\zz ;q)
\rangle\, \chi _{\lambda }\cdot E^{\sigma , -}_{\rr , \alpha }(\zz
;q).
\end{equation}

(ii) By Remark~\ref{rem:coefficient-rings}, the right hand side of
\eqref{e:LLT-series} belongs to $\ZZ [q^{-1}]$.  The $q^{-1}$ on the
left hand side of \eqref{e:LLT-series} serves to give $\Lcal ^{\sigma
}_{\rr ,\beta /\alpha }(\zz ;q)$ coefficients in $\ZZ [q]$, while the
indexing with $\sigma ^{-1}$ instead of $\sigma $ allows us to
formulate the connection between LLT series and LLT polynomials more
naturally in \S \ref{ss:LLT-series-vs-pol}.
\end{remark}

The next proposition gives a formula for $\Lcal ^{\sigma }_{\rr ,\beta
/\alpha }(\zz ;q)$, generalizing
\cite[Proposition~4.4.2]{BHMPS-paths}.  To state it we need the
$q$-symmetrization operator
\begin{equation}\label{e:Hq-r}
\Hbold _{q}^{\rr }(f) = \sigmabold \Big( \frac{f(\zz )}{\prod _{\alpha
\in R_{+}\setminus R_{+}(\GL _{\rr })} (1 - q\, \zz ^{\alpha })} \Big).
\end{equation}
Here $f(z_{1},\ldots,z_{l})$ is a Laurent polynomial, and
\eqref{e:Hq-r} is to be interpreted as a formal infinite linear
combination of irreducible $\GL _{l}$ characters by expanding each
factor $(1- q\, \zz ^{\alpha }) = 1 + q\, \zz ^{\alpha } +\cdots $ as
a geometric series before applying the Weyl symmetrization operator
$\sigmabold $.

Although we won't use it, we mention that when $f$ is $S_{\rr }$
invariant, $\Hbold ^{\rr }_{q}$ is a $q$-analog of induction from $\GL
_{\rr }$ characters to $\GL _{l}$ characters.  When $f$ is a product
of Schur functions $\prod _{i}s_{\lambda _{(i)}}(Z_{i})$ in blocks of
variables $Z_{1} = z_{1},\ldots,z_{r_{1}}$, $Z_{2} =
z_{r_{1}+1},\ldots,z_{r_{1}+r_{2}}$, etc., $\Hbold _{q}^{\rr }(f)$ is
a $q$-analog of $\prod _{i}s_{\lambda _{(i)}}(\zz )$, whose Schur
expansion yields the generalized Kostka polynomials studied by
Shimozono, Weyman and Zabrocki in \cite{ShimWeym00,ShimZabr01}.

\begin{prop}\label{prop:L-via-Hq}
For $\rr $, $\sigma $, $\alpha $, $\beta $ as in
Definition~\ref{def:LLT-series}, we have
\begin{equation}\label{e:L-via-Hq}
\Lcal ^{\sigma }_{\rr ,\beta /\alpha }(\zz ;q) = \Hbold
^{w_{0}^{k}(\rr )}_{q} \big( w_{0}(F^{\sigma ^{-1}}_{\rr ,\beta }\,
\overline{E^{\sigma ^{-1}}_{\rr ,\alpha }}) \big).
\end{equation}
\end{prop}

\begin{proof}
By Proposition~\ref{prop:orthogonality},
\begin{equation}\label{e:L-via-Hq-proof-I}
\begin{aligned}
\langle \chi _{\lambda } \rangle\, \Lcal ^{\sigma }_{\rr ,\beta
/\alpha }(\zz ;q) & = \langle\, \zz ^{\rho _{\rr }} F^{\sigma
^{-1}}_{\rr ,\beta }(\overline{\zz }; q),\, \zz ^{-\rho _{\rr }}\chi
_{\lambda }\cdot E^{\sigma ^{-1} }_{\rr , \alpha }(\zz ;q^{-1})\,
\rangle ^{\rr }_{q^{-1}}
\\
& = \langle \zz ^{0} \rangle\, \chi _{\lambda }\, F^{\sigma
^{-1}}_{\rr ,\beta }(\overline{\zz }; q)\, E^{\sigma ^{-1} }_{\rr ,
\alpha }(\zz ;q^{-1}) \, \frac{\prod _{\alpha \in R_{+}}(1 - \zz
^{\alpha })}{\prod _{\alpha \in R_{+}\setminus R_{+}(\GL _{\rr })}(1 -
q \, \zz ^{\alpha })}.
\end{aligned}
\end{equation}
We can invert the variables $z_{i}$ and apply $w_{0}$ without changing
the constant term, so the above is equal to
\begin{equation}\label{e:L-via-Hq-proof-II}
\langle \zz ^{0} \rangle\,\overline{\chi _{\lambda }}
\,w_{0} \big( F^{\sigma ^{-1}}_{\rr ,\beta }(\zz ; q)\, \overline{ E^{\sigma
^{-1} }_{\rr , \alpha }(\zz ;q)} \big) \, \frac{\prod _{\alpha \in R_{+}}(1
- \zz ^{\alpha })}{\prod _{\alpha \in R_{+}\setminus R_{+}(\GL
_{w_{0}^{k}(\rr ) })}(1 - q\, \zz ^{\alpha })}.
\end{equation}
Using \eqref{e:1-G-in-Weyl} for $\GL _{l}$, this is the same as
\begin{equation}\label{e:L-via-Hq-proof-III}
\langle 1_{\GL _{l}} \rangle\,\overline{\chi _{\lambda }}\, \sigmabold
\bigg( \frac{w_{0} (F^{\sigma ^{-1}}_{\rr ,\beta }\, \overline{
E^{\sigma ^{-1} }_{\rr , \alpha }})}{\prod _{\alpha \in R_{+}\setminus
R_{+}(\GL _{w_{0}^{k}(\rr ) })}(1 - q\, \zz ^{\alpha })} \bigg) =
\langle \chi _{\lambda } \rangle \, \Hbold ^{w_{0}^{k}(\rr )}_{q}
\big( w_{0}(F^{\sigma ^{-1}}_{\rr ,\beta }\, \overline{E^{\sigma
^{-1}}_{\rr ,\alpha }}) \big).
\qedhere
\end{equation}
\end{proof}

\subsection{Relation between LLT series and LLT polynomials}
\label{ss:LLT-series-vs-pol}

Here we derive an identity relating the polynomial part $\Lcal
^{\sigma }_{\rr ,\beta /\alpha }(\zz ;q)_{\pol }$ of an LLT series to
an LLT polynomial $\Gcal _{\nubold }(X;q)$, extending the treatment of
the case $\rr = (1^{l})$ in \cite[\S 4.5]{BHMPS-paths}.

\begin{defn}\label{def:triple}
Let $\nubold = \betabold /\alphabold = (\beta _{(1)}/\alpha
_{(1)},\ldots,\beta _{(k)}/\alpha _{(k)})$ be a tuple of skew
diagrams, and let $\sigma \in S_{k}$ be a permutation.  A {\em $\sigma
$-triple} in $ \betabold /\alphabold$ is an ordered triple of boxes
$(a,b,c)$ such that
\begin{itemize}
\item [(i)] $b$ is a box of $\nu _{(i)} = \beta _{(i)}/\alpha _{(i)}$
for some $i$;
\item [(ii)] $a$ is either in or immediately left of a row of $\nu
_{(j)}$ and $c$ is either in or immediately right of the same row, for
some $j>i$; 
\item [(iii)] $a$ and $c$ are adjacent with $a$ left of $c$; and
\item [(iv)] $b$ has the same content as $c$ if $\sigma (i)<\sigma
(j)$, or the same content as $a$ if $\sigma (i)>\sigma (j)$.
\end{itemize}
More precisely, (ii) and (iii) mean that if $\alpha _{(j)} = (\alpha
_{1},\ldots,\alpha _{m})$ and $\beta _{(j)} = (\beta _{1},\ldots,\beta
_{m})$, then we have $a = (x, y)$ and $c = (x+1,y)$ for some $1\leq
y\leq m$ and $\alpha _{y}\leq x\leq \beta _{y}$.  In particular, $a$
and $c$ can be the boxes left and right of an empty row with $\beta
_{y} = \alpha _{y}$.  The set of triples thus depends on the
presentation of $\nubold $ as $\betabold /\alphabold $ and not just on
the set of boxes in $\nubold $.
\end{defn}

\begin{figure}
\[
\begin{tikzpicture}[scale=.5]
\draw (6,1) grid (9,2);
\draw (4,2) grid (6,4);
\draw (6,1) -- (6,2);
\node [above left] at (4,4) {$\nu _{(i)}$};
\draw (8,7) -- (8,5) -- (10,5) -- (10,7) -- (7,7) -- (7,6) -- (10,6);
\draw (9,7) -- (9,5); \draw (10,4) -- (10,5);
\draw (10,3) rectangle (11,4);
\node [above left] at (7,7) {$\nu _{(j)}$};
\node at (7.5,6.5) {$a_{1}$};
\node at (5.5,3.5) {$b_{1}$};
\node at (8.5,6.5) {$c_{1}$};
\node at (8.5,5.5) {$a_{2}$};
\node at (5.5,2.5) {$b_{2}$};
\node at (9.5,5.5) {$c_{2}$};
\node at (9.5,4.5) {$a_{3}$};
\node at (7.5,1.5) {$b_{3}$};
\node at (10.5,4.5) {$c_{3}$};
\node at (10.5,3.5) {$a_{4}$};
\node at (8.5,1.5) {$b_{4}$};
\node at (11.5,3.5) {$c_{4}$};

% content diagonals
\draw [dashed,gray] (3,3) -- (8,8);
\draw [dashed,gray] (3,2) -- (4,3); \draw [dashed,gray] (5,4) -- (7,6); 
\draw [dashed,gray] (8,7) -- (9,8);
\draw [dashed,gray] (3,1) -- (4,2); \draw [dashed,gray] (6,4) -- (8,6); 
\draw [dashed,gray] (9,7) -- (10,8); 
\draw [dashed,gray] (3,0) -- (5,2); 
\draw [dashed,gray] (6,3) -- (8,5); \draw [dashed,gray] (10,7) -- (11,8); 
\draw [dashed,gray] (4,0) -- (9,5); \draw [dashed,gray] (10,6) -- (12,8);
\draw [dashed,gray] (5,0) -- (6,1); \draw [dashed,gray] (7,2) -- (9,4);
\draw [dashed,gray] (10,5) -- (13,8);
\draw [dashed,gray] (6,0) -- (7,1); \draw [dashed,gray] (8,2) -- (10,4);
\draw [dashed,gray] (11,5) -- (13,7);
\draw [dashed,gray] (7,0) -- (8,1); \draw [dashed,gray] (9,2) -- (10,3);
\draw [dashed,gray] (11,4) -- (13,6);
\draw [dashed,gray] (8,0) -- (11,3); \draw [dashed,gray] (12,4) -- (13,5);
\end{tikzpicture}
\]
\caption{\label{fig:triples}%
Examples of $\sigma $-triples.  The dashed lines show boxes of equal
content in $\nu _{(i)}$ and $\nu _{(j)}$, with $i<j$.  If $\sigma
(i)<\sigma (j)$, then $(a_{1},b_{1},c_{1})$ and $(a_{3},b_{3},c_{3})$
are $\sigma $-triples. If $\sigma (i)>\sigma (j)$, then
$(a_{2},b_{2},c_{2})$ and $(a_{4},b_{4},c_{4})$ are $\sigma $-triples.
Triples $(a_{3},b_{3},c_{3})$ and $(a_{4},b_{4},c_{4})$ illustrate the
point that $a$ and/or $c$ may be just outside a (possibly empty) row
of $\nu _{(j)}$.}
\end{figure}
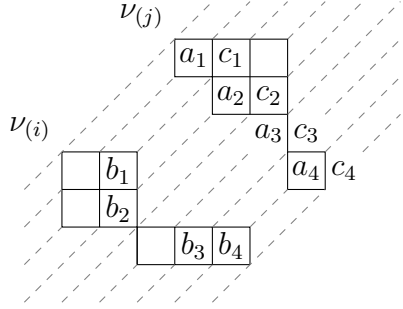

Strictly speaking, the indices $i$ and $j$ are part of the data of a
triple, in keeping with our understanding that the set of boxes of
$\nubold $ is the disjoint union of the sets of boxes of the $\nu
_{(i)}$.  Figure~\ref{fig:triples} illustrates the definition.

\begin{defn}\label{def:N-beta-alpha}
(i) Given a tuple of skew diagrams $\betabold /\alphabold = (\beta
_{(1)}/\alpha _{(1)},\ldots,\beta _{(k)}/\alpha _{(k)})$ and $\sigma
\in S_{k}$, an {\em increasing $\sigma $-triple} in a negative tableau
$T\in \SSYT _{-}(\betabold /\alphabold )$ is a $\sigma $-triple
$(a,b,c)$ such that $T(a)<T(b)<T(c)$, with the convention $T(a) =
-\infty $, $T(c)= \infty $ if $a$ or $c$ is not a box of $\betabold
/\alphabold $.

(ii) We define the generating function
\begin{equation}\label{e:N-beta-alpha}
N^{\sigma }_{\betabold /\alphabold }(X;q) = \sum _{T\in \SSYT _{-}
(\betabold /\alphabold )} q^{h_{\sigma }(T)} \xx ^{T},
\end{equation}
where $h_{\sigma }(T)$ is the number of increasing $\sigma $-triples
in $T$.
\end{defn}

Although it is not obvious a priori, the next proposition implies that
$N^{\sigma }_{\betabold /\alphabold }(X;q)$ is symmetric.

\begin{prop}\label{prop:N-beta-alpha}
We have the identity
\begin{equation}\label{e:N-vs-G}
N^{\sigma }_{\betabold /\alphabold }(X;q) = q^{h_{\sigma }(\betabold
/\alphabold )} \omega \, \Gcal _{\sigma (\betabold /\alphabold
)}(X;q^{-1}),
\end{equation}
where $h_{\sigma }(\betabold
/\alphabold )$ is the number of $\sigma $-triples in $\betabold
/\alphabold $.
\end{prop}

\begin{proof}
Let $\nubold $ be the tuple of skew diagrams given by $\betabold/\alphabold$ 
and note that $\sigma (\nubold
)$ is the rearrangement of $\nubold $ with $\sigma (\nubold )_{(\sigma
(i))} = \nu _{(i)}$.  For $a\in \nu _{(i)}$, let $\sigma (a)$ denote
the corresponding box of $\sigma (\nubold )_{(\sigma (i))}$.  We also
use this notation for boxes adjacent to $\nu _{(i)}$, which may
occur in a $\sigma $-triple. 

For a tableau $T\in \SSYT _{-}(\nubold )$ let $\sigma (T)$ denote its
image in $\SSYT _{-}(\sigma (\nubold ))$, defined by $\sigma
(T)(\sigma (a))= T(a)$.

By Proposition~\ref{prop:omega-G-nu}, the right hand side of
\eqref{e:N-vs-G} can be written
\begin{equation}\label{e:G-side}
\sum _{T\in \SSYT _{-}(\betabold /\alphabold )} q^{h_{\sigma
}(\betabold /\alphabold ) - \inv (\sigma (T))} \xx ^{\sigma (T)}.
\end{equation}
Since $\xx ^{\sigma (T)} = \xx ^{T}$, \eqref{e:N-vs-G} will follow if we
show that 
\begin{equation}\label{e:h-sigma(T)}
h_{\sigma }(T) = h_{\sigma }(\betabold /\alphabold ) - \inv
(\sigma (T)).
\end{equation}

Consider the image $(a',b',c') = (\sigma (a),\sigma (b),\sigma (c))$
of a $\sigma $-triple $(a,b,c)$ in $\betabold /\alphabold$.  The
definition of $\sigma $-triple implies that $(a',b')$ is an attacking
pair in $\sigma (\nubold )$ if both boxes $a$ and $b$ are in $\nubold
$, and similarly for $(b',c')$.  One also sees that every attacking
pair in $\sigma (\nubold )$ belongs in this way to the image of a
unique triple.  Since $a,c$ are in the same row, we have $T(a)<T(c)$
for every negative tableau $T$.  This holds even if $a$ or $c$ is not
in $\nubold $, by the convention that $T(a) = -\infty $, $T(c) =
\infty $ in these cases.  Hence, at most one of the pairs $(a',b')$,
$(b',c')$ is an attacking inversion in $\sigma (T)$, since we would
have $T(a)\geq T(b)\geq T(c)$ if they both were.  Moreover, $(a,b,c)$
is an increasing $\sigma $-triple if and only if neither $(a',b')$ nor
$(b',c')$ is an attacking inversion.  Hence, the number of increasing
$\sigma $-triples is the total number of $\sigma $-triples minus the
number of attacking pairs in $\sigma (T)$.
\end{proof}

Let $\rr = (r_{1},\ldots,r_{k})$ be a composition of $l$, and let
$\alpha ,\beta \in X^{++}(\GL _{\rr })$ be dominant regular weights
for $\GL _{\rr }$ such that $\alpha _{i}\leq \beta _{i}$ for all
$1\leq i\leq l$.  To these data we associate a tuple of skew diagrams
$\betabold /\alphabold = (\beta
_{(1)}/\alpha _{(1)},\ldots,\beta _{(k)}/\alpha _{(k)})$ by defining
\begin{equation}\label{e:nu(r,beta/alpha)}
\begin{aligned}
(\alpha  _{(i)})_{j} & = \alpha  _{r_{1}+\cdots +r_{i-1}+j} + j,\\
(\beta _{(i)})_{j} & = \beta _{r_{1}+\cdots +r_{i-1}+j} + j
\end{aligned}
\end{equation}
for $1 \leq i \leq k$ and $1 \leq j \leq r_{i}$.
In other words,
$(\alpha _{(1)}|\cdots |\alpha _{(k)}) = \alpha - \rho _{\rr }$ and
$(\beta _{(1)}|\cdots |\beta _{(k)}) =\beta -\rho _{\rr }$, where
$(\cdot |\cdots |\cdot )$ denotes concatenation and $\rho _{\rr } =
-((1,\ldots,r_{1})|\cdots |(1,\ldots,r_{k}))$.

This construction has the combinatorially natural feature that the
contents of the boxes in the $j$-th row of $\beta _{(i)}/\alpha
_{(i)}$ are $\alpha _{m}+1,\ldots,\beta _{m}$, where $m = r_{1}+\cdots
+r_{i-1}+j$ is the index corresponding to the $j$-th position in the
$i$-th block of the partition of $[l]$ into intervals of lengths
$r_{i}$.

\begin{thm}\label{thm:L-vs-G}
Given a composition $\rr = (r_{1},\ldots,r_{k})$ of $l$, $\sigma \in
S_{k}$, and weights $\alpha ,\beta \in X^{++}(\GL _{\rr })$, we have
\begin{equation}\label{e:L-vs-G}
\Lcal ^{\sigma }_{\rr ,\beta /\alpha }(\zz ;q)_{\pol } = 
\begin{cases}
q^{h_{\sigma }(\betabold /\alphabold )} \Gcal _{\sigma
(\betabold/\alphabold )}(z_{1},\ldots,z_{l};q^{-1})
  & \text{if $\alpha _{i}\leq \beta _{i}$ for all $1\leq i\leq l$},\\
0 & \text{otherwise},
\end{cases}
\end{equation}
where $\betabold /\alphabold $ is the tuple associated to $\rr $,
$\alpha $, $\beta $ by the construction above.
\end{thm}

The proof will be based on Proposition~\ref{prop:N-beta-alpha} and the
following lemma.

\begin{lemma}\label{lem:em-on-E-r}
For $\alpha ,\beta \in X^{++}(\GL _{\rr })$ and $\sigma \in S_{k}$,
we have
\begin{equation*}
\langle E^{\sigma , -}_{\rr , \beta }(\zz ;q) \rangle\, e_{m}(\zz )\,
E^{\sigma , -}_{\rr , \alpha }(\zz ;q) =
\begin{cases}
q^{-|{\Inv (\beta +\epsilon\,
\tau })\setminus \Inv (\alpha +\epsilon \, \tau )|} & \text{if $\beta
- \alpha = \varepsilon _{I}$ for $I\subseteq [l]$, $|I| = m$},\\
0 & otherwise,
\end{cases}
\end{equation*}
where $e_{m}$ is an elementary symmetric function, $\tau =\sigmahat \,
^{-1}$, and $\varepsilon _{I}$ is the $0$-$1$ vector with $1$'s in
positions $i\in I$.
\end{lemma}

\begin{proof}
Using \eqref{e:E-minus-alt}, the coefficient 
$\langle E^{\sigma , -}_{\rr , \beta }\rangle f$ of
$ E^{\sigma , -}_{\rr , \beta }$
in any $S_{\rr }$-antisymmetric function $f$ is equal to
$\langle E^{\sigmahat }_{\beta } \rangle f$.
Applying this with $f = e_{m}\,
E^{\sigma , -}_{\rr , \alpha } $,
\begin{equation}\label{e:coefficient-reduction}
\langle E^{\sigma , -}_{\rr , \beta } \rangle\, e_{m}\,
E^{\sigma , -}_{\rr , \alpha }=
\langle E^{\sigmahat }_{\beta } \rangle\, e_{m}\,
E^{\sigma , -}_{\rr , \alpha } = \sum _{w\in S_{\rr }}
(-1)^{\ell (w)} \langle E^{\sigmahat }_{\beta } \rangle\,
e_{m}\, E^{\sigmahat }_{w(\alpha )}.
\end{equation}
By \cite[Lemma~4.5.1]{BHMPS-paths} the coefficient $ \langle
E^{\sigmahat }_{\beta }\rangle\, e_{m}\, E^{\sigmahat
}_{w(\alpha )}$ vanishes unless $\beta = w(\alpha
)+\varepsilon _{I}$ for some $I$.  Since $w\in S_{\rr }$ and $\alpha
\in X^{++}(\GL _{\rr })$, $w(\alpha )$ is not dominant for $\GL _{\rr
}$ if $w\not =1$.  In that case there is an index $i$ such that
$s_{i}\in S_{\rr }$ and $w(\alpha )_{i}<w(\alpha )_{i+1}$.  Since
$\beta \in X^{++}(\GL _{\rr })$ we then have $\beta _{i}-w(\alpha
)_{i} - (\beta _{i+1} - w(\alpha )_{i+1})\geq 2$, and therefore $\beta
-w(\alpha )$ is not of the form $\varepsilon _{I}$.  This shows that
the terms for $w\not =1$ on the right hand side of
\eqref{e:coefficient-reduction} vanish, leaving
\begin{equation}\label{e:coefficient-reduction-bis}
\langle E^{\sigma , -}_{\rr , \beta } \rangle\, e_{m}\,
E^{\sigma , -}_{\rr , \alpha } = 
 \langle E^{\sigmahat }_{\beta } \rangle\,
e_{m}\, E^{\sigmahat }_{\alpha }.
\end{equation}
The lemma now follows from \cite[Lemma~4.5.1]{BHMPS-paths}.
\end{proof}

\begin{proof}[Proof of Theorem~\ref{thm:L-vs-G}]

Let $L^{\sigma }_{\rr ,\beta /\alpha }(X;q)$ 
be the unique linear
combination of Schur functions $s_{\lambda }(X)$ with $\ell (\lambda
)\leq l$ that specializes in $l$ variables $z_{1},\ldots,z_{l}$ to
$L^{\sigma }_{\rr ,\beta /\alpha }(\zz ;q) = \Lcal ^{\sigma }_{\rr
,\beta /\alpha }(\zz ;q)_{\pol }$.  We will prove that
\begin{equation}\label{e:L-vs-G-bis}
L ^{\sigma }_{\rr ,\beta /\alpha }(X ;q) = 
\begin{cases}
q^{h_{\sigma }(\betabold /\alphabold )} \Gcal _{\sigma
(\betabold/\alphabold )}(X;q^{-1})
  & \text{if $\alpha _{i}\leq \beta _{i}$ for all $1\leq i\leq l$},\\
0 & \text{otherwise}.
\end{cases}
\end{equation}
Clearly this implies \eqref{e:L-vs-G} (actually, \eqref{e:L-vs-G} and
\eqref{e:L-vs-G-bis} are equivalent, by
Corollary~\ref{cor:G-nu-width}).

By \eqref{e:LLT-series-alt}, we have
\begin{equation}\label{e:Hall-s-lambda-L}
\langle \, s_{\lambda },\, L^{\sigma }_{\rr ,\beta /\alpha
}(X;q^{-1})\, \rangle = \langle E^{\sigma ^{-1},-}_{\rr ,\beta }
\rangle\, s_{\lambda }(\zz )\, E^{\sigma ^{-1}, -}_{\rr ,\alpha },
\end{equation}
where $\langle -,- \rangle$ is the Hall inner product on symmetric
functions.  Note that \eqref{e:Hall-s-lambda-L} holds even if $\ell
(\lambda )>l$, since both sides are zero in that case.  By linearity,
\eqref{e:Hall-s-lambda-L} therefore holds with any symmetric function $f$ in
place of $s_{\lambda }$.  In particular, taking $f = e_{\mu } = e_{\mu
_{1}}\cdots e_{\mu _{n}}$, we have
\begin{equation}\label{e:Hall-e-mu-L}
\langle \, e_{\mu },\, L^{\sigma }_{\rr ,\beta /\alpha }(X;q^{-1})\,
\rangle = \langle E^{\sigma ^{-1},-}_{\rr ,\beta } \rangle\, e_{\mu
}(\zz )\, E^{\sigma ^{-1}, -}_{\rr ,\alpha }.
\end{equation}
Evaluating this last expression by using Lemma~\ref{lem:em-on-E-r} to
multiply by $e_{\mu _{1}}$ through $e_{\mu _{n}}$ in succession gives
a sum over chains of weights
\begin{equation}\label{e:weight-chain}
\alpha = \alpha ^{(0)}, \alpha ^{(1)}, \ldots, \alpha ^{(n)} =
\beta\qquad (\alpha ^{(i)}\in X^{++}(\GL _{\rr }))
\end{equation}
such that for each $i = 1,\ldots,n$ we have $\alpha ^{(i)} = \alpha
^{(i-1)} + \varepsilon _{I_{i}}$ for some index set $I_{i}\subseteq
[l]$ of size $|I_{i}| = \mu _{i}$. 
In particular, this shows that
$L^{\sigma }_{\rr ,\beta /\alpha }(X;q) = 0$ if the condition $\alpha
_{i}\leq \beta _{i}$ for all $i$ does not hold, so we assume from now
on that it does.

To the weights in \eqref{e:weight-chain} we can now associate a chain
of tuples of skew diagrams
\begin{equation}\label{e:diagram-chain}
\varnothing = \alphabold ^{(0)}/\alphabold \subseteq \alphabold
^{(1)}/\alphabold \subseteq \cdots \subseteq \alphabold
^{(n)}/\alphabold = \betabold /\alphabold
\end{equation}
by the construction in \eqref{e:nu(r,beta/alpha)}.  The condition on
the $\alpha ^{(i)}$ means that each $\alphabold ^{(i)}/\alphabold
^{(i-1)}$ is a tuple of vertical strips of size $|\alphabold
^{(i)}/\alphabold ^{(i-1)}| = \mu _{i}$.  In other words, $\alphabold
^{(i)}/\alphabold^{(i-1)}$ is the set of boxes $a$ with $T(a) = \overline{i}$
in a negative tableau $T\in \SSYT _{-}(\betabold /\alphabold )$ of
weight $\xx ^{T} = \xx ^{\mu }$. 

From Lemma~\ref{lem:em-on-E-r}, the term in \eqref{e:Hall-e-mu-L}
corresponding to $T$, or to the weight sequence
\eqref{e:weight-chain}, is $q^{-h(T)}$, where
\begin{equation}\label{e:L-Hall-term}
 h(T) = \sum _{i=1}^{n} |{\Inv (\alpha ^{(i)}+\epsilon \, \tau )
\setminus \Inv (\alpha ^{(i-1)} + \epsilon \, \tau )}|,
\end{equation}
with $\tau = (\widehat{\sigma ^{-1}}) ^{-1}$. 
We claim that $h(T) = h_{\sigma }(T)$ is the number of increasing $\sigma
$-triples in $T$.
Granting the claim, we then have
\begin{equation}\label{e:L-coefficient}
\langle e_{\mu },\, L^{\sigma }_{\rr ,\beta /\alpha }(X;q) \rangle = 
\sum _{\substack{T\in \SSYT _{-}(\betabold /\alphabold )\\
\xx ^{T} = \xx ^{\mu }}} q^{h_{\sigma }(T)}.
\end{equation}
By definition, the sum on the right is the coefficient $\langle \xx
^{\mu } \rangle\, N^{\sigma }_{\betabold /\alphabold } = \langle
e_{\mu },\, \omega \, N^{\sigma }_{\betabold /\alphabold } \rangle$.
Using Proposition~\ref{prop:N-beta-alpha}, this implies
\begin{equation}\label{e:L-vs-G-wrapup}
L^{\sigma }_{\rr ,\beta /\alpha }(X;q) = \omega \, N^{\sigma
}_{\betabold /\alphabold }(X;q) = q^{h_{\sigma }(\betabold /\alphabold
)} \Gcal _{\sigma (\betabold /\alphabold )}(X;q^{-1}).
\end{equation}

It remains only to verify that $h(T) = h_{\sigma }(T)$.  Let $[l] =
J_{1}\coprod \cdots \coprod J_{k}$ be the partition of $[l]$ into
intervals of lengths $|J_{j}| = r_{j}$.  The weights $\alpha ^{(i)}$
are strictly decreasing on each block $J_{j}$, so $|{\Inv (\alpha
^{(i)}+\epsilon \, \tau ) \setminus \Inv (\alpha ^{(i-1)} + \epsilon
\, \tau )}|$ only counts inversions between distinct blocks.

Now, $\tau = (\widehat{\sigma ^{-1}}) ^{-1}$ carries the blocks
$J_{j}$ to intervals of lengths $\sigma (\rr )$ in the order given by
$\sigma $; in other words, for $s\in J_{j}$, $s'\in J_{j'}$ with
$j<j'$, we have $\tau (s)<\tau (s')$ if and only if $\sigma (j)<\sigma
(j')$.  Thus, if $\alpha ^{(i)}_{s} = \alpha ^{(i)}_{s'}$, we have
$(s,s')\in \Inv (\alpha ^{(i)} + \epsilon \, \tau )$ if and only if
$\sigma (j) > \sigma (j')$.

By construction, $\alpha ^{(i)}_{s}$ is the content of the last box in
the row of $\alphabold ^{(i)}/\alphabold$ corresponding to the index
$s$, or of the box immediately left of an empty row.  For $s\in
J_{j}$, $s'\in J_{j'}$ with $j<j'$, it follows that $(s,s')\in \Inv
(\alpha ^{(i)}+\epsilon \, \tau ) \setminus \Inv (\alpha ^{(i-1)} +
\epsilon \, \tau )$ if and only if $\alphabold ^{(i)}/\alphabold
^{(i-1)}$ has a box $b$ in the row corresponding to $s$, and one of
the following two conditions holds, where $a$ is the last box in the
row of $\alphabold ^{(i-1)}/\alphabold $ corresponding to $s'$, or the
box immediately to the left if this row is empty:
\begin{itemize}
\item [(i)]
$\sigma (j)<\sigma (j')$ and $c(b) = c(a)+1$ and the box $c$ with
content $c(a)+1 = c(b)$ in the same row as $a$ is not in $\alphabold
^{(i)}/\alphabold $; or

\item [(ii)] $\sigma (j) > \sigma (j')$ and $c(b) = c(a)$ and the box
$c$ with content $c(a)+1 = c(b)+1$ in the same row as $a$ is not in
$\alphabold ^{(i)}/\alphabold $.
\end{itemize}
These conditions are equivalent to $(a,b,c)$ forming an increasing
$\sigma $-triple in $T$ with $T(b) = \overline{i}$.
Since a triple is determined by the box $b$ and the index of the row
containing $a$ and $c$, we see that $|{\Inv (\alpha ^{(i)}+\epsilon \,
\tau ) \setminus \Inv (\alpha ^{(i-1)} + \epsilon \, \tau )}|$ counts
$\sigma $-triples in $T$ such that $T(b) = \overline{i}$. Summing over
$i$ yields $h(T) = h_{\sigma }(T)$, as claimed.
\end{proof}

\section{Cauchy identity and winding permutations}
\label{s:Cauchy-winding}

As in \cite{BHMPS-paths}, the infinite series form of our main theorem
will follow by combining a Cauchy identity for Hall-Littlewood
polynomials---semi-symmetric Hall-Littlewood polynomials $E^{\sigma
}_{\rr ,\mu }$, $F^{\sigma }_{\sS ,\mu }$, in this case---with an
identity that allows us to change the `twist' $\sigma $ when the latter
has a special form.  In this section we establish the two identities
that we need.

\subsection{Cauchy identity}
\label{ss:Cauchy}

Our next theorem generalizes the Cauchy identity for non-symmetric
Hall-Littlewood polynomials \cite[Theorem~5.1.1]{BHMPS-paths}.  A new
feature that appears in the semi-symmetric case is that different
compositions $\rr $, $\sS $ may govern the blocks of variables in the
functions $E^{\sigma }_{\rr ,\mu }(\xx ,q)$ and $F^{\sigma }_{\sS ,\mu
}(\yy ,q)$ that play a role in the identity, subject to some
conditions which we now define.

\begin{defn}\label{def:almost-incr-decr}
Given $\sigma \in S_{k}$, a sequence $(m_{1},\ldots,m_{k})\in \ZZ
^{k}$ is {\em $\sigma $-almost decreasing} if
\begin{equation}\label{e:almost-decr}
m_{i}\geq m_{j} - \chi (\sigma ^{-1}(i) > \sigma ^{-1}(j)) \quad
\text{for all $i<j$},
\end{equation}
and {\em $\sigma $-almost increasing} if
\begin{equation}\label{e:almost-incr}
m_{i}\leq m_{j} + \chi (\sigma ^{-1}(i)<\sigma ^{-1}(j)) \quad
\text{for all $i<j$}.
\end{equation}
\end{defn}

If $\rr = (r_{1},\ldots, r_{k})$ is a strict composition, then any
choice of $\rho _{\rr }$ satisfying \eqref{e:rho-r} determines two
sequences $(M_{1}, \ldots, M_{k})$ and $(m_{1}, \ldots, m_{k})$ such
that $\rho _{\rr }$ is the concatenation of blocks $(M_{i}, M_{i}-1,
\ldots,m_{i})$ of length $r _{i}$.  We extend this to weak
compositions as follows.

\begin{convention}\label{conv:block-max-min}
Let $\rr = (r_{1},\ldots, r_{k})$ be a weak composition.  Whenever we
choose $\rho _{\rr }$ satisfying \eqref{e:rho-r}, we also choose
sequences of integers $(M_{1}, \ldots, M_{k})$ and $(m_{1}, \ldots,
m_{k})$ such that for $r_{i} > 0$, the corresponding block of $\rho
_{\rr }$ is $(M_{i}, M_{i}-1, \ldots,m_{i})$, and for $r_{i} = 0$ we
have $M_{i}<m_{i}$.  We refer to $M_{i}$ and $m_{i}$ as the {\em block
maxima and minima of $\rho _{\rr }$}, including any artificial maxima
and minima $M_{i}<m_{i}$ that we may have ascribed to empty blocks.
\end{convention}

\begin{thm}\label{thm:Cauchy}
Suppose we are given weak compositions $\rr = (r_{1},\ldots, r_{k})$,
$\sS = (s_{1},\ldots,s_{k})$, a permutation $\sigma \in S_{k}$, and a
choice of $\rho _{\rr }$, $\rho _{\sS }$ with associated block maxima
and minima, in keeping with Convention~\ref{conv:block-max-min}.
Assume that $\rho _{\rr }$ and $\rho _{\sS }$ have the same block
maxima $M_{i}$, and let $m_{i}$ and $n_{i}$ be their respective block
minima.  Assume further that $(m_{1},\ldots,m_{k})$ is $\sigma
$-almost decreasing and $(n_{1},\ldots,n_{k})$ is $\sigma $-almost
increasing.

Then, using notation explained below, we have the identity
\begin{equation}\label{e:Cauchy}
\frac{\prod _{i<j} \Omega [-q\, t\, X_{i} Y_{j}]}{\prod _{i\leq j}
\Omega [- t\, X_{i} Y_{j}]} = \sum _{\lambda } t^{|\lambda |}\, \xx
^{-\rho _{\rr }} E^{\sigma }_{\rr ,\lambda +\rho _{\rr }}(\xx ;q^{-1})
\, \yy ^{-\rho _{\sS}} F^{\sigma }_{\sS ,\lambda +\rho _{\sS }}(\yy
;q).
\end{equation}
The variables on the right hand side are $\xx =x_{1},\ldots,x_{r}$,
$\yy =y_{1},\ldots,y_{s}$, where $r = |\rr | = \sum _{i}r_{i}$ and $s
= |\sS | = \sum _{i}s_{i}$.  The index $\lambda $ ranges over tuples
of partitions $(\lambda _{(1)},\ldots,\lambda _{(k)})$ such that $\ell
(\lambda _{(i)})\leq \min (r_{i},s_{i})$, with $|\lambda | = \sum
_{i}|\lambda _{(i)}|$ denoting the sum of all the parts.  In the
expression $\lambda +\rho _{\rr }$, we interpret $\lambda $ as a
weight in $X^{+}(\GL _{\rr })$ by padding $\lambda _{(i)}$ with zeroes
to length $r_{i}$ and concatenating; $\lambda +\rho _{\sS }$ is
interpreted in a similar way.

The expressions involving $\Omega $ on the left hand side are as
defined in \S \ref{ss:symmetric-functions}, with plethystic alphabets
$X_{1},\ldots,X_{k}$ constructed from blocks
of length $r_{i}$ among the variables $\xx $, by the rule
\begin{equation}\label{e:Xi}
X_{i} = x_{i,1}+\cdots +x_{i,r_{i}}, \quad \text{where $x_{i,j} =
x_{r_{1}+\cdots + r_{i-1}+j}$}.
\end{equation}
The alphabets $Y_{1},\ldots,Y_{k}$ are similarly constructed from
blocks of length $s_{i}$ among the $\yy $ variables.
\end{thm}

\begin{remark}\label{rem:LHS}
(i) The left hand side of \eqref{e:Cauchy} expands to
\begin{equation}\label{e:Cauchy-expanded}
\frac{\prod _{i<j} \Omega [-q\, t\, X_{i} Y_{j}]}{\prod _{i\leq j}
\Omega [- t\, X_{i} Y_{j}]} = \frac{\prod _{i< j} \prod _{a=1}^{r_{i}}
\prod _{b=1}^{s_{j}} (1 - q\, t\, x_{i,a}y_{j,b}) }{ \prod _{i\leq j}
\prod _{a=1}^{r_{i}} \prod _{b=1}^{s_{j}} (1 - t\, x_{i,a}y_{j,b}) }.
\end{equation}

(ii) If $\rr =\sS =(1^{l})$ and $\rho _{\rr } = \rho _{\sS }$ is
constant, the theorem reduces to \cite[Theorem~5.1.1]{BHMPS-paths}. 
A little more generally, if $\rr =\sS =(1^{l})$, the hypotheses on the
block maxima and minima are satisfied when $\rho _{\rr } = \rho _{\sS
}$ has the form $\sigma (1,1,\ldots,1,0,0,\ldots,0) +
\text{(constant)}$.

(iii) If $k=1$, so $\rr =(r)$, $\sS =(s)$ and $\sigma =1\in S_{1}$, we
have $\xx ^{-\rho _{\rr }} E^{\sigma }_{\rr ,\lambda +\rho _{\rr
}}(\xx ;q^{-1}) = s_{\lambda }(x_{1},\ldots,x_{r})$ and $\yy ^{-\rho
_{\sS }} F^{\sigma }_{\sS ,\lambda +\rho _{\sS }}(\yy ;q) = s_{\lambda
}(y_{1},\ldots,y_{s})$ by Example~\ref{ex:special-semi-syms} (ii) and
\eqref{e:F-vs-E-bar}.
In this case the theorem reduces to the classical Cauchy identity for
Schur functions.

(iv) Adding a constant vector to $\rho _{\rr }$ or $\rho _{\sS}$
does not change the conclusion, so the hypothesis that $\rho _{\rr }$
and $\rho _{\sS}$ have the same block maxima can be weakened to having
block maxima that differ by a constant.  The hypotheses with this
weakening seem to be essentially as general as possible.

(v) The artificial maxima and minima ascribed to any zero-length
blocks in $\rho _{\rr }$ or $\rho _{\sS }$ have no effect on the
identity \eqref{e:Cauchy}.  Nevertheless, to conclude that the
identity holds, we require that such maxima and minima can be chosen
satisfying the hypotheses of the theorem.

\end{remark}

Before proving Theorem~\ref{thm:Cauchy}, we develop a series of
lemmas.

\begin{lemma}\label{lem:m<n-almost-decr}
Given $\sigma \in S_{l}$, let $m, n\in \ZZ ^{l}$ be $\sigma $-almost
decreasing sequences such that $m\leq n$ coordinate-wise, and $m\not
=n$.  Then there is an index $k\in [l]$ such that $m_{k}<n_{k}$,
and $m+\varepsilon _{k}$ is $\sigma $-almost decreasing, where
$\varepsilon _{k}$ is the $k$-th unit coordinate vector.
\end{lemma}

\begin{proof}
In fact, let $k$ be the smallest index such that $m_{k}<n_{k}$.  In
positions $i<j$, condition \eqref{e:almost-decr} for $m+\varepsilon
_{k}$ is the same as the condition on $m$ if $k\not \in \{i,j \}$, and
is weaker than the condition on $m$ if $i=k$.  If $j = k$, then since
$m_{i} = n_{i}$ for $i<k$ by assumption, the condition on
$m+\varepsilon _{k}$ in positions $i<j$ is weaker than the condition
on $n$.
\end{proof}

\begin{lemma}\label{lem:sk-almost-decr}
Given $\sigma \in S_{l}$ and $1\leq k <l$, if $m =
(m_{1},\ldots,m_{l})$ is $\sigma $-almost decreasing, then the
sequence
\begin{equation}\label{e:sk-almost-decr}
m' = (m_{1},\ldots,m_{k-1},m'_{k},m'_{k+1},m_{k+2},\ldots,m_{l})
\end{equation}
is $(s_{k}\sigma )$-almost decreasing, where
\begin{equation}\label{e:sk-almost-decr-detail}
\begin{aligned}
m'_{k} & = \max \big(m_{k}-\chi (\sigma ^{-1}(k+1) > \sigma ^{-1}(k)),
\, m_{k+1} \big)\\
m'_{k+1} & = m_{k}.
\end{aligned}
\end{equation}
\end{lemma}

\begin{proof}
We have $m'_{i} = m_{s_{k}(i)}$ for $i\not =k$, so if $i<j$ and $k\not
\in \{i,j \}$, then $s_{k}(i) < s_{k}(j)$ and the condition
$m'_{i}\geq m'_{j} - {\chi ((s_{k}\sigma )^{-1}(i)>(s_{k}\sigma
)^{-1}(j))}$ becomes $m_{s_{k}(i)}\geq m_{s_{k}(j)} - {\chi (\sigma
^{-1}s_{k}(i)>\sigma ^{-1}s_{k}(j))}$, which holds by hypothesis.

For $i = k$ and $j>k+1$, the condition becomes $m'_{k}\geq m_{j} -
{\chi (\sigma ^{-1}(k+1) > \sigma ^{-1}(j))}$, which follows from
$m_{k+1}\geq m_{j} - {\chi (\sigma ^{-1}(k+1) > \sigma ^{-1}(j))}$ and
$m'_{k}\geq m_{k+1}$.  For $(i,j) = (k,k+1)$ the condition becomes
$m'_{k}\geq m'_{k+1} - {\chi (\sigma ^{-1}(k+1) > \sigma ^{-1}(k))}$,
which follows from the definition of $m'_{k}$ and $m'_{k+1}$.

For $i<k = j$, the condition becomes $m_{i}\geq m'_{k} - {\chi (\sigma
^{-1}(i) > \sigma ^{-1}(k+1))}$.  By the definition of $m'_{k}$, this
is the conjunction of $m_{i}\geq m_{k+1} - {\chi (\sigma ^{-1}(i) >
\sigma ^{-1}(k+1))}$ and $m_{i}\geq m_{k} - {\chi (\sigma ^{-1}(k+1) >
\sigma ^{-1}(k))} - {\chi (\sigma ^{-1}(i) > \sigma ^{-1}(k+1))}$.
The first of these holds by hypothesis.  The second follows from the
hypothesis $m_{i}\geq m_{k} - {\chi (\sigma ^{-1}(i) > \sigma
^{-1}(k))}$ and the inequality
\begin{equation}\label{e:chi-transitivity}
\chi (\sigma ^{-1}(i) > \sigma ^{-1}(k)) \leq \chi (\sigma ^{-1}(i) >
\sigma ^{-1}(k+1)) + \chi (\sigma ^{-1}(k+1) > \sigma ^{-1}(k)),
\end{equation}
which is logically equivalent to the transitive law $(\sigma ^{-1}(i)
< \sigma ^{-1}(k+1)) \wedge (\sigma ^{-1}(k+1) < \sigma
^{-1}(k))\Rightarrow (\sigma ^{-1}(i) < \sigma ^{-1}(k))$.
\end{proof}

\begin{lemma}\label{lem:Tk-properties}
The Demazure-Lusztig operators $T_{k}$ in \eqref{e:Ti} have the
following properties, where $(r,s)\geq (a,b)$ means $r\geq a$ and
$s\geq b$.

\noindent (i) If $(r,s)\geq (c,c+1)$, then every term
$z_{k}^{u}z_{k+1}^{v}$ in $T_{k}(z_{k}^{r}z_{k+1}^{s})$ has $(u,v)\geq
(c+1,c)$.

\noindent (ii) If $(r,s)\geq (c+1,c)$, then every term
$z_{k}^{u}z_{k+1}^{v}$ in $T_{k}^{-1}(z_{k}^{r}z_{k+1}^{s})$ has
$(u,v)\geq (c,c+1)$.

\noindent (iii) If $(r,s)\geq (c,c)$, then every term
$z_{k}^{u}z_{k+1}^{v}$ in $T_{k}(z_{k}^{r}z_{k+1}^{s})$ or
$T_{k}^{-1}(z_{k}^{r}z_{k+1}^{s})$ has $(u,v)\geq (c,c)$.

\noindent (iv) If $r\geq s$, then $q^{-1}T_{k}(z_{k}^{r}z_{k+1}^{s}) =
z_{k}^{s}z_{k+1}^{r} +O(z_{k}^{s+1})$.

\noindent (v) If $r > s$, then $T_{k}^{-1}(z_{k}^{r}z_{k+1}^{s}) =
z_{k}^{s}z_{k+1}^{r} +O(z_{k}^{s+1})$.
\end{lemma}

\begin{proof}
All the properties follow from the explicit formulas
\begin{equation}\label{e:Tk-formulas}
\begin{aligned}
T_{k} (z_{k}^{r}z_{k+1}^{s}) & =
\begin{cases}
q\, z_{k}^{s}z_{k+1}^{r} + (q-1) (z_{k}^{r}z_{k+1}^{s} +
z_{k}^{r-1}z_{k+1}^{s+1} + \ldots + z_{k}^{s+1}z_{k+1}^{r-1}) & r>s,
\\
q\, z_{k}^{s}z_{k+1}^{r} &	r=s,
\\
z_{k}^{s}z_{k+1}^{r} + (1-q) (z_{k}^{s-1}z_{k+1}^{r+1} +
z_{k}^{s-2}z_{k+1}^{r+2} + \ldots + z_{k}^{r+1}z_{k+1}^{s-1}) & r< s,
\end{cases}
\\
T_{k}^{-1}(z_{k}^{r}z_{k+1}^{s}) & =
\begin{cases}
q^{-1}\, z_{k}^{s}z_{k+1}^{r} + (q^{-1}-1) (z_{k}^{r}z_{k+1}^{s} +
z_{k}^{r+1}z_{k+1}^{s-1} + \ldots + z_{k}^{s-1}z_{k+1}^{r+1}) & r <
s,
\\
q^{-1}\, z_{k}^{s}z_{k+1}^{r}&	r=s,
\\
z_{k}^{s}z_{k+1}^{r} + (1-q^{-1}) (z_{k}^{s+1}z_{k+1}^{r-1} +
z_{k}^{s+2}z_{k+1}^{r+2} + \ldots + z_{k}^{r-1}z_{k+1}^{s+1}) & r > s.
\end{cases}
\end{aligned}
\end{equation}
\end{proof}

\begin{lemma}\label{lem:Tk-on-spaces}
Let  $\sigma \in S_{l},$  $1\leq k <l$, and $m,m'\in \ZZ^{l}$ 
be as in Lemma~\ref{lem:sk-almost-decr}, and suppose $\nu \in \ZZ ^{l}$ is such
that $\nu \geq m' $ coordinate-wise.
Let $V_{m} = \ZZ [q^{\pm 1}]\{\zz ^{\mu } \mid \mu \geq m
\}$.  Then the operator $T_{k}$ satisfies

\noindent (i) $T_{k}(\zz ^{\nu }) \in V_{m}$ if $\sigma ^{-1}(k) <
\sigma ^{-1}(k+1)$;

\noindent (ii) $T_{k}^{-1}(\zz ^{\nu }) \in V_{m}$ if $\sigma ^{-1}(k)
> \sigma ^{-1}(k+1)$.
\end{lemma}

\begin{proof}
Since $T_{k}$ acts only on the variables $z_{k}$, $z_{k+1}$, we need
only consider the exponents $\mu _{k}$, $\mu _{k+1}$ of terms $\zz
^{\mu }$ occurring in $T_{k}^{\pm 1}(\zz ^{\nu })$.

In case (i), we have $(\nu _{k},\nu _{k+1})\geq (\max
(m_{k}-1,m_{k+1}),m_{k}) $ with $m_{k+1}\leq m_{k}$ and want to show
that $(\mu _{k},\mu _{k+1})\geq (m_{k},m_{k+1})$ for every term $\zz
^{\mu }$ in $T_{k}(\zz ^{\nu })$.  If $m_{k+1} < m_{k}$, this follows
from Lemma~\ref{lem:Tk-properties} (i) with $c = m_{k}-1$.  If $m_{k} =
m_{k+1}$, it follows from Lemma~\ref{lem:Tk-properties} (iii).

In case (ii), we have $(\nu _{k},\nu _{k+1})\geq (\max
(m_{k},m_{k+1}), m_{k})$ with $m_{k+1}\leq m_{k}+1$ and want to show
that $(\mu _{k},\mu _{k+1})\geq (m_{k},m_{k+1})$ for every term $\zz
^{\mu }$ in $T_{k}^{-1}(\zz ^{\nu })$.  If $m_{k+1} = m_{k} + 1$, this
follows from Lemma~\ref{lem:Tk-properties} (ii) with $c = m_{k}$.  If
$m_{k+1}\leq m_{k}$, it follows from Lemma~\ref{lem:Tk-properties}
(iii).
\end{proof}

\begin{lemma}\label{lem:E-term-bound}
Given $\sigma \in S_{l}$ and $m,\lambda \in \ZZ ^{l}$, if $m$ is
$\sigma $-almost decreasing and $\lambda \geq m$ coordinate-wise, then
for every term $\zz ^{\mu }$ with non-zero coefficient in $E^{\sigma
}_{\lambda }(\zz ;q)$, we have $\mu \geq m$ coordinate-wise.
\end{lemma}

\begin{proof}
If $\lambda $ is dominant, then $\zz ^{\lambda }$ is the only term and
the result is a tautology.  Otherwise, pick an index $k$ such that
$\lambda _{k}<\lambda _{k+1}$.  Then the recurrence \eqref{e:E-recurrence} 
gives
\begin{equation}\label{e:E-term-bound-recurrence}
E^{\sigma }_{\lambda } =
\begin{cases}
q^{-1}\, T_{k}\, E^{s_{k}\sigma}_{s_{k}(\lambda )},
 &	  \sigma^{-1}(k) < \sigma ^{-1}(k+1),
\\
T_{k}^{-1}\, E^{s_{k}\sigma }_{s_{k}(\lambda )},
 & \sigma^{-1}(k) > \sigma ^{-1}(k+1).
\end{cases}
\end{equation}
Let $m'$ be the $(s_{k}\sigma )$-almost decreasing sequence given by
Lemma~\ref{lem:sk-almost-decr} for this $\sigma $, $m$ and $k$.  Then
$m'_{k} \leq \max (m_{k},m_{k+1})\leq \max (\lambda _{k},\lambda
_{k+1}) = \lambda _{k+1}$.  Since $m'$ and $s_{k}(m)$ agree in all but
the $k$-th position, this shows that $s_{k}(\lambda )\geq m' $.  By
induction on $|{\Inv (-\lambda )}|$, we can assume that all terms $\zz
^{\nu }$ in $E^{s_{k}\sigma }_{s_{k}(\lambda )}$ satisfy $ \nu \geq
m'$.  The result now follows from Lemma~\ref{lem:Tk-on-spaces}.
\end{proof}

\begin{remark}\label{rem:term-bound-necessity}
Suppose $\lambda $ itself is $\sigma $-almost decreasing.  Then
Lemma~\ref{lem:E-term-bound} with $m = \lambda $ implies that $\mu
\geq \lambda $ coordinate-wise for every term $\zz ^{\mu }$ in
$E^{\sigma }_{\lambda }$.  Since $E^{\sigma}_{\lambda }$ is
homogeneous of degree $|\lambda |$, $\mu \geq \lambda $ implies $\mu
=\lambda $.  Hence, $E^{\sigma }_{\lambda } = \zz ^{\lambda }$.  In
fact, it can be shown that $E^{\sigma }_{\lambda } = \zz ^{\lambda }$
if and only if $\lambda $ is $\sigma $-almost decreasing.
\end{remark}

\begin{lemma}\label{lem:ns-E-specialized}
Given $\sigma $, $m$ and $\lambda $ as in
Lemma~\ref{lem:E-term-bound}, suppose that $m+\varepsilon _{j}$ is
also $\sigma $-almost decreasing, where $\varepsilon _{j}$ is the
$j$-th unit coordinate vector.  Then the coefficient of
$z_{j}^{m_{j}}$ in $E^{\sigma }_{\lambda }$ is given by
\begin{equation}\label{e:ns-E-specialized}
\langle z_{j}^{m_{j}} \rangle\, E^{\sigma }_{\lambda }(\zz ;q) = 
\begin{cases}
E^{\tau }_{\hat{\lambda}}(z_{1},\ldots,z_{j-1},z_{j+1},\ldots,z_{l};q) 
 & \lambda _{j} = m_{j},
\\
0 & \lambda _{j}>m_{j},
\end{cases}
\end{equation}
where $\hat{\lambda } = (\lambda _{1},\ldots,\lambda _{j-1},\lambda
_{j+1},\ldots, \lambda _{l})$ and $\tau \in S_{l-1}$ is the
permutation such that $\tau ^{-1}(1),\allowbreak \ldots,\allowbreak
\tau ^{-1}(l-1)$ are in the same relative order as $\sigma
^{-1}(1),\allowbreak \ldots,\allowbreak {\sigma
^{-1}(j-1)},\allowbreak {\sigma ^{-1}(j+1)},\allowbreak
\ldots,\allowbreak \sigma ^{-1}(l)$.
\end{lemma}

\begin{proof}
If $\lambda _{j}>m_{j}$, we have $\lambda \geq m+\varepsilon _{j}$.
Since $m+\varepsilon _{j}$ is assumed to be $\sigma $-almost
decreasing, the result follows from Lemma~\ref{lem:E-term-bound} in
this case.  Now assume that $\lambda _{j} = m_{j}$.

If $\lambda $ is dominant, so is $\hat{\lambda }$.  Then both sides of
\eqref{e:ns-E-specialized} reduce to $z_{1}^{\lambda _{1}}\cdots
z_{j-1}^{\lambda _{j-1}}z_{j+1}^{\lambda _{j+1}}\cdots z_{l}^{\lambda
^{l}}$.  If $\lambda $ is not dominant, we proceed by induction on
$|{\Inv (-\lambda )}|$, again using the recurrence
\eqref{e:E-term-bound-recurrence} for some index $k$ such that
$\lambda _{k}<\lambda _{k+1}$.  

If $k\not \in \{j-1,j \}$, then using Lemma~\ref{lem:sk-almost-decr}
for this $\sigma $, $m$ and $k$, and using it again with
$m+\varepsilon _{j}$ in place of $m$, we get a sequence $m'$ such that
both $m'$ and $m'+\varepsilon _{j}$ are $(s_{k}\sigma )$-almost
decreasing.  We also have $s_{k}(\lambda )\geq m'$, as in the proof of
Lemma~\ref{lem:E-term-bound}.  
Then we have
\begin{equation}\label{e:ns-E-spec-induction}
\langle z_{j}^{m_{j}} \rangle\, E^{s_{k}\sigma }_{s_{k}(\lambda )}(\zz ;q) = 
E^{s_{k'}\tau }_{s_{k'}(\hat{\lambda
})}(z_{1},\ldots,z_{j-1},z_{j+1},\ldots,z_{l};q) 
\end{equation}
by induction, where $k' = k$ if $k<j-1$ and $k' = k-1$ if $k>j$.  Note
that $s_{k'}\tau $ and $s_{k'}(\hat{\lambda })$ are to $s_{k}\sigma $
and $s_{k}(\lambda )$ as $\tau $ and $\hat{\lambda }$ are to $\sigma $
and $\lambda $.

Letting $\hat{\zz }$ denote the variables with $z_{j}$ omitted, the
recurrence for $E^{\tau }_{\hat{\lambda }}$ takes the form
\begin{equation}\label{e:hat-recurrence}
E^{\tau } _{\hat{\lambda }}(\hat{\zz };q)  =
\begin{cases}
q^{-1}\, T_{k}\, E^{s_{k'}\tau }_{s_{k'}(\hat{\lambda
})}(\hat{\zz };q), 
 &	  \sigma^{-1}(k) < \sigma ^{-1}(k+1),
\\
T_{k}^{-1}\, E^{s_{k'}\tau }_{s_{k'}(\hat{\lambda} )}(\hat{\zz };q),
 & \sigma^{-1}(k) > \sigma ^{-1}(k+1).
\end{cases}
\end{equation}
To see this, observe that the variables $z_{k}$, $z_{k+1}$ on which
$T_{k}$ acts are in positions $k'$, $k'+1$ in $\hat{\zz }$, and that
$\tau ^{-1}(k') < \tau ^{-1}(k'+1)$ if and only if $\sigma^{-1}(k) <
\sigma ^{-1}(k+1)$.  For $k\not \in \{j-1,j \}$, taking the
coefficient of $z_{j}^{m_{j}}$ commutes with $T_{k}$, so
\eqref{e:E-term-bound-recurrence}, \eqref{e:ns-E-spec-induction}, and
\eqref{e:hat-recurrence} imply $\langle z_{j}^{m_{j}} \rangle\,
E^{\sigma }_{\lambda }(\zz ;q) = E^{\tau } _{\hat{\lambda }}(\hat{\zz
};q)$, as desired.

Since $m+\varepsilon _{j}$ is $\sigma $-almost decreasing, we have
$\lambda _{j-1}\geq m_{j-1}\geq m_{j} = \lambda _{j}$, so we never
have $k=j-1$.  This leaves the case $k=j$.  In this case, let
\begin{equation}\label{e:m'-for-k=j}
\begin{gathered}
m' = (m_{1},\ldots, m_{j-1},m'_{j}, m'_{j+1}, m_{j+2},\ldots,m_{l}),\\
\end{gathered}
\end{equation}
where
\begin{equation}\label{e:m'-for-k=j-bis}
\begin{aligned}
m'_{j} & = m_{j}+\chi (\sigma ^{-1}(j) > \sigma ^{-1}(j+1)),\\
m'_{j+1} & = m_{j}.
\end{aligned}
\end{equation}
Then $s_{j}(\lambda )\geq m'$, since $\lambda _{j+1}> \lambda_{j}$
implies $ \lambda _{j+1} \geq m_{j}+1\geq m'_{j}$.  Provided that $m'$
and $m'+\varepsilon _{j+1}$ are $(s_{j}\sigma )$-almost decreasing,
the inductive hypothesis implies
\begin{equation}\label{e:by-induction}
\langle
z_{j+1}^{m'_{j+1}} \rangle \, E^{s_{j}\sigma}_{s_{j}(\lambda )}(\zz
;q) = E^{\tau }_{\hat{\lambda
}}(z_{1},\ldots,z_{j},z_{j+2},\ldots,z_{l};q),
\end{equation}
since deleting position $j+1$ from $s_{j}\sigma $ and $s_{j}(\lambda
)$ gives the same $\tau $ and $\hat{\lambda}$ as deleting position $j$
from $\sigma $ and $\lambda $.  We now verify that $m'$ and
$m'+\varepsilon _{j+1}$ are indeed $(s_{j}\sigma )$-almost decreasing.
We then complete the proof by showing that $\langle z_{j}^{m_{j}}
\rangle E^{\sigma }_{\lambda }(\zz ;q) = s_{j} \big( \langle
z_{j+1}^{m_{j}} \rangle E^{s_{j}\sigma }_{s_{j}(\lambda )}(\zz ;q)
\big)$, which is equal to $E^{\tau }_{\hat{\lambda }}(\hat{\zz };q)$
by \eqref{e:by-induction} and the fact that $m_{j} = m'_{j+1}$.

For $j\not \in \{r,s \}$, the conditions on entries in positions $r<s$
for $m'$ and $m'+\varepsilon _{j+1}$ to be $(s_{j}\sigma )$-almost
decreasing reduce to the conditions in positions $s_{j}(r)<s_{j}(s)$
for $m$ and $m+\varepsilon _{j}$ to be $\sigma $-almost decreasing.

In positions $r<j$, the required condition for both $m'$ and
$m'+\varepsilon _{j+1}$ is $m_{r}\geq m'_{j} - {\chi (\sigma ^{-1}(r)
> \sigma ^{-1}(j+1))}$.  This follows from $m_{r}\geq m_{j} +1 - {\chi
(\sigma ^{-1}(r) > \sigma ^{-1}(j))}$, which holds because
$m+\varepsilon _{j}$ is $\sigma $-almost decreasing, and $1 - {\chi
(\sigma ^{-1}(r) > \sigma ^{-1}(j))} \geq {\chi (\sigma ^{-1}(j) >
\sigma ^{-1}(j+1))} - {\chi (\sigma ^{-1}(r) > \sigma ^{-1}(j+1))}$,
which is logically equivalent to $(\sigma ^{-1}(r) < \sigma
^{-1}(j+1)) \wedge (\sigma ^{-1}(j) > \sigma ^{-1}(j+1)) \Rightarrow
(\sigma ^{-1}(r) < \sigma ^{-1}(j))$.

In positions $j$ and $s>j+1$ the required condition is $m'_{j}\geq
m_{s} - {\chi (\sigma ^{-1}(j+1) > \sigma ^{-1}(s))}$.  This follows
from $m_{j}\geq m_{s} - {\chi (\sigma ^{-1}(j) > \sigma ^{-1}(s))}$,
which holds because $m$ is $\sigma $-almost decreasing, and $ {\chi
(\sigma ^{-1}(j) > \sigma ^{-1}(s))} \leq \chi (\sigma ^{-1}(j) >
\sigma ^{-1}(j+1)) + {\chi (\sigma ^{-1}(j+1) > \sigma ^{-1}(s))}$,
which is logically equivalent to $ (\sigma ^{-1}(j) < \sigma
^{-1}(j+1)) \wedge (\sigma ^{-1}(j+1) < \sigma ^{-1}(s)) \Rightarrow
(\sigma ^{-1}(j) < \sigma ^{-1}(s))$.

Finally, in positions $j$, $j+1$, the condition on $m'+\varepsilon
_{j+1}$, which is stronger than the one on $m'$, is $m'_{j}\geq
m'_{j+1} +1 - {\chi (\sigma ^{-1}(j+1)>\sigma ^{-1}(j))}$.  This holds
with equality by the definition of $m'_{j}$, $m'_{j+1}$.

We have left to prove that $\langle z_{j}^{m_{j}} \rangle E^{\sigma
}_{\lambda }(\zz ;q) = s_{j} \big( \langle z_{j+1}^{m_{j}} \rangle
E^{s_{j}\sigma }_{s_{j}(\lambda )}(\zz ;q) \big)$.  We do this by
using the expression for $ E^{\sigma }_{\lambda }$ in terms of
$E^{s_{j}\sigma }_{s_{j}(\lambda )}$ given by the recurrence
\eqref{e:E-term-bound-recurrence} with $k = j$.

By Lemma~\ref{lem:E-term-bound}, every term $\zz ^{\nu }$ of
$E^{s_{j}\sigma }_{s_{j}(\lambda )}$ satisfies $\nu \geq m'$.

In the case $\sigma ^{-1}(j) < \sigma ^{-1}(j+1)$, this gives
$(\nu _{j},\nu _{j+1})\geq (m_{j},m_{j})$.  If $\nu _{j+1} = m_{j}$,
then we have $\langle z_{j}^{m_{j}} \rangle\, q^{-1} T_{j} \zz ^{\nu }
= s_{j}( \langle z_{j+1}^{m_{j}} \rangle\, \zz ^{\nu })$ by
Lemma~\ref{lem:Tk-properties} (iv).  If $\nu _{j+1}>m_{j}$, then
$\langle z_{j}^{m_{j}} \rangle T_{j}\zz ^{\nu } = 0 = \langle
z_{j+1}^{m_{j}} \rangle \zz ^{\nu }$ by Lemma~\ref{lem:Tk-properties}
(i), and we again have $\langle z_{j}^{m_{j}} \rangle\, q^{-1} T_{j}
\zz ^{\nu } = s_{j}( \langle z_{j+1}^{m_{j}} \rangle\, \zz ^{\nu })$.

In the case $\sigma ^{-1}(j) > \sigma ^{-1}(j+1)$, we have $(\nu
_{j},\nu _{j+1})\geq (m_{j}+1,m_{j})$.  If $\nu _{j+1} = m_{j}$, then
$\langle z_{j}^{m_{j}} \rangle\, T_{j}^{-1} \zz ^{\nu } = s_{j}(
\langle z_{j+1}^{m_{j}} \rangle\, \zz ^{\nu })$ by
Lemma~\ref{lem:Tk-properties} (v).  If $\nu _{j+1}>m_{j}$, then
$\langle z_{j}^{m_{j}} \rangle T_{j}\zz ^{\nu } = 0 = \langle
z_{j+1}^{m_{j}} \rangle \zz ^{\nu }$ by Lemma~\ref{lem:Tk-properties}
(iii), again giving $\langle z_{j}^{m_{j}} \rangle\, T_{j}^{-1} \zz
^{\nu } = s_{j}( \langle z_{j+1}^{m_{j}} \rangle\, \zz ^{\nu })$.

In either case, \eqref{e:E-term-bound-recurrence} yields $\langle
z_{j}^{m_{j}} \rangle E^{\sigma }_{\lambda } = s_{j} \big( \langle
z_{j+1}^{m_{j}} \rangle E^{s_{j}\sigma }_{s_{j}(\lambda )} \big) $, as
desired.
\end{proof}

The next two lemmas will be stated and proved for a composition
$\rr $ that we implicitly assume is strict.  However, both lemmas 
generalize immediately to weak compositions by
Remark~\ref{rem:sigma-hat-and-weak-r} (ii).

\begin{lemma}\label{lem:F-vs-E-with-rho}
Given a composition $\rr = (r_{1},\ldots,r_{k})$, a permutation
$\sigma \in S_{k}$, and a choice of $\rho _{\rr }$, let $\rr
'=w_{0}^{k}(\rr )$ and $\rho _{\rr '} = w_{0} w_{0}^{\rr } (\rho _{\rr
})$; i.e., the 
blocks of $\rho _{\rr '}$ are those of $\rho _{\rr }$
in reverse order.  Then for any tuple of partitions $ \lambda =
(\lambda _{(1)},\ldots,\lambda _{(k)})$ such that $\ell (\lambda
_{(i)})\leq r_{i}$, we have
\begin{equation}\label{e:F-vs-E-with-rho}
\zz ^{-\rho _{\rr }} F^{\sigma }_{\rr ,\lambda +\rho _{\rr }}(\zz ;q)
= w_{0}\big(\zz ^{-\rho _{\rr '}} E^{w_{0}^{k}\sigma }_{\rr ', \lambda
'+\rho _{\rr '}}(\zz ;q^{-1})\big),
\end{equation}
with notation $\lambda +\rho _{\rr }$ and $\lambda '+\rho _{\rr '}$ as
in Theorem~\ref{thm:Cauchy}, and $\lambda' = (\lambda
_{(k)},\ldots,\lambda _{(1)})$.
\end{lemma}

\begin{proof}
Follows from \eqref{e:F-vs-E-w0}.
\end{proof}

\begin{lemma}\label{lem:E-basis}
Given a composition $\rr = (r_{1},\ldots,r_{k})$ and a permutation
$\sigma \in S_{k}$, fix a choice of $\rho _{\rr }$ with $\sigma
$-almost decreasing block minima.  Then the Laurent polynomials $\zz
^{-\rho _{\rr }} E^{\sigma }_{\rr ,\lambda +\rho _{\rr }}(\zz ;q)$ for
weights $\lambda \in \NN ^{|\rr |}\cap X^{+}(\GL _{\rr })$ are
polynomials in $\zz $, and form a basis of the ring $\kk [\zz
]^{S_{\rr }}$ of $S_{\rr }$ invariant polynomials.
\end{lemma}

\begin{remark}\label{rem:E-basis}
The condition $\lambda \in \NN ^{|\rr |}\cap X^{+}(\GL _{\rr })$ means
that $\lambda $ is a concatenation of partitions $\lambda
_{(1)},\ldots,\lambda _{(k)}$ of lengths $\ell (\lambda _{(i)})\leq
r_{i}$, where $\lambda _{(i)}$ is padded with zeroes to length
$r_{i}$.
\end{remark}

\begin{proof}[Proof of Lemma~\ref{lem:E-basis}] Let $\widehat{m} =
(m_{1}^{r_{1}},\ldots,m_{k}^{r_{k}})$ be the concatenation of constant
blocks $(m_{i}^{r_{i}})$, where $m_{1},\ldots,m_{k}$ are the block
minima of $\rho _{\rr }$.  Then $\widehat{m}$ is $\sigmahat $-almost
decreasing and 
$\lambda +\rho _{\rr }\geq \widehat{m}$ for every $\lambda \in \NN
^{|\rr |}\cap X^{+}(\GL _{\rr })$, so
Lemma~\ref{lem:E-term-bound} implies that $\zz
^{-\widehat{m}}E^{\sigmahat }_{\lambda +\rho _{\rr }}(\zz ;q)$ is a
polynomial in $\zz $.  Note that $\rho' _{\rr } = \rho _{\rr } -
\widehat{m}$ is the weight satisfying \eqref{e:rho-r} whose block
minima are equal to zero.  Since $\semis _{\rr }$ commutes with
multiplication by the $S_{\rr }$ invariant monomial $\zz
^{\widehat{m}}$, we have $\zz ^{-\rho _{\rr }} E^{\sigma }_{\rr
,\lambda +\rho _{\rr }} = \zz ^{-\rho '_{\rr }} \semis _{\rr }\, \zz
^{-\widehat{m}} E^{\sigmahat }_{\lambda +\rho _{\rr }}$.

We now check that the $S_{\rr}$ invariant Laurent polynomial $\zz
^{-\rho '_{\rr }} \semis _{\rr }\, \zz ^{-\widehat{m}} E^{\sigmahat
}_{\lambda +\rho_{\rr }}$ is in fact a polynomial in $\zz$.  Consider
a monomial $\zz ^{\nu }$ appearing in $\zz^{-\widehat{m}}E^{\sigmahat
}_{\lambda +\rho _{\rr }}(\zz ;q)$, which must satisfy $\nu \in \NN
^{|\rr |}$ by the previous paragraph.  From
\eqref{e:delta_r-vs-sigma_r} we obtain
\begin{equation}\label{e:delta_r-z^nu}
\zz ^{-\rho' _{\rr }} \semis _{\rr } (\zz ^{\nu }) = \begin{cases} \pm
\chi _{\nu _{+} - \rho' _{\rr }}(\GL _{\rr }) & \text{if $\nu $ is
$\GL _{\rr }$ regular,}
\\
0 & \text{otherwise},
\end{cases}
\end{equation}
where $\nu _{+}$ is the $\GL _{\rr }$ dominant weight in the $S_{\rr
}$ orbit of $\nu $.  If $\nu \in \NN ^{|\rr |}$ is $\GL _{\rr }$ regular,
then $\nu _{+}\in \NN ^{|\rr |}\cap X^{++}(\GL _{\rr })$ satisfies $\nu
_{+} \geq \rho '_{\rr }$ coordinate-wise, so $\chi _{\nu _{+}-\rho
'_{\rr }}$ is a polynomial character.

Since the polynomial characters $\chi _{\lambda }(\GL _{\rr })$ for
$\lambda \in \NN ^{|\rr |}\cap X^{+}(\GL _{\rr })$ are a basis of
$\kk [\zz ]^{S_{\rr }}$, it follows from \eqref{e:semi-sym-triangular}
that the polynomials $\zz ^{-\rho _{\rr }} E^{\sigma }_{\rr ,\lambda
+\rho _{\rr }}$ are also a basis.
\end{proof}

\begin{cor}\label{cor:F-basis}
If, instead, the block minima of $\rho _{\rr }$ are $\sigma $-almost
increasing, Lemma~\ref{lem:E-basis} holds with $\zz ^{-\rho _{\rr }}
F^{\sigma }_{\rr ,\lambda +\rho _{\rr }}(\zz ;q)$ in place of $\zz
^{-\rho _{\rr }} E^{\sigma }_{\rr ,\lambda +\rho _{\rr }}(\zz ;q)$
\end{cor}

\begin{proof}
The reverse of a $\sigma $-almost increasing sequence is $(w_{0}\sigma
)$-almost decreasing.  In particular, if the block minima of $\rho
_{\rr }$ are $\sigma $-almost increasing, then those of $\rho _{\rr
'}$ in Lemma~\ref{lem:F-vs-E-with-rho} are $(w_{0}^{k} \sigma
)$-almost decreasing.  The corollary now follows from
Lemmas~\ref{lem:F-vs-E-with-rho} and \ref{lem:E-basis}.
\end{proof}

\begin{lemma}\label{lem:E-specialized}
Suppose we are given weak compositions $\rr = (r_{1},\ldots,r_{k})$,
$\rr ' = (r'_{1},\ldots,r'_{k})$ such that $r_{i}\leq r'_{i}$ for all
$i$, a permutation $\sigma \in S_{k}$, and a choice of $\rho _{\rr }$,
$\rho _{\rr '}$ with associated block maxima and minima, in keeping
with Convention~\ref{conv:block-max-min}.  Assume that $\rho _{\rr }$
and $\rho _{\rr '}$ have the same block maxima $M_{i}$, and that the
block minima $m_{i}$ of $\rho _{\rr '}$ and $n_{i}$ of $\rho _{\rr }$
are $\sigma $-almost decreasing and satisfy $m_{i}\leq n_{i}$ for all
$i$ (note that $r_{i}\leq r'_{i}$ already implies $m_{i}\leq n_{i}$
for $r'_{i}>0$).

Let $\zz $ (resp.\ $\zz '$) be a list of $|\rr |$ (resp.\ $|\rr '|$)
variables, subdivided into blocks $Z_{1},\ldots,Z_{k}$ (resp.\
$Z'_{1},\ldots,Z'_{k}$) of lengths $r_{1},\ldots,r_{k}$ (resp.\
$r'_{1},\ldots,r'_{k}$).  Let $\lambda _{(1)},\ldots,\lambda _{(k)}$
be partitions with $\ell (\lambda _{(i)})\leq r'_{i}$, and define
$\lambda +\rho _{\rr '}$ as in Theorem~\ref{thm:Cauchy}, so that by
Lemma~\ref{lem:E-basis}, $(\zz ')^{-\rho _{\rr '}}E^{\sigma }_{\rr '
,\lambda +\rho _{\rr '}}(\zz ';q)$ is a symmetric polynomial in each
block of variables $Z'_{i}$, and if $\ell (\lambda _{(i)})\leq r_{i}$
for all $i$, then $\zz ^{-\rho _{\rr }}E^{\sigma }_{\rr ,\lambda +\rho
_{\rr }}(\zz ;q)$ is a symmetric polynomial in each block of variables
$Z_{i}$.

Upon specializing $r_{i}$ of the variables in each block $Z'_{i}$ to $Z_{i}$
and setting the other variables to zero, we then have
\begin{equation}\label{e:specialization}
(\zz ')^{-\rho _{\rr
'}}E^{\sigma }_{\rr ' ,\lambda +\rho _{\rr '}}(\zz ';q)
\big|_{Z'_{i}\mapsto Z_{i}} = \begin{cases}
\zz ^{-\rho _{\rr }}E^{\sigma }_{\rr ,\lambda +\rho _{\rr }}(\zz
;q) &	    \text{if $\ell (\lambda _{(i)})\leq r_{i}$ for all $i$,}\\
0 & \text{otherwise}.
\end{cases}
\end{equation}
\end{lemma}

\begin{proof}
Observe that the specialization property \eqref{e:specialization} is
`transitive'---that is, given $\rr \leq \rr '\leq \rr ''$
coordinate-wise, if \eqref{e:specialization} holds for specialization
from $\zz ''$ to $\zz '$ and from $\zz '$ to $\zz $, then it holds for
specialization from $\zz ''$ to $\zz $.  This is true even if $\rr $,
$\rr '$, $\rr ''$ are weak compositions, with some blocks of variables
empty.  Using this and Lemma~\ref{lem:m<n-almost-decr}, we can reduce
to the case that the block minima of $\rho_\rr$ and $\rho_{\rr'}$
differ by a unit coordinate vector $\varepsilon _{j}$.

If $r '_{j} = 0$, then $\rr = \rr '$ and the specialization property
is trivial, so we assume that $r '_{j}>0$.  Then $r _{j} = r'_{j} -1$,
and $r _{i} = r '_{i}$ for all $i\not =j$.  Thus, we are specializing
one variable in $Z'_{j}$ to zero and leaving the other blocks $Z'_{i}$
unchanged, apart from renaming the variables.  We can also assume for
simplicity that $r _{i}>0$ for $i\not =j$, since deleting any parts $r
_{i} = r '_{i} = 0$ preserves the hypotheses and does not change the
conclusion.

Let $m = (m_{1},\ldots,m_{k})$ be the sequence of block minima of
$\rho _{\rr '}$.
That of $\rho _{\rr }$ is then $m+\varepsilon _{j}$.
Let $\widehat{m}$ be the concatenation of constant blocks
$(m_{i}^{r_{i}})$, except for the $j$-th block, which we take to be
$((m_{j}+1)^{r_{j}},m_{j})$, so its length is $r'_{j} = r_{j}+1$.  The
assumption that $m$ and $m+\varepsilon _{j}$ are $\sigma $-almost
decreasing implies that $\widehat{m}$ and $\widehat{m}+ \varepsilon
_{\jhat }$ are $\sigmahat $-almost decreasing, where $\jhat =
r_{1}+\cdots +r_{j}+1$ is the last index in the $j$-th block.

Now consider the right hand side of \eqref{e:specialization}.  By
Lemma~\ref{lem:ns-E-specialized}, since $\lambda +\rho _{\rr '}\geq
\widehat{m}$, the coefficient $\langle (z'_{\jhat })^{m_{j}} \rangle
\, E^{\sigmahat }_{\lambda +\rho _{\rr '}}(\zz ';q)$ becomes
$E^{\sigmahat }_{\lambda +\rho _{\rr }}(\zz ;q)$ after renaming the
variables if $\ell (\lambda _{(j)})\leq r_{j}$, or zero if $\ell
(\lambda _{(j)}) = r'_{j}$.  Note that $\sigmahat $ in $E^{\sigmahat
}_{\lambda +\rho _{\rr }}$ is defined with respect to $\rr $, and that
this is the $\tau $ in Lemma~\ref{lem:ns-E-specialized} when we take
$\sigma $ there to be $\sigmahat $ defined with respect to $\rr '$,
and $j$ to be $\jhat $.  Hence, the right hand side of
\eqref{e:specialization} is given by
\begin{equation}\label{e:spec-RHS}
\zz ^{-\rho _{\rr }} \semis _{\rr } \Big( \big( \langle (z'_{\jhat
})^{m_{j}} \rangle \, E^{\sigmahat }_{\lambda +\rho _{\rr '}}(\zz ';q)
\big) \big|_{\hat{\zz}'\rightarrow \zz } \Big),
\end{equation}
where $\hat{\zz}'$ stands for $\zz '$ with $z'_{\jhat }$ omitted.  

To complete the proof we need to show that \eqref{e:spec-RHS} is equal
to the left hand side of \eqref{e:specialization}, which by definition
is
\begin{equation}\label{e:spec-LHS}
(\zz ')^{-\rho _{\rr '}} \semis _{\rr '}\, E^{\sigmahat }_{\lambda
+\rho _{\rr '}}(\zz ';q) \big|_{Z'_{i}\mapsto Z_{i}}.
\end{equation}
 For this it suffices to show that
\begin{equation}\label{e:spec-LHS=RHS}
(\zz ')^{-\rho _{\rr '}} \semis _{\rr '}\, (\zz ')^{\nu }
\big|_{Z'_{i}\mapsto Z_{i}} = \zz ^{-\rho _{\rr }} \semis _{\rr }
\Big( \big( \langle (z'_{\jhat })^{m_{j}} \rangle \, (\zz ')^{\nu }
\big) \big|_{\hat{\zz}'\rightarrow \zz } \Big)
\end{equation}
for all terms $(\zz ')^{\nu}$ occurring in $E^{\sigmahat }_{\lambda
+\rho _{\rr '}} (\zz ';q)$.  Thus, by Lemma~\ref{lem:E-term-bound}, we
can assume $\nu \geq \widehat{m}$.

If $\nu $ is not $\GL _{\rr '}$ regular, then $(\zz ')^{-\rho _{\rr
'}} \semis _{\rr '}\, (\zz ')^{\nu } = 0$.  If $\nu $ is
$\GL _{\rr '}$ regular, then $(\zz ')^{-\rho _{\rr '}} \semis _{\rr
'}\, (\zz ')^{\nu } = \pm \chi _{\mu}(\GL _{\rr '})$, where $\nu _{+}
= \mu +\rho _{\rr '}$ is the dominant $\GL _{\rr '}$ weight in the
$S_{\rr '}$ orbit of $\nu $, that is, $\nu _{+}$ is the weight
obtained by sorting each block of $\nu $ into decreasing order.  As in
the proof of Lemma~\ref{lem:E-basis}, only polynomial characters $\chi
_{\mu }(\GL _{\rr '})$ arise, so we can regard the index $\mu $ as a
tuple of partitions $(\mu _{(1)},\ldots,\mu _{(k)})$ with $\ell (\mu
_{(i)})\leq r'_{i}$.

If $\nu _{\jhat } > m_{j}$, the right hand side of
\eqref{e:spec-LHS=RHS} vanishes.  Since $\nu \geq \widehat{m}$, all
entries of $\nu $ in the $j$-th block (of length $r'_{j} = r_{j}+1)$
are greater than $m_{j}$, and the same holds for $\nu _{+}$.  If $\nu
$ is not $\GL _{\rr '}$ regular, the left side of
\eqref{e:spec-LHS=RHS} vanishes immediately.  If $\nu $ is $\GL _{\rr
'}$ regular, then the corresponding character $\chi _{\mu }(\GL _{\rr
'})$ has $\ell (\mu _{(j)}) = r'_{j}$.  In that case $\chi _{\mu }(\GL
_{\rr '})$ vanishes upon specializing one of the variables $Z'_{j}$ to
zero.  Hence the left hand side of \eqref{e:spec-LHS=RHS} vanishes in
either case.

If $\nu _{\jhat } = m_{j}$, then $ \big( \langle (z'_{\jhat })^{m_{j}}
\rangle \, (\zz ')^{\nu } \big) \big|_{\hat{\zz}'\rightarrow \zz } =
\zz ^{\kappa }$, where $\kappa $ is just $\nu $ with the entry in
position $\jhat $ deleted.  In this case $\nu _{\jhat }$ is strictly
less than all other entries in the $j$-th block.  Hence, $\nu $ is $\GL
_{\rr '}$ regular if and only if $\kappa $ is $\GL _{\rr }$ regular,
and when these hold, the permutation $v\in S_{r_{j}} \subseteq S_{\rr
}$ such that $\kappa = v(\kappa _{+})$ has the same length as the
permutation $w\in S_{r'_{j}}\subseteq S_{\rr '}$ such that $\nu =w(\nu
_{+})$.  Indeed, the two permutations agree on the interval $[\jhat -
r_{j}, \jhat - 1]$, and $w$ fixes $\jhat $.

Hence, we have $(\zz ')^{-\rho _{\rr '}} \semis _{\rr '}\, (\zz
')^{\nu} = \pm \chi _{\mu}(\GL _{\rr '})$ and $\zz ^{-\rho _{\rr }}
\semis _{\rr }\, \zz ^{\kappa } = \pm \chi _{\mu}(\GL _{\rr })$
for the same tuple of partitions $\mu $ such that $\ell (\mu
_{(i)})\leq r_{i}$, and with the same sign.  Now
\eqref{e:spec-LHS=RHS} follows because $\chi _{\mu}(\GL _{\rr '})$
specializes to $\chi _{\mu}(\GL _{\rr })$ upon setting one variable in
the $j$-th block $Z'_{j}$ to zero.
\end{proof}

\begin{cor}\label{cor:F-specialized}
If the block minima of $\rho _{\rr }$ and $\rho _{\rr '}$ are $\sigma
$-almost increasing instead of $\sigma $-almost decreasing, then
Lemma~\ref{lem:E-specialized} holds with $(\zz ')^{-\rho _{\rr
'}}F^{\sigma }_{\rr ' ,\lambda +\rho _{\rr '}}(\zz ';q)$ and $\zz
^{-\rho _{\rr }}F^{\sigma }_{\rr ,\lambda +\rho _{\rr }}(\zz ;q)$ in
place of $(\zz ')^{-\rho _{\rr '}}E^{\sigma }_{\rr ' ,\lambda +\rho
_{\rr '}}(\zz ';q)$ and $\zz ^{-\rho _{\rr }}E^{\sigma }_{\rr ,\lambda
+\rho _{\rr }}(\zz ;q)$.
\end{cor}

\begin{proof}
Follows from Lemma~\ref{lem:E-specialized} using
Lemma~\ref{lem:F-vs-E-with-rho} in the same way that
Corollary~\ref{cor:F-basis} follows from Lemma~\ref{lem:E-basis}.
\end{proof}

\begin{proof}[Proof of Theorem~\ref{thm:Cauchy}] First we consider the
case when $\rr = \sS $.  Inverting the variables $y_{i}$ and $q$, we
are to prove
\begin{equation}\label{e:Cauchy-bar}
\frac{\prod _{i<j} \Omega [-q^{-1}\, t\, X_{i}
\overline{Y_{j}}]}{\prod _{i\leq j} \Omega [- t\, X_{i}
\overline{Y_{j}}]} = \sum _{\lambda } t^{|\lambda |}\, \xx ^{-\rho
_{\rr }} E^{\sigma }_{\rr ,\lambda +\rho _{\rr }}(\xx ;q) \,
\overline{ \yy ^{-\rho _{\rr}}F^{\sigma }_{\rr ,\lambda +\rho _{\rr
}}(\yy ;q)},
\end{equation}
where the index $\lambda $ ranges over all $k$-tuples of partitions
with $\ell (\lambda _{(i)})\leq r_{i}$.  By
Corollary~\ref{cor:F-basis}, the functions $\overline{ \yy ^{-\rho
_{\rr}}F^{\sigma }_{\rr ,\lambda +\rho _{\rr }}(\yy ;q)}$ appearing in
the sum form a basis of $\kk [y_{1}^{-1},\ldots,y_{r}^{-1}]^{S_{\rr
}}$.  We need to show that the coefficient of each such basis element
in the product on the left hand side is $t^{|\lambda |}\, \xx ^{-\rho
_{\rr }} E^{\sigma }_{\rr ,\lambda +\rho _{\rr }}(\xx ;q)$.  By
Proposition~\ref{prop:orthogonality}, the coefficient in question is
given by
\begin{equation}\label{e:F-coefficient}
\big\langle \, \yy ^{-\rho _{\rr }} E^{\sigma }_{\rr ,\lambda +\rho
_{\rr }}(\yy ;q),\, \frac{\prod _{i<j} \Omega [-q^{-1}\, t\, X_{i}
\overline{Y_{j}}]}{\prod _{i\leq j} \Omega [- t\, X_{i}
\overline{Y_{j}}]}\, \big\rangle ^{\rr }_{q},
\end{equation}
where the inner product is in the $\yy $ variables.  Note that it is
permissible to perform such operations term by term in the power
series in $t$ on each side of \eqref{e:Cauchy-bar}.

By Lemma~\ref{lem:E-basis}, $\yy ^{-\rho _{\rr }} E^{\sigma }_{\rr
,\lambda +\rho _{\rr }}(\yy ;q)$ is in fact a polynomial in $\yy $ and
not just a Laurent polynomial.  Moreover, it is homogeneous of degree
$|\lambda |$.  Hence, the result follows if we show that
\begin{equation}\label{e:reproducing}
 f(t\, \xx ) = \big\langle \, f(\yy ),\, \frac{\prod _{i<j} \Omega
[-q^{-1}\, t\, X_{i} \overline{Y_{j}}]}{\prod _{i\leq j} \Omega [- t\,
X_{i} \overline{Y_{j}}]}\, \big\rangle ^{\rr }_{q}
\end{equation}
for every $S_{\rr }$ invariant polynomial $f\in \kk [\yy ]^{S_{\rr
}}$.  Let $X_{i}$ and $Y_{i}$ outside the plethystic brackets stand
for the list of variables in each block, so that $f(\yy ) =
f(Y_{1},\ldots,Y_{k})$, for instance.  Using the definition
\eqref{e:r-inner-prod} of the inner product, the right hand side of
\eqref{e:reproducing} can then be written
\begin{equation}\label{e:reproducing-RHS}
\langle 1_{\GL _{\rr }} \rangle\, f(Y_{1},\ldots,Y_{k})\, \frac{\prod
_{i<j} \Omega [-q^{-1}\, t\, X_{i} \overline{Y_{j}}]}{\prod _{i\leq j}
\Omega [- t\, X_{i} \overline{Y_{j}}]}\, \prod _{i<j} \frac{\Omega [-
Y_{i} \overline{Y_{j}}]}{\Omega [- q^{-1} Y_{i} \overline{Y_{j}}]},
\end{equation}
where the coefficient $\langle 1_{\GL _{\rr }} \rangle$ is taken in
the $\yy $ variables.  Note that this is the same as taking the
coefficient of $1_{\GL _{r_{i}}}(Y_{i})$ in each block of variables
separately.

The only part of \eqref{e:reproducing-RHS} that involves
$\overline{Y_{1}}$ is the factor $\Omega [-t
X_{1}\overline{Y_{1}}]^{-1} = \Omega [t X_{1}\overline{Y_{1}}]$;
everything else involves only symmetric polynomials in $Y_{1}$.  The
classical Cauchy identity implies
\begin{equation}\label{e:classical-Cauchy}
\langle 1_{\GL _{r_{1}}} \rangle g(Y_{1}) \, \Omega [t X_{1}
\overline{Y_{1}}] = g(t X_{1})
\end{equation}
for every symmetric polynomial $g(Y_{1})$, by reducing to the case
that $g$ is a Schur function.  Taking the coefficient $\langle 1_{\GL
_{\rr }} \rangle$ by starting with $\GL _{r_{1}}$ and using
\eqref{e:classical-Cauchy} reduces \eqref{e:reproducing-RHS} to
\begin{equation}\label{e:reproducing-RHS-bis}
\langle 1_{\GL _{(r_{2},\ldots,r_{k})}} \rangle\, f(t\, X_{1},
Y_{2},\ldots,Y_{k})\, \frac{\prod _{1<i<j} \Omega [-q^{-1}\, t\, X_{i}
\overline{Y_{j}}]}{\prod _{1<i\leq j} \Omega [- t\, X_{i}
\overline{Y_{j}}]}\, \prod _{1<i<j} \frac{\Omega [- Y_{i}
\overline{Y_{j}}]}{\Omega [- q^{-1} Y_{i} \overline{Y_{j}}]},
\end{equation}
once we observe that after removing the factor $\Omega [t\,
X_{1}\overline{Y_{1}}]$ and setting $Y_{1} = t\, X_{1}$ in the rest,
all factors with index $i=1$ cancel.  We can assume by induction on
$k$ that \eqref{e:reproducing-RHS-bis} reduces to $f(t\, \xx )$.

For the general case, choose an integer $N$ less than or equal to all
$M_{i}$, $m_{i}$ and $n_{i}$.  Define $\rr ' =(r'_{1},\ldots, r'_{k})$
by $r'_{i} = M_{i}-N+1$ and choose $\rho _{\rr '}$ to have block
maxima $M_{i}$ and (hence) constant block minima $m'_{i} = N$.  The
constant sequence $(N,\ldots,N)$ is both $\sigma $-almost decreasing
and $\sigma $-almost increasing, so the case of the theorem with equal
compositions holds for $\rr '$ and $\rho _{\rr '}$, by what was shown
above (note that since we chose $N\leq M_{i}$ for all $i$, $\rho _{\rr
'}$ has no artificial zero-length blocks).  Denote the blocks of
variables in this case by $X'_{i}$, $Y'_{i}$.  By the choice of $N$,
we have $r_{i},s_{i}\leq r'_{i}$, so there are at least as many
variables in each block $X'_{i}$, $Y'_{i}$ as in $X_{i}$, $Y_{i}$.

Specializing $r_{i}$ of the variables in each $X'_{i}$ to $X_{i}$ and
$s_{i}$ of the variables in $Y'_{i}$ to $Y_{i}$, and setting the other
variables to zero, the left hand side of \eqref{e:Cauchy} with both
compositions equal to $\rr '$ reduces to the left hand side for
compositions $\rr $ and $\sS $.  By Lemma~\ref{lem:E-specialized} and
Corollary~\ref{cor:F-specialized}, the same thing happens on the right
hand side.  Thus, the general case follows from the case already
proven.
\end{proof}

\subsection{Winding permutations}
\label{ss:winding}

We use the following notions from \cite[Definition~5.2.1]{BHMPS-paths}.

\begin{defn}\label{def:winding}
A permutation $\sigma \in S_{k}$ is a {\em winding permutation} if
$\sigma (1),\ldots,\sigma (k)$ are in the same relative order as
$c_{1},\ldots,c_{k}$, where $c_{i} = \fp{y + x\, i}$ are the
fractional parts of an arithmetic progression, for any real $x$, $y$
with $x$ assumed irrational, so the $c_{i}$ are distinct.

The {\em descent indicator} of $\sigma $ is the $\{0,1 \}$-valued
vector $(\eta _{1},\ldots,\eta _{k-1})$ defined by
\begin{equation}\label{e:descent-indicator}
\eta _{i} = \chi (\sigma (i)>\sigma (i+1)).
\end{equation}

The {\em head} and {\em tail} of $\sigma $ are the permutations $\tau
, \theta \in S_{k-1}$ such that $\tau (1),\ldots,\tau (k-1)$ are in
the same relative order as $\sigma (1),\ldots,\sigma (k-1)$ and
$\theta (1),\ldots,\theta (k-1)$ are in the same relative order as
$\sigma (2),\ldots,\sigma (k)$.
\end{defn}

Proposition~\ref{prop:winding}, below, is the counterpart of
\cite[Proposition~5.2.2]{BHMPS-paths} for semi-symmetric
Hall-Littlewood polynomials.  We start with a more general identity.

\begin{lemma}\label{lem:general-winding}
Let $\tau ,\theta \in S_{k}$ and $\eta \in \ZZ ^{k}$ be such that
$|\eta _{i} - \eta _{j}|\leq 1$ for all $i,j$, and
\begin{itemize}
\item [(a)] $\eta _{i} = \eta _{j}$ implies $\tau (i)<\tau
(j)\Leftrightarrow \theta (i)<\theta (j)$, i.e., $\tau $ and $\theta $
are in the same relative order in positions $i,j$;
\item [(b)] $\eta _{i} - \eta _{j} = 1$ implies $\tau (i)>\tau (j)$
and $\theta (i)<\theta (j)$, i.e., $\theta (\eta )$ is dominant and
$\tau (\eta )$ is antidominant.
\end{itemize}
Given a composition $\rr = (r_{1},\ldots,r_{k})$, let $\widehat{\eta }
= (\eta _{1}^{r_{1}},\ldots,\eta _{k}^{r_{k}})$ be the concatenation
of constant blocks $(\eta _{i}^{r_{i}})$.  Then for every $\mu \in
X^{++}(\GL _{\rr })$ we have the identities
\begin{align}\label{e:winding-1}
E^{\theta ^{-1}}_{\rr ,\mu }(\zz ;q) & = \zz ^{\widehat{\eta}}\,
E^{\tau ^{-1}}_{\rr ,\mu - \widehat{\eta}}(\zz ;q),
\\
\label{e:winding-2} F^{\theta ^{-1}}_{\rr ,\mu }(\zz ;q) & = \zz
^{\widehat{\eta}}\, F^{\tau ^{-1}}_{\rr , \mu - \widehat{\eta}}(\zz
;q).
\end{align}
\end{lemma}

\begin{proof}
First we show that \eqref{e:winding-1} implies \eqref{e:winding-2}.
Note that $w_{0}\tau $, $w_{0}\theta $ and $-\eta $ satisfy the same
hypothesis as $\tau $, $\theta $ and $\eta $.  Hence, assuming the
validity of \eqref{e:winding-1}, we have
\begin{multline}\label{e:winding-1->2}
F^{\theta ^{-1}}_{\rr ,\mu }(\zz ;q) =
\zz ^{\rho _{\rr } - w_{0}^{\rr
}(\rho _{\rr })} \overline{E^{( w_{0} \theta )^{-1}}_{\rr ,-w_{0}^{\rr
}(\mu )}(\zz ;q)}\\
 = \zz ^{\rho _{\rr } - w_{0}^{\rr }(\rho _{\rr })} \overline{\zz
^{-\widehat{\eta }} E^{( w_{0} \tau )^{-1}}_{\rr ,-w_{0}^{\rr }(\mu )
+ \widehat{\eta }}(\zz ;q)} = \zz ^{\widehat{\eta}}\, F^{\tau
^{-1}}_{\rr , \mu - \widehat{\eta}}(\zz ;q),
\end{multline}
using \eqref{e:F-vs-E-bar} and \eqref{e:shift}.  It remains to prove
\eqref{e:winding-1}.

Let $\sim $ denote equality up to a non-zero scalar factor.  Since
$\widehat{\eta }$ is $S_{\rr }$ invariant,
\eqref{e:semi-sym-triangular} implies that both sides of
\eqref{e:winding-1} have leading term $\zz ^{\rho _{\rr }}\chi _{\mu
-\rho _{\rr }}$ with coefficient $1$.  Hence, if \eqref{e:winding-1}
holds up to $\sim $ equivalence, then it holds with equality.

Let $l = |\rr |$, and let $\lambda _{+}$ denote the dominant weight in
the $S_{l}$ orbit of any $\GL _{l}$ weight $\lambda \in \ZZ ^{l}$.
The entries of $\eta $ and $\widehat{\eta }$ take at most two values
$c$ and $c+1$.  As in the proof of
\cite[Proposition~5.2.2]{BHMPS-paths}, this implies that there is a
$w\in S_{l}$ such that $\mu =w(\mu _{+})$ and $\mu -\widehat{\eta } =
w((\mu -\widehat{\eta })_{+})$.  Note that $(\mu -\widehat{\eta })_{+}
= \mu _{+} - w^{-1}(\widehat{\eta })$.  Hence, up to $\sim $
equivalence, \eqref{e:winding-1} can be written
\begin{equation}\label{e:winding-1-op}
\semis _{\rr }\, T_{\widetilde{\theta }}^{-1} \, T_{\widetilde{\theta
} w}(\zz ^{\mu _{+}}) \sim \zz ^{\widehat{\eta}}\, \semis _{\rr }\,
T_{\widetilde{\tau }}^{-1} \, T_{\widetilde{\tau } w}(\zz
^{-w^{-1}(\widehat{\eta })}\, \zz ^{\mu _{+}}),
\end{equation}
where $\widetilde{\theta } = (\widehat{\theta ^{-1}})^{-1}$,
$\widetilde{\tau } = (\widehat{\tau ^{-1}})^{-1}$.  Since
multiplication by the $S_{\rr }$ invariant monomial $\zz
^{\widehat{\eta }}$ commutes with $\semis _{\rr }$,
\eqref{e:winding-1-op} follows if we prove the operator identity
$T_{\widetilde{\theta }}^{-1} \, T_{\widetilde{\theta } w} \sim \zz
^{\widehat{\eta}}\, T_{\widetilde{\tau }}^{-1} \, T_{\widetilde{\tau }
w}\, \zz ^{-w^{-1}(\widehat{\eta })}$, or equivalently
\begin{equation}\label{e:operator-identity}
T_{\widetilde{\theta } w}^{-1}\, T_{\widetilde{\theta }}\, \, \zz
^{\widehat{\eta}}\, T_{\widetilde{\tau }}^{-1} \, T_{\widetilde{\tau }
w} \sim \zz ^{w^{-1}(\widehat{\eta })}.
\end{equation}

We now prove \eqref{e:operator-identity} for all $w\in S_{l}$ by
induction on $\ell (w)$.  For this we use the well-known operator
identities (the same as \cite[(115)]{BHMPS-paths})
\begin{align}\label{e:affine-Hecke-0}
T_{i}\, \zz ^{\mu }\, T_{i}^{-1} = T_{i}^{-1}\, \zz ^{\mu }\, T_{i} =
\zz ^{\mu } = \zz ^{s_{i}\mu }
&\quad
\text{ if } \langle \alpha _{i}^{\vee }, \mu \rangle \defeq  \mu _{i} - \mu
_{i+1} = 0,
\\
\label{e:affine-Hecke-1}
T_{i}\, \zz ^{\mu }\, T_{i} =q\, \zz ^{s_{i}\mu }
&\quad 
\text{ if } \langle \alpha _{i}^{\vee }, \mu \rangle = -1,
\\
\label{e:affine-Hecke+1}
T_{i}^{-1}\, \zz ^{\mu }\, T_{i}^{-1} = q^{-1}\, \zz ^{s_{i}\mu }
&\quad 
\text{ if } \langle \alpha _{i}^{\vee }, \mu \rangle = 1,
\end{align}
which follow directly from the definition of $T_{i}$.

The base case $w=1$ of the induction is trivial.  Otherwise, let $w = v
s_{i} > v$, and assume by induction that
\begin{equation}\label{e:operator-identity-induct}
T_{\widetilde{\theta } v}^{-1}\, T_{\widetilde{\theta }} \,\, \zz
^{\widehat{\eta}}\, T_{\widetilde{\tau }}^{-1} \, T_{\widetilde{\tau }
v} \sim \zz ^{v^{-1}(\widehat{\eta })}.
\end{equation}
We have $T_{\widetilde{\theta } w}^{-1} = T_{i}^{e_{1}}
T_{\widetilde{\theta } v}^{-1}$ and $T_{\widetilde{\tau } w} =
T_{\widetilde{\tau } v} T_{i}^{e_{2}}$, where
\begin{equation}\label{e:Ti-signs}
e_{1} = \begin{cases}
-1 & \widetilde{\theta } v s_{i} > \widetilde{\theta } v\\
1  & \widetilde{\theta } v s_{i} < \widetilde{\theta } v
\end{cases}\, ,
\qquad 
e_{2} = \begin{cases}
1 & \widetilde{\tau } v s_{i} > \widetilde{\tau } v\\
-1  & \widetilde{\tau } v s_{i} < \widetilde{\tau } v
\end{cases}\, .
\end{equation}
Then \eqref{e:operator-identity-induct} implies
\eqref{e:operator-identity}, provided we show that
\begin{equation}\label{e:op-ident-step}
T_{i}^{e_{1}} \zz ^{v^{-1}(\widehat{\eta })} T_{i}^{e_{2}} \sim \zz
^{s_{i} v^{-1}(\widehat{\eta })} = \zz ^{w^{-1}(\widehat{\eta })}.
\end{equation}
Set $a = v(i)$, $b = v(i+1)$, and let $a',b'\in [k]$ be the indices of
the blocks containing $a$ and $b$ in the partition of $[l]$ into
intervals of lengths $r_{1},\ldots,r_{k}$.  Note that $vs_{i}>v$
implies $a<b$, and that, in the same way, we have $\widehat{\tau
}vs_{i}> \widehat{\tau }v \Leftrightarrow \widehat{\tau
}(a)<\widehat{\tau }(b)$ and $\widehat{\theta }vs_{i}> \widehat{\theta
}v \Leftrightarrow \widehat{\theta }(a)<\widehat{\theta }(b)$.

{\em Case I:} $\widehat{\eta }_{a} = \widehat{\eta }_{b}$.  One way
this can happen is if $a' = b'$, so $a$ and $b$ are in the same block.
Since $\widetilde{\tau }$ and $\widetilde{\theta }$ are increasing on
each block, we then have $\widetilde{\theta } v s_{i} >
\widetilde{\theta } v$, $\widetilde{\tau } v s_{i} > \widetilde{\tau }
v$, $e_{1} = -1$, $e_{2} = 1$.  Otherwise, if $a'\not =b'$, we have
$\eta _{a'} = \eta _{b'}$.  By hypothesis, we then have $\tau
(a')<\tau (b') \Leftrightarrow \theta (a')<\theta (b')$.  By
construction, this implies $\widetilde{\tau }(a)<\widetilde{\tau }(b)
\Leftrightarrow \widetilde{\theta}(a)< \widetilde{\theta}(b)$, or
equivalently $\widetilde{\tau }v s_{i} > \widetilde{\tau }v
\Leftrightarrow \widetilde{\theta }v s_{i} > \widetilde{\theta }v $,
hence $e_{2} = -e_{1}$.  Thus, we have $e_{2} = -e_{1}$ either way, and
since $\langle \alpha _{i}^{\vee }, v^{-1}(\widehat{\eta }) \rangle =
\widehat{\eta }_{a}-\widehat{\eta }_{b} = 0$, \eqref{e:op-ident-step}
reduces to \eqref{e:affine-Hecke-0}.

{\em Case II:} $\widehat{\eta }_{a} - \widehat{\eta }_{b} = \eta _{a'}
- \eta _{b'} = 1$.  Then $\tau (a')>\tau (b')$ and $\theta (a')<\theta
(b')$ by hypothesis, which implies $\widetilde{\tau
}(a)>\widetilde{\tau }(b)$ and $\widetilde{\theta
}(a)<\widetilde{\theta }(b)$ by construction.  In other words,
$\widetilde{\tau }v s_{i}>\widetilde{\tau }v$ and $\widetilde{\theta
}v s_{i}<\widetilde{\theta }v$, so $e_{1} = e_{2} = -1$.  In this
case, $\langle \alpha _{i}^{\vee }, v^{-1}(\widehat{\eta }) \rangle =
1$, so \eqref{e:op-ident-step} reduces to \eqref{e:affine-Hecke+1}.

{\em Case III:} $\widehat{\eta }_{a} - \widehat{\eta }_{b} = \eta
_{a'} - \eta _{b'} = -1$.  The reasoning in Case II with $a$ and $b$
exchanged gives $e_{1} = e_{2} =1$.  Since $\langle \alpha _{i}^{\vee
}, v^{-1}(\widehat{\eta }) \rangle = -1$, \eqref{e:op-ident-step}
reduces to \eqref{e:affine-Hecke-1}.
\end{proof}

\begin{remark}\label{rem:winding-proof}
In the proof of Lemma~\ref{lem:general-winding}, we implicitly assumed
that $\rr $ was a strict composition.  However, with the conventions
in Remark~\ref{rem:sigma-hat-and-weak-r} (ii), the weak composition
case follows from the strict composition case.
\end{remark}

\begin{prop}\label{prop:winding}
If $\sigma \in S_{k+1}$ is a winding permutation, then its descent
indicator $\eta = (\eta _{1},\ldots,\eta _{k})$ and head and tail
$\tau ,\theta \in S_{k}$ satisfy the hypotheses of
Lemma~\ref{lem:general-winding}.  Hence, given a (strict or weak)
composition $\rr =(r_{1},\ldots,r_{k})$, identities
(\ref{e:winding-1}--\ref{e:winding-2}) hold for every $\mu \in
X^{++}(\GL _{\rr })$.
\end{prop}

\begin{proof}
By definition, $\eta $ is $\{0,1 \}$-valued, so $|\eta _{i} - \eta
_{j}|<1$ for all $i,j$.  Let ${\mathbf c} = (c_{1},\ldots,c_{k+1})$ be
the sequence $c_{i} = \fp{y + xi}$ in the definition of the winding
permutation $\sigma $.  Since adding an integer to $x$ does not change
${\mathbf c}$, we can assume $0<x<1$.  Then, since $\sigma $ and
${\mathbf c}$ are in the same relative order, we have
\begin{equation}\label{e:descent-rules}
\begin{aligned}
\eta _{i} & = 1 \ \Leftrightarrow\ \sigma (i)>\sigma (i+1) \
\Leftrightarrow \ c_{i}>c_{i+1} \ \Leftrightarrow \ c_{i+1}=
c_{i}+x-1,
\\
\eta _{i} & = 0 \ \Leftrightarrow \ \sigma (i)<\sigma (i+1) \
\Leftrightarrow \ c_{i}<c_{i+1} \ \Leftrightarrow \ c_{i+1}= c_{i}+x.
\end{aligned}
\end{equation}

If $\eta _{i} = \eta _{j}$, then $c_{i+1}-c_{i} = c_{j+1}-c_{j}$, so
$c_{i+1}-c_{j+1} = c_{i}-c_{j}$ and $c_{i}<c_{j} \Leftrightarrow
c_{i+1}<c_{j+1}$.  Then $\sigma (i)<\sigma (j) \Leftrightarrow \sigma
(i+1)<\sigma (j+1)$, or equivalently, $\tau (i)<\tau (j)
\Leftrightarrow \theta (i)<\theta (j)$.  This shows that hypothesis (a) in
Lemma~\ref{lem:general-winding} is satisfied.

If $\eta _{i} - \eta _{j} = 1$, that is, if $\eta _{i} = 1$ and $\eta
_{j} = 0$, then $c_{i+1} = c_{i}+x-1$ and $c_{j+1}= c_{j} +x$ imply
$c_{i+1} - c_{j+1} = c_{i} - c_{j} -1$.  Since $|c_{i} - c_{j}|<1$ and
$|c_{i+1} - c_{j+1}|<1$, we must have $c_{i} - c_{j} > 0$ and $c_{i+1}
- c_{j+1} < 0$.  Then $\sigma (i)>\sigma (j)$ and $\sigma (i+1)<\sigma
(j+1)$, or equivalently, $\tau (i)>\tau (j)$ and $\theta (i)<\theta
(j)$.  This shows that hypothesis (b) in
Lemma~\ref{lem:general-winding} is satisfied.

By Remark~\ref{rem:winding-proof}, we can conclude that
(\ref{e:winding-1}--\ref{e:winding-2}) hold even if $\rr $ is a weak
composition.
\end{proof}

\section{Proof of the main results}
\label{s:main-theorems}

\subsection{Stable form of the main theorem}
\label{ss:stable-main}

In \S \ref{ss:nest-identity-proof}, we prove the combinatorial
Theorem~\ref{thm:nest-identity} by restricting to the polynomial part
of a stronger infinite series identity, given by the following
theorem, which expresses the full Catalanimal in \eqref{e:H-for-den}
in terms of LLT series.

\begin{thm}\label{thm:stable-main}
Given a positive integer $h$ and real numbers $s$, $p$ with $p$
irrational, let 
\begin{equation}\label{e:b-and-c}
b_{i} = \lfloor s - p\, (i-1) \rfloor - \lfloor s - p\, i \rfloor, \qquad 
c_{i} = \fp{s - p\, (i-1)}
\end{equation}
for $i=1,\ldots,h+1$, where $\fp{a} = a - \lfloor a \rfloor$ denotes
the fractional part of $a$.  Let $\sigma \in S_{h+1}$ be the
permutation such that $\sigma (1),\ldots,\sigma (h+1)$ are in the same
relative order as $c_{1},\ldots,c_{h+1}$, and let $\tau ,\theta \in
S_{h}$ be its head and tail.  Let $(u_{1},\ldots,u_{h})$,
$(v_{1},\ldots,v_{h})$ be integer sequences which are respectively
$\theta ^{-1}$-almost decreasing and $\tau ^{-1}$-almost increasing,
that is,
\begin{equation}\label{e:u-and-v}
u_{i}\geq u_{j} - \chi (\theta (i)>\theta (j)),\quad v_{i}\leq v_{j} +
\chi (\tau (i)<\tau (j))\qquad \text{for all $i<j$},
\end{equation}
and let $\gammabold = (\gamma _{1},\ldots,\gamma _{h})\in \ZZ
_{+}^{h}$ be a sequence of positive integers with first differences
\begin{equation}\label{e:gamma-diff}
\gamma _{i+1} - \gamma _{i} = u_{i} - v_{i+1}.
\end{equation}
Fix $\rho _{w_{0}(\gammabold )}, \rho '_{w_{0}(\gammabold )}$
satisfying \eqref{e:rho-r} for $\GL _{w_{0}(\gammabold ) }$, with
block minima $u_{h+1-i}$ for $\rho _{w_{0}(\gammabold )}$ and
$v_{h+1-i}$ for $\rho '_{w_{0}(\gammabold )}$.  Then
\begin{multline}\label{e:stable-main}
H(R_{q},R_{t},R_{qt},((u_{1}-v_{1}+b_{1})^{\gamma _{1}}, \ldots,
(u_{h}-v_{h}+b_{h})^{\gamma _{h}})) \\
 = \sum _{\lambda }t^{|\lambda |} \Lcal ^{\tau
w_{0}}_{w_{0}(\gammabold ),\, ((0;\lambda ) +(b_{h}^{\gamma
_{h}},\ldots,b_{1}^{\gamma _{1}})+\rho _{w_{0}(\gammabold )
})/((\lambda;0) +\rho' _{w_{0}(\gammabold ) })}(\zz ;q),
\end{multline}
where $\lambda $ ranges over tuples of partitions $(\lambda
_{(h-1)},\ldots,\lambda _{(1)})$ such that $\ell (\lambda _{(i)})\leq
\min (\gamma _{i}, \gamma _{i+1})$, we form $(0;\lambda )$ and
$(\lambda ;0)$ by prepending or appending an empty partition to
$\lambda $, and we interpret these as weights in $X^{+}(\GL
_{w_{0}(\gammabold ) })$ by padding the $i$-th component to length
$\gamma _{h+1-i}$ and concatenating, as in Theorem~\ref{thm:Cauchy}.

The root sets are defined by $\alpha _{ij}\in R_{q} = R_{t}$ if $i<j$
are in distinct blocks of the partition of $\{1,\ldots,|\gammabold |
\}$ into intervals of length $\gamma _{i}$, and $\alpha _{ij}\in
R_{qt}$ if $i<j$ are in distinct, non-adjacent blocks; in other words,
$R_{q} = R_{t} = R_{+}\setminus R_{+}(\GL _{\gammabold })$ and $R_{qt}
= [R_{q},R_{t}]$.
\end{thm}

\begin{example}\label{ex:stable-main}
Before turning to the proof, we give an example to clarify the
notation in Theorem~\ref{thm:stable-main}.  Let $h = 4$, $p \approx
.67$, and $s = 6.5$.  Then $(b_{1},\ldots,b_{5}) = (1,0,1,1,0)$ and
$(c_{1},\ldots,c_{5})\approx (.5,.83,.16,.49,.82)$.  The permutation
$\sigma $ with the same relative order as the numbers $c_{i}$ is
$\sigma =(3,5,1,2,4)$ in one-line notation.  Its head and tail are
$\tau =(3,4,1,2)$ and $\theta =(4,1,2,3)$.

The sequences $(u_{1},\ldots,u_{4}) = (1,2,0,0)$ and
$(v_{1},\ldots,v_{4}) = (1,0,3,2)$ satisfy the required
almost-decreasing/increasing conditions.  The sequence $\gammabold
=(3,4,3,1)$ has first differences
$(u_{1},u_{2},u_{3})-(v_{2},v_{3},v_{4}) = (1,-1,-2)$.

On the left hand side of \eqref{e:stable-main} we have the Catalanimal
\[
H(R_{q},R_{t},R_{qt},(1,1,1,2,2,2,2,-2,-2,-2,-1))
\]
in $|\gammabold | = 11$ variables $\zz = z_{1},\ldots,z_{11}$.  A root
$\alpha _{ij}$ ($i<j$) belongs to $R_{q} = R_{t}$ if $i$ and $j$ are
in distinct blocks of the partition $\{\{1,2,3 \}, \{4,5,6,7
\},\{8,9,10 \},\{11 \} \}$, and to $R_{qt}$ if $i$ and $j$ are in
non-adjacent blocks.

The terms on the right hand side of \eqref{e:stable-main}
are indexed
by triples of partitions $\lambda _{(3)}=(\lambda _{3,1})$, $\lambda
_{(2)} = (\lambda _{2,1},\lambda _{2,2},\lambda _{2,3})$, $\lambda
_{(1)} = (\lambda _{1,1},\lambda _{1,2},\lambda _{1,3})$ of lengths at
most $1$, $3$ and $3$, respectively.  More explicitly, the term indexed by a
given triple is
\[
\Lcal ^{(2,1,4,3)}_{(1,3,4,3),\beta /\alpha }(\zz ;q),\quad \text{where}\quad  
\arraycolsep=.3ex
\begin{array}{rcccccccccccl}
\beta \ = \ (& 1,& 1,& 1,& 1,& 0,& 0,& 0,& 0,& 1,& 1,& 1& )\\
        +\, (& 0,& 2,& 1,& 0,& 5,& 4,& 3,& 2,& 3,& 2,& 1& )\\
        +\, (& 0,& \lambda _{3,1},& 0,& 0,& \lambda _{2,1},& \lambda _{2,2},&
      \lambda _{2,3},& 0,& \lambda _{1,1},& \lambda _{1,2},& \lambda _{1,3}& )
\\[1.5ex]
\alpha\  = \ (& 2,& 5,& 4,& 3,& 3,& 2,& 1,& 0,& 3,& 2,& 1& )\\
   +\, (& \lambda _{3,1},& \lambda _{2,1},& \lambda _{2,2},& \lambda _{2,3},&
   \lambda _{1,1},& \lambda _{1,2},& \lambda _{1,3},& 0,& 0,& 0,& 0& ).
\end{array}
\]
Note that for $\lambda = \varnothing$, $\beta -\alpha $ reduces to the
reverse of the weight vector in the Catalanimal.
\end{example}

\begin{proof}[Proof of Theorem~\ref{thm:stable-main}]
We will prove a stronger identity
\begin{multline}\label{e:unstraightened}
\zz ^{((u_{1}-v_{1}+b_{1})^{\gamma _{1}}, \ldots,
(u_{h}-v_{h}+b_{h})^{\gamma _{h}})} \frac{\prod _{\alpha \in R_{qt}}
(1 - q\, t\, \zz ^{\alpha })}{\prod _{\alpha \in R_{t}} (1- t\, \zz
^{\alpha })} \\
= \sum _{\lambda }t^{|\lambda |} w_{0} \Big( F^{(\tau
w_{0})^{-1}}_{w_{0}(\gammabold ),\, ((0;\lambda ) +(b_{h}^{\gamma
_{h}},\ldots,b_{1}^{\gamma _{1}})+\rho _{w_{0}(\gammabold )})}(\zz
;q)\, \overline{E^{(\tau w_{0})^{-1}}_{w_{0}(\gammabold ),((\lambda;0)
+\rho' _{w_{0}(\gammabold )})}(\zz ;q)} \Big).
\end{multline}
Then \eqref{e:stable-main} follows from
Proposition~\ref{prop:L-via-Hq} after applying the operator $\Hbold
_{q}^{\gammabold }$ in \eqref{e:Hq-r} to both sides of
\eqref{e:unstraightened}.

By construction, we have $\gamma _{i}+u_{i} = \gamma _{i+1}+v_{i+1}$.
Define
\begin{equation}\label{e:Mi}
M_{1} = \gamma _{1}+v_{1}-1, \, M_{2} = \gamma _{2}+v_{2}-1 = \gamma
_{1}+u_{1}-1, \, \ldots, \, M_{h+1} = \gamma _{h}+u_{h}-1.
\end{equation}
We now apply our Cauchy identity, Theorem~\ref{thm:Cauchy}, with
$k=h+1$, taking the $\sigma $ there to be $(\sigma w_{0})^{-1}$, 
setting
\begin{equation}\label{e:Cauchy-inputs}
\rr = (0,\gamma _{h},\ldots,\gamma _{1}),\quad \sS = (\gamma
_{h},\ldots,\gamma _{1}, 0),
\end{equation}
and taking $\rho _{\rr }$ and $\rho _{\sS }$ to have block maxima
$M_{h+1},\ldots,M_{1}$.  We ascribe arbitrarily large artificial minima
denoted by $\infty $ to the empty first block of $\rho _{\rr }$ and
last block of $\rho _{\sS}$.  The block minima are then
\begin{equation}\label{e:Cauchy-input-minima}
\rho _{\rr }\colon (\infty ,v_{h},\ldots,v_{1}),\quad \rho _{\sS
}\colon (u_{h},\ldots,u_{1},\infty ).
\end{equation}
By \eqref{e:u-and-v}, these sequences are $(\sigma w_{0})^{-1}$-almost
decreasing and $(\sigma w_{0})^{-1}$-almost increasing, respectively,
so the hypotheses of Theorem~\ref{thm:Cauchy} are satisfied.  Note
that \eqref{e:Cauchy-input-minima} also implies that $\rho _{\rr } =
\rho' _{w_{0}(\gammabold )}$ and $\rho _{\sS } = \rho
_{w_{0}(\gammabold )}$.

We set the variables $\xx $ and $\yy $ in \eqref{e:Cauchy} to $\xx =
w_{0}(\overline{\zz})$, $\yy = w_{0}(\zz )$, that is, $x_{i} =
z_{l+1-i}^{-1}$, $y_{i} = z_{l+1-i}$, where $l = |\gammabold |$.
Grouping the variables $\xx $ into blocks $X_{i}$ of size $r_{i}$, the
variables $\yy $ into blocks $Y_{i}$ of size $s_{i}$, and the
variables $\zz $ into blocks $Z_{i}$ of size $\gamma _{i}$, we then
have $(X_{1}, \ldots, X_{h+1}) = (0, \overline{Z_{h}}, \ldots,
\overline{Z_{1}})$ and $(Y_{1}, \ldots, Y_{h+1}) = (Z_{h}, \ldots,
Z_{1}, 0)$.  Theorem~\ref{thm:Cauchy} now yields
\begin{multline}\label{e:Cauchy-output}
\frac{\prod _{\alpha \in R_{qt}}
(1 - q\, t\, \zz ^{\alpha })}{\prod _{\alpha \in R_{t}} (1- t\, \zz
^{\alpha })} = 
\frac{\prod _{i+1<j}\Omega [-q\, t\, Z_{i}\overline{Z_{j}}]}{\prod
_{i< j}\Omega [-t\, Z_{i}\overline{Z_{j}}]}\\
= \zz ^{((v_{1}-u_{1})^{\gamma _{1}}, \ldots, (v_{h}-u_{h})^{\gamma
_{h}})} \sum _{\lambda } t^{|\lambda |} w_{0} \Big( F^{(\sigma
w_{0})^{-1}}_{\sS ,(0;\lambda )+\rho _{w_{0}(\gammabold )}}(\zz ;q)
\overline{E^{(\sigma w_{0})^{-1}}_{\rr ,(\lambda;0) +\rho'
_{w_{0}(\gammabold )}}(\zz ;q)} \Big),
\end{multline}
where we used $\rho _{\rr } = \rho '_{w_{0}(\gammabold )}$, $\rho
_{\sS } = \rho _{w_{0}(\gammabold )}$, and
\begin{equation}\label{e:combine-rho}
\yy ^{-\rho _{\sS }} \xx ^{-\rho _{\rr }} = w_{0}(\zz ^{\rho'
_{w_{0}(\gammabold )} - \rho _{w_{0}(\gammabold )}}) =\zz
^{((v_{1}-u_{1})^{\gamma _{1}}, \ldots, (v_{h}-u_{h})^{\gamma _{h}})}.
\end{equation}
Because of the zero-length first block in $\rr $ and last block in
$\sS $, the sum in \eqref{e:Cauchy} is over tuples of partitions of
the form $(0,\lambda _{(h-1)},\ldots,\lambda _{(1)},0)$, where $\ell
(\lambda _{(i)})\leq \min (\gamma _{i},\gamma _{i+1})$.  These become
$(\lambda ;0)$ and $(0;\lambda )$ when interpreted as weights for $\GL
_{\rr }$ and $\GL _{\sS }$, respectively.

Dropping zero-length blocks as in
Remark~\ref{rem:sigma-hat-and-weak-r} (ii) gives
\begin{align}\label{e:drop-zeroes-E}
E^{(\sigma w_{0})^{-1}}_{\rr ,(\lambda ;0) +\rho' _{w_{0}(\gammabold
)}}(\zz ;q) & = E^{(\tau w_{0})^{-1}}_{w_{0}(\gammabold ),(\lambda ;0)
+\rho' _{w_{0}(\gammabold )}}(\zz ;q),
\\
\label{e:drop-zeroes-F}
F^{(\sigma w_{0})^{-1}}_{\sS ,(0; \lambda )+\rho _{w_{0}(\gammabold
)}}(\zz ;q) & = F^{(\theta w_{0})^{-1}}_{w_{0}(\gammabold ),(0;
\lambda ) +\rho _{w_{0}(\gammabold )}}(\zz ;q).
\end{align}

To complete the proof, we observe that
\begin{equation}\label{e:b-characterization}
b_{i} = p + c_{i+1}-c_{i} = \begin{cases}
\lfloor p \rfloor, & c_{i}>c_{i+1},\\
\lceil p \rceil , & c_{i}<c_{i+1}.
\end{cases}
\end{equation}
Hence $(b_{h},\ldots,b_{1}) = \eta + \lfloor p \rfloor (1^{h})$, where
$\eta $ is the descent indicator of the permutation $\sigma w_{0}\in
S_{h+1}$.  Since $\sigma w_{0}$
is a winding permutation with head $\theta w_{0}^{h}$ and tail $\tau
w_{0}^{h}$, \eqref{e:shift} and \eqref{e:winding-2} imply
\begin{equation}\label{e:F-change}
F^{(\theta w_{0})^{-1}}_{w_{0}(\gammabold ),(0; \lambda ) +\rho
_{w_{0}(\gammabold )}}(\zz ;q) = \zz ^{-(b_{h}^{\gamma
_{h}},\ldots,b_{1}^{\gamma _{1}})} F^{(\tau
w_{0})^{-1}}_{w_{0}(\gammabold ),(0; \lambda ) + (b_{h}^{\gamma
_{h}},\ldots,b_{1}^{\gamma _{1}}) + \rho _{w_{0}(\gammabold )}}(\zz
;q).
\end{equation}
Combining \eqref{e:Cauchy-output},
(\ref{e:drop-zeroes-E}--\ref{e:drop-zeroes-F}),
and \eqref{e:F-change}
gives \eqref{e:unstraightened}.
\end{proof}

\subsection{Proof of the main combinatorial theorem}
\label{ss:nest-identity-proof}

We now prove Theorem~\ref{thm:nest-identity}, using
Theorem~\ref{thm:stable-main}.
Before giving the full proof, we outline
the argument.  Theorem \ref{thm:stable-main} gives an expansion
\[
H = \sum _{\lambda } t^{|\lambda|}\Lcal ^{\tau w_{0}}_{w_{0}(\gammabold ),\, 
\beta /\alpha }(\zz ;q)
\]
for the Catalanimal $H$ in Theorem~\ref{thm:nest-identity}, in which $\beta$
and $\alpha $ depend on $\lambda $. Applying
Theorem~\ref{thm:L-vs-G}, we find that the surviving terms in $H_{\pol
}$ belong to those $\lambda $ for which we have $\beta \geq \alpha$
coordinate-wise, and these correspond one-to-one with nests in the
given den by Lemma \ref{lem:parametrization}.
Theorem~\ref{thm:L-vs-G} also gives us an expression
$t^{|\lambda|}\Lcal ^{\tau w_{0}}_{w_{0}(\gammabold ),\, \beta /\alpha
}(\zz ;q)_{\pol} = \sum t^{|\lambda|} q ^{h_{\tau w_{0}}(\betabold
/\alphabold )} \Gcal _{\tau w_{0} (\betabold/\alphabold )}(\zz ;
q^{-1})$ for these terms.  We then complete the proof by using the
results of \S \S \ref{ss:nest-stats}--\ref{ss:nest-LLT} to verify that
for the corresponding nest $\pi $, we have $a(\pi) = |\lambda|$,
$\dinv _{p}(\pi ) = h_{\tau w_{0}}(\betabold /\alphabold )$, and $\tau
w_0(\betabold /\alphabold) = \nubold (\pi)$.

\begin{proof}[Proof of Theorem~\ref{thm:nest-identity}]
Let $(h,p,\dd ,\ee )$ be the given den, with $\gG $ as in
\eqref{e:g-vector}.  We fix $s$ such that the line $y +p\, x = s$
passes above the heads and feet, as in \S \ref{ss:nest-LLT}, and
define $f_{i} = \lfloor s-p\, i \rfloor$ to be the $y$-coordinate of
the highest lattice point below the line $y +p\, x = s$ at $x = i$.
For $i = 1,\ldots,h$, we set
\begin{equation}\label{e:combinatorial-u-and-v}
u_{i} = f_{i} - e_{i},\quad v_{i} = f_{i-1} - d_{i-1}.
\end{equation}
We will apply Theorem~\ref{thm:stable-main} with these values of $s$,
$p$, $u_{i}$ and $v_{i}$.  The numbers in \eqref{e:b-and-c} are then
given by $b_{i} = f_{i-1}-f_{i}$ and $c_{i+1} = s-p\, i-f_{i}$.

We start by verifying that the hypotheses in \eqref{e:u-and-v} are
satisfied.  The hypothesis on the $u_{i}$ can be restated as
$u_{i}\geq u_{j} - \chi (c_{i+1}>c_{j+1})$ for $i<j$.  Since the
$c_{i}$ are distinct and $c_{i}\in [0,1)$ for all $i$, this is
equivalent to $u_{i}+c_{i+1}>u_{j}+c_{j+1}-1$, or to $u_{i}-p\, i -
f_{i} >u_{j}-p\, j-f_{j} - 1$.  The latter reduces to condition
\eqref{e:den-e-incr} in the definition of a den.  Similarly, the
hypothesis on the $v_{i}$ is equivalent to $v_{i}+c_{i}<v_{j}+c_{j}+1$
for $i<j$, which reduces to \eqref{e:den-d-decr}.

In \eqref{e:gamma-diff}, we have $u_{i} - v_{i+1} = d_{i}-e_{i}$, so
we can take $\gamma _{i} = g_{i}$.  We also have $u_{i}-v_{i}+b_{i} =
d_{i-1}-e_{i}$.  Hence, the Catalanimal in \eqref{e:stable-main}
coincides with the Catalanimal $H$ in \eqref{e:H-for-den}.

Turning to the right hand side of \eqref{e:stable-main}, we use
Theorem~\ref{thm:L-vs-G} to evaluate the polynomial part of $\Lcal
^{\tau w_{0}}_{w_{0}(\gammabold ),\, ((0;\lambda ) +(b_{h}^{\gamma
_{h}},\ldots,b_{1}^{\gamma _{1}})+\rho _{w_{0}(\gammabold )
})/((\lambda;0) +\rho' _{w_{0}(\gammabold ) })}(\zz ;q) = \Lcal ^{\tau
w_{0}}_{w_{0}(\gammabold ),\, \beta /\alpha }(\zz ;q)$, where 
\begin{align}\label{e:beta-and-alpha}
\beta & = (0;\lambda ) +(b_{h}^{\gamma _{h}},\ldots,b_{1}^{\gamma
_{1}})+\rho _{w_{0}(\gammabold ) }\\
\alpha & = (\lambda;0) +\rho' _{w_{0}(\gammabold ) }.
\end{align}
The weight $\beta $ is the concatenation of blocks
\begin{equation}\label{e:beta-block}
\lambda _{(k)} + ((f_{k-1}-e_{k})^{g_{k}})+ (g_{k}-1, g_{k}-2, \cdots
, 0)
\end{equation}
in the order $k = h,h-1,\ldots,1$, if we set $\lambda _{(h)} =
\varnothing $.  Similarly, $\alpha $ is the concatenation of blocks
\begin{equation}\label{e:alpha-block}
\lambda _{(k-1)} + ((f_{k-1}-d_{k-1})^{g_{k}}) + (g_{k}-1, g_{k}-2,
\cdots , 0)
\end{equation}
in the same order, with $\lambda _{(0)} = \varnothing $.  The $\lambda
_{(k)}$ for $1\leq k\leq h-1$ vary over partitions of length $\ell
(\lambda _{(k)})\leq \min (g_{k},g_{k+1})$.  By
Theorem~\ref{thm:L-vs-G}, $\Lcal ^{\tau w_{0}}_{w_{0}(\gammabold ),\,
\beta /\alpha }(\zz ;q)_{\pol } = 0$ unless $\alpha \leq \beta $
coordinate-wise.  From \eqref{e:beta-block} and \eqref{e:alpha-block},
we see that $\alpha \leq \beta $ if and only if $(\lambda
_{(k-1)})_{i} - d_{k-1}\leq (\lambda _{(k)})_{i} - e_{k}$ for $k =
1,\ldots, h$ and $i\leq g_{k}$.  Since $d_{k-1} = e_{k-1} + g_{k} -
g_{k-1}$, this is equivalent to the condition $e_{k}-g_{k} - (\lambda
_{(k)})_{i}\leq e_{k-1}-g_{k-1}-(\lambda _{(k-1)})_{i}$ in
Lemma~\ref{lem:parametrization}.  By that lemma, the indices $\lambda
$ for which $\Lcal ^{\tau w_{0}}_{w_{0}(\gammabold ),\, \beta /\alpha
}(\zz ;q)_{\pol } \not = 0$ correspond to nests $\pi $ in the given
den.

For these indices, Theorem~\ref{thm:L-vs-G} gives 
\begin{equation}\label{e:L-to-G}
\Lcal ^{\tau w_{0}}_{w_{0}(\gammabold ),\, \beta /\alpha }(\zz
;q)_{\pol } = q ^{h_{\tau w_{0}}(\betabold /\alphabold )} \Gcal _{\tau
w_{0} (\betabold/\alphabold )}(z_{1},\ldots,z_{l};q^{-1}),
\end{equation}
where $l=|\gG |$, and $\betabold/\alphabold $ is related to $\alpha $
and $\beta $ by the recipe in \eqref{e:nu(r,beta/alpha)} for blocks of
lengths $g_{h},\ldots,g_{1}$.  Writing $\betabold/\alphabold = (\beta
_{(h)}/\alpha _{(h)},\ldots,\beta _{(1)}/\alpha _{(1)})$ with
decreasing indices and using \eqref{e:beta-block} and
\eqref{e:alpha-block}, this recipe gives $\beta _{(k)} =
((f_{k-1}-e_{k} + g_{k})^{g_{k}}) + \lambda _{(k)}$ and $\alpha _{(k)}
= ((f_{k-1}-d_{k-1} + g_{k})^{g_{k}}) + \lambda _{(k-1)} =
((f_{k-1}-e_{k-1} + g_{k-1})^{g_{k}}) + \lambda _{(k-1)}$.  The
permutation $\sigma $ in Definition~\ref{def:nu(pi)} for the given den
and choice of $s$ is the same as $\tau $ in \eqref{e:L-to-G}. Hence, by
Remark~\ref{rem:nu(pi)-vs-lambda}, $\beta _{(k)}/\alpha _{(k)} =
\nubold (\pi )_{\tau (k)}$,
so $\tau w_{0}(\betabold /\alphabold
) = \tau (\beta _{(1)}/\alpha _{(1)},\ldots,\beta _{(h)}/\alpha
_{(h)}) = \nubold (\pi )$.

As noted in Definition~\ref{def:a(pi)}, we have $a(\pi ) = |\lambda |$
for the nest $\pi $ corresponding to $\lambda $.

We now show that $\dinv _{p}(\pi ) = h_{\tau w_{0}}(\betabold
/\alphabold )$.  Because the components $\beta _{(i)}/\alpha _{(i)}$
of $\betabold /\alphabold$ are indexed in decreasing order, a $(\tau
w_{0})$-triple $(a,b,c)$ in $\betabold /\alphabold$ has $a,c$ in
$\beta _{(j)}/\alpha _{(j)}$ and $b$ in $\beta _{(i)}/\alpha _{(i)}$
for $i>j$, with content $c(b)$ equal to $c(a)$ if $\tau (i) > \tau
(j)$,
or to $c(a)+1$ if $\tau (i) < \tau (j)$.  As in \S \ref{ss:nest-LLT},
the box $b\in \beta _{(i)}/\alpha _{(i)}$ corresponds to an element
$(S,k)\in \Sbold (\pi )$ with $S$ on the line $x = i-1$, with
$c(b)+c_{i}$ equal to the vertical distance between the line $y +p\, x
= s$ and the south endpoint of $S$.  Similarly, even though boxes $a$
and $c$ need not actually be in $\beta _{(j)}/\alpha _{(j)}$, the
boundary between them corresponds to a non-sink lattice point $P$ at
$x = j-1$ on some path $\pi _{k'}$ in $\pi $, with $c(a)+c_{j}$ equal
to the vertical distance between $y + p\, x = s$ and $P$.

The tuple $(P,k',S,k)$ is counted by $\dinv _{p}(\pi )$ if and only if
$0< c(b) - c(a) + c_{i} - c_{j} < 1$.  Since $c_{i},c_{j}\in [0,1)$,
we have $|c_{i} - c_{j}|<1$, and by the definition of $\tau $, we have
$\tau (i)<\tau (j)$ if and only if $c_{i}<c_{j}$.  If $\tau (i) > \tau
(j)$, it follows that $0< c(b) - c(a) + c_{i} - c_{j} < 1$ if and only
if $c(b) = c(a)$, while if $\tau (i) < \tau (j)$, it follows that $0<
c(b) - c(a) + c_{i} - c_{j} < 1$ if and only if $c(b) = c(a)+1$.
Hence, tuples $(P,k',S,k)$ counted by $\dinv _{p}(\pi )$ are in
bijective correspondence with $(\tau w_{0})$-triples in $\betabold
/\alphabold $, giving $\dinv _{p}(\pi ) = h_{\tau w_{0}}(\betabold
/\alphabold )$.

From \eqref{e:L-to-G} and the expressions for $a(\pi )$ and $\dinv
_{p}(\pi )$ we see that the polynomial part of the series on the right
hand side of \eqref{e:stable-main} is equal to
\begin{equation}\label{e:stable-main-pol}
\sum _{\pi } t^{a(\pi )} q^{\dinv _{p}(\pi )} \Gcal _{\nubold (\pi
)}(z_{1},\ldots,z_{l},q^{-1}),
\end{equation}
where the sum is over nests $\pi $ in the given den.  Since the
Catalanimal on the left hand side of \eqref{e:stable-main} is equal to
$H$, this proves Theorem~\ref{thm:nest-identity}.
\end{proof}

\bibliographystyle{hamsplain}
\bibliography{references}

\providecommand{\bysame}{\leavevmode\hbox to3em{\hrulefill}\thinspace}
\begin{thebibliography}{10}

\bibitem{BakDunFor00}
T.~H. Baker, C.~F. Dunkl, and P.~J. Forrester, \emph{Polynomial eigenfunctions
  of the {C}alogero-{S}utherland-{M}oser models with exchange terms},
  Calogero-{M}oser-{S}utherland models ({M}ontr\'{e}al, {QC}, 1997), CRM Ser.
  Math. Phys., Springer, New York, 2000, pp.~37--51.

\bibitem{BeGaHaTe99}
F.~Bergeron, A.~M. Garsia, M.~Haiman, and G.~Tesler, \emph{Identities and
  positivity conjectures for some remarkable operators in the theory of
  symmetric functions}, Methods Appl. Anal. \textbf{6} (1999), no.~3, 363--420,
  Dedicated to Richard A.\ Askey on the occasion of his 65th birthday, Part
  III.

\bibitem{BeGaSeXi16}
Francois Bergeron, Adriano Garsia, Emily Sergel~Leven, and Guoce Xin,
  \emph{Compositional {$(km,kn)$}-shuffle conjectures}, Int. Math. Res. Not.
  IMRN (2016), no.~14, 4229--4270.

\bibitem{BHMPS-llt}
Jonah Blasiak, Mark Haiman, Jennifer Morse, Anna Pun, and George~H. Seelinger,
  \emph{{LLT} polynomials in the {S}chiffmann algebra}, 2021,
  \mbox{arXiv:2112.07063 [math.CO]}.

\bibitem{BHMPS-delta}
\bysame, \emph{A proof of the {E}xtended {D}elta {C}onjecture}, 2021,
  \mbox{arXiv:2102.08815 [math.CO]}.

\bibitem{BHMPS-paths}
\bysame, \emph{A shuffle theorem for paths under any line}, 2021,
  \mbox{arXiv:2102.07931 [math.CO]}.

\bibitem{BurbSchi12}
Igor Burban and Olivier Schiffmann, \emph{On the {H}all algebra of an elliptic
  curve, {I}}, Duke Math. J. \textbf{161} (2012), no.~7, 1171--1231.

\bibitem{CarlMell18}
Erik Carlsson and Anton Mellit, \emph{A proof of the shuffle conjecture}, J.
  Amer. Math. Soc. \textbf{31} (2018), no.~3, 661--697.

\bibitem{FeHaHoShYa09}
B.~Feigin, K.~Hashizume, A.~Hoshino, J.~Shiraishi, and S.~Yanagida, \emph{A
  commutative algebra on degenerate {$\mathbb{CP}^1$} and {M}acdonald
  polynomials}, J. Math. Phys. \textbf{50} (2009), no.~9, 095215, 42.

\bibitem{FeigTsym11}
B.~L. Feigin and A.~I. Tsymbaliuk, \emph{Equivariant {$K$}-theory of {H}ilbert
  schemes via shuffle algebra}, Kyoto J. Math. \textbf{51} (2011), no.~4,
  831--854.

\bibitem{GarsHaim96}
A.~M. Garsia and M.~Haiman, \emph{Some natural bigraded {$S\sb n$}-modules and
  {$q,t$}-{K}ostka coefficients}, Electron. J. Combin. \textbf{3} (1996),
  no.~2, Research Paper 24, approx.\ 60 pp.\ (electronic), The Foata
  Festschrift.

\bibitem{GrojHaim07}
I.~Grojnowski and M.~Haiman, \emph{Affine {H}ecke algebras and positivity of
  {LLT} and {M}acdonald polynomials}, Unpublished manuscript, 2007.

\bibitem{HaHaLoReUl05}
J.~Haglund, M.~Haiman, N.~Loehr, J.~B. Remmel, and A.~Ulyanov, \emph{A
  combinatorial formula for the character of the diagonal coinvariants}, Duke
  Math. J. \textbf{126} (2005), no.~2, 195--232.

\bibitem{LaLeTh97}
Alain Lascoux, Bernard Leclerc, and Jean-Yves Thibon, \emph{Ribbon tableaux,
  {H}all-{L}ittlewood functions, quantum affine algebras, and unipotent
  varieties}, J. Math. Phys. \textbf{38} (1997), no.~2, 1041--1068.

\bibitem{LoehWarr08}
Nicholas~A. Loehr and Gregory~S. Warrington, \emph{Nested quantum {D}yck paths
  and {$\nabla(s_\lambda)$}}, Int. Math. Res. Not. IMRN (2008), no.~5, Art. ID
  rnm 157, 29.

\bibitem{Macdonald72}
I.~G. Macdonald, \emph{The {P}oincar\'{e} series of a {C}oxeter group}, Math.
  Ann. \textbf{199} (1972), 161--174.

\bibitem{Macdonald95}
\bysame, \emph{Symmetric functions and {H}all polynomials}, second ed., The
  Clarendon Press, Oxford University Press, New York, 1995, With contributions
  by A.~Zelevinsky, Oxford Science Publications.

\bibitem{Mellit16}
Anton Mellit, \emph{Toric braids and $(m,n)$-parking functions}, 2016,
  \mbox{arXiv:1604.07456 [math.CO]}.

\bibitem{Negut14}
Andrei Negut, \emph{The shuffle algebra revisited}, Int. Math. Res. Not. IMRN
  (2014), no.~22, 6242--6275.

\bibitem{SchiVass13}
Olivier Schiffmann and Eric Vasserot, \emph{The elliptic {H}all algebra and the
  {$K$}-theory of the {H}ilbert scheme of {$\mathbb A^2$}}, Duke Math. J.
  \textbf{162} (2013), no.~2, 279--366.

\bibitem{ShimWeym00}
Mark Shimozono and Jerzy Weyman, \emph{Graded characters of modules supported
  in the closure of a nilpotent conjugacy class}, European J. Combin.
  \textbf{21} (2000), no.~2, 257--288.

\bibitem{ShimZabr01}
Mark Shimozono and Mike Zabrocki, \emph{Hall-{L}ittlewood vertex operators and
  generalized {K}ostka polynomials}, Adv. Math. \textbf{158} (2001), no.~1,
  66--85.

\end{thebibliography}

\end{document}